\documentclass[a4paper,10pt,leqno]{amsart}
        \title[On crossed product rings with twisted involutions, \ldots]
        {On crossed product rings with twisted involutions,
        their module categories and $L$-theory}
        \author{Arthur Bartels} 
             \address{Imperial College London\\
               Huxley Building\\
               London SW7 2AZ, UK}
             \email{a.bartels@imperial.ac.uk}
             \urladdr{http://ma.ic.ac.uk/~abartels/}
        \author{Wolfgang L\"uck}
             \address{Westf\"alische Wilhelms-Universit\"at M\"unster\\
               Mathematisches Institut\\
               Einsteinstr.~62,
               D-48149 M\"unster, Germany}
             \email{lueck@math.uni-muenster.de}
             \urladdr{http://www.math.uni-muenster.de/u/lueck}
               \date{October, 2007}

\usepackage{hyperref}
\usepackage{calc}
\usepackage{enumerate,amssymb}
\usepackage[arrow,curve,matrix,tips,2cell]{xy}
  \SelectTips{eu}{10} \UseTips
  \UseAllTwocells

\DeclareMathAlphabet{\matheurm}{U}{eur}{m}{n}
\newcommand{\addcat}{\matheurm{Add\text{-}Cat}}
\newcommand{\addcatinv}{\matheurm{Add\text{-}Cat_{\text{inv}}}}
\newcommand{\addGcatinv}{\matheurm{Add\text{-}G\text{-}Cat_{\text{inv}}}}
\newcommand{\Gaddcat}{(G,v)\text{-}\matheurm{Add}\text{-}\matheurm{Cat}}
\newcommand{\sGaddcat}
{\matheurm{strict}\text{-}(G,v)\text{-}\matheurm{Add}\text{-}\matheurm{Cat}}
\newcommand{\Zcat}{\matheurm{\IZ\text{-}Cat}}
\newcommand{\Zcatinv}{\matheurm{\IZ\text{-}Cat_{\text{inv}}}}
\newcommand{\ZGcatinv}{\matheurm{\IZ\text{-}G\text{-}Cat_{\text{inv}}}}

\newcommand{\Groupoids}{\matheurm{Groupoids}}
\newcommand{\Or}{\matheurm{Or}}
\newcommand{\Spectra}{\matheurm{Spectra}}

\DeclareMathOperator{\asmb}{asmb}
\DeclareMathOperator{\aut}{aut}

\DeclareMathOperator{\cent}{cent}

\DeclareMathOperator{\colim}{colim}
\DeclareMathOperator{\func}{func}

\DeclareMathOperator{\forget}{forget}

\DeclareMathOperator{\id}{id}

\DeclareMathOperator{\ind}{ind}

\DeclareMathOperator{\map}{map}
\DeclareMathOperator{\mor}{mor}
\DeclareMathOperator{\pt}{pt}
\DeclareMathOperator{\pr}{pr}
\DeclareMathOperator{\res}{res}

\newcommand{\VCyc}{{\mathcal{VC}\text{yc}}}

  \newcommand{\IZ}{\mathbb{Z}}

  \newcommand{\cala}{\mathcal{A}}
  \newcommand{\calb}{\mathcal{B}}
  \newcommand{\calc}{\mathcal{C}}

  \newcommand{\calg}{\mathcal{G}}
  \newcommand{\calh}{\mathcal{H}}

  \newcommand{\cals}{\mathcal{S}}

  \newcommand{\bfE}{{\mathbf E}}

  \newcommand{\bfK}{{\mathbf K}}
  \newcommand{\bfL}{{\mathbf L}}

\theoremstyle{plain}
\newtheorem{theorem}{Theorem}[section]

\newtheorem{lemma}[theorem]{Lemma}
\newtheorem{corollary}[theorem]{Corollary}

\theoremstyle{definition}
\newtheorem{definition}[theorem]{Definition}

\newtheorem{example}[theorem]{Example}

\newtheorem{remark}[theorem]{Remark}

\theoremstyle{remark}

\makeatletter\let\c@equation=\c@theorem\makeatother

\newcommand{\version}[1]                       %marks the date of last editing and compilation
{\begin{center} last edited on #1\\
last compiled on \today
\end{center}
}

\newcommand{\EGF}[2]{E_{#2}(#1)}
\newcommand{\FGP}[1]{#1\text{-}\matheurm{FGP}}
\newcommand{\FGF}[1]{#1\text{-}\matheurm{FGF}}

\newcommand{\Groupoidsover}[1]{\Groupoids\downarrow #1}
\newcommand{\intgf}[2]{\int_{#1} #2}

\newcommand{\OrG}[1]{\matheurm{Or}(#1)}

\begin{document}

\maketitle

\begin{abstract}
We study the Farrell-Jones Conjecture with coefficients in an
additive $G$-category with involution. This is a variant of
the $L$-theoretic Farrell-Jones Conjecture which originally deals with
group rings with the standard involution.
We show that this formulation of the conjecture can be applied
to crossed product rings $R \ast G$ equipped with twisted involutions
and automatically implies the a priori more general fibered version.
\end{abstract}

\newlength{\origlabelwidth}
\setlength\origlabelwidth\labelwidth

%%%%%%%%%%%%%%%%%%%%%%%%%%%%%%%%%%%%%%%%%%%%%%%%%%%%%%%%%%%%%%%%%%%%%%%%%%%%%%%%%%%%%%%%%
%%%%%%%%%%%%%%%%%%%%%%%%%%%%%%%%%%%%%%%%%%%%%%% Introduction %%%%%%%%%%%%%%%%%%%%%%%%%%%%
%%%%%%%%%%%%%%%%%%%%%%%%%%%%%%%%%%%%%%%%%%%%%%%%%%%%%%%%%%%%%%%%%%%%%%%%%%%%%%%%%%%%%%%%%

\typeout{----------------------------  Introduction ------------------------------------}
\section*{Introduction}
\label{sec:introduction}

The \emph{Farrell-Jones Conjecture for algebraic $L$-theory} predicts for a
group $G$ and a ring $R$ with involution $r \mapsto \overline{r}$ that
the so called \emph{assembly map}
\begin{eqnarray}
&
\asmb^{G,R}_n \colon H_n^G\bigl(\EGF{G}{\VCyc};\bfL_R^{\langle - \infty\rangle}\bigr)
\to L_n^{\langle - \infty\rangle}(RG)
&
\label{ass_for_RG}
\end{eqnarray}
is bijective for all $n \in \IZ$. Here the target is the \emph{$L$-theory}
of the group ring $RG$ with the standard involution sending
$\sum_{g \in G} r_g \cdot g$ to $\sum_{g \in G}  \overline{r_g} \cdot
g^{-1}$.  This is the group one wants to understand. It is a crucial
ingredient in the surgery program for the classification of closed
manifolds. The source is a much  easier
to handle term, namely, a $G$-homology theory  applied
to the \emph{the classifying space $\EGF{G}{\VCyc}$
of the family $\VCyc$ of virtually cyclic subgroups of $G$}.
There is also a $K$-theory version of the Farrell-Jones Conjecture.
The original source for the (Fibered) Farrell-Jones Conjecture is the paper by
Farrell-Jones~\cite[1.6 on page~257 and~1.7 on
page~262]{Farrell-Jones(1993a)}. More information can be found
for instance in the survey article~\cite{Lueck-Reich(2005)}.

In this paper we study the Farrell-Jones Conjecture with
coefficients in an additive $G$-category with involution.
We show that this more general formulation of the conjecture allows to
consider instead of the group ring $RG$ the crossed product ring with involution
$R\ast _{c,\tau,w} G$
(see Section~\ref{sec:Crossed_product_rings_and_involutions}),
which is a generalization of the twisted group ring, and
to use more general involutions, for instance the one given by twisting
the standard involution with a group homomorphism $w_1 \colon G \to
\{\pm 1\}$. The data describing $R\ast _{c,\tau} G$ and more general
involutions are pretty complicated. It turns out that it is
convenient to put these into a more general but easier to handle
context, where the coefficients are given by an additive
$G$-categories $\cala$ with involution
(see Definition~\ref{def:additive_G-category_with_involution}).

\begin{definition}[$L$-theoretic Farrell-Jones Conjecture]
\label{def:L-theoretic_Farrell-Jones_Conjecture}
  A group $G$ together with an additive $G$-category with involution
  $\cala$ satisfy the \emph{$L$-theoretic Farrell-Jones Conjecture
    with coefficients in $\cala$} if the assembly map
\begin{eqnarray}
& \asmb^{G,\cala}_n \colon H_n^G\bigl(\EGF{G}{\VCyc};\bfL_{\cala}^{\langle - \infty\rangle}\bigr) \to
 H_n^G\bigl(\pt;\bfL_{\cala}^{\langle - \infty\rangle}\bigr)
= L_n^{\langle - \infty\rangle}\left(\intgf{G}{\cala}\right).
&
\label{ass_for_cala}
\end{eqnarray}
induced by the projection $\EGF{G}{\VCyc} \to \pt$ is bijective for all $n \in \IZ$..

A group $G$ satisfies the \emph{$L$-theoretic Farrell-Jones
  Conjecture} if for any additive $G$-category with involution $\cala$
the \emph{$L$-theoretic Farrell-Jones Conjecture with coefficients in
  $\cala$} is true.
\end{definition}
Here $\intgf{G}{\cala}$ is a certain homotopy colimit which yields an
additive category with involution and we use the $L$-theory associated
to an additive category with involution due to Ranicki
(see~\cite{Ranicki(1988)}, \cite{Ranicki(1992)}
and~\cite{Ranicki(1992a)}). The $G$-homology theory
$H_n^G\bigl(-;\bfL_{\cala}^{\langle - \infty\rangle}\bigr)$ is briefly
recalled in Section~\ref{sec:G-homology_theories}. If $R$ is a ring
with involution, $\cala$ is the additive category with involution
given by finitely generated free $R$-modules and we equip $\cala$ with
the trivial $G$-action, then the assembly map~\eqref{ass_for_cala}
agrees with the one for $RG$ in~\eqref{ass_for_RG} (see
Theorem~\ref{the:FJC_for_crossed_products} below).  This general setup
is also a very useful framework when one is dealing with categories
appearing in controlled topology, which is an important tool for
proving the Farrell-Jones Conjecture for certain groups.

Next we state the main results of this paper.

\begin{theorem} \label{the:FJC_for_crossed_products}
Suppose that $G$ satisfies the $L$-theoretic Farrell-Jones
Conjecture in the sense of
Definition~\ref{def:L-theoretic_Farrell-Jones_Conjecture}. Let $R$
be ring  with the data $(c,\tau,w)$ and $R \ast_{c,\tau,w} G$ be
the associated crossed product ring with involution as explained
in Section~\ref{sec:Crossed_product_rings_and_involutions}. Then
the assembly map
\begin{eqnarray}
&
\asmb^{G,R_{c,\tau,w}}_n \colon
H_n^G\bigl(\EGF{G}{\VCyc};\bfL_{R,c,\tau,w}^{\langle - \infty\rangle}\bigr)
\to L_n^{\langle - \infty\rangle}(R \ast_{c,\tau,w} G)
&
\end{eqnarray}
is bijective.

Here $\bfL_{R,c,\tau,w}^{\langle - \infty\rangle}$ is a functor from
the orbit category $\Or(G)$ to the category of spectra such that
$\pi_n\bigl(\bfL_{R,c,\tau,w}^{\langle - \infty\rangle}(G/H)\bigr)$
for $H \subseteq G$ agrees with $L_n^{\langle -\infty \rangle}(R
\ast_{c|_H,\tau|_N,w|_H} H)$.
\end{theorem}

Another important feature is that in this setting
the (unfibered) Farrell-Jones Conjecture does already imply the
fibered version.

\begin{definition}[Fibered $L$-theoretic Farrell-Jones Conjecture]
\label{def:fibered_L-theoretic_Farrell-Jones_Conjecture}
A group $G$ satisfies the \emph{fibered $L$-theoretic Farrell-Jones
  Conjecture} if for any group homomorphism $\phi \colon K \to G$ and
additive $G$-category with involution $\cala$ the assembly map
\begin{eqnarray*}
& \asmb^{\phi,\cala}_n \colon
H_*^K\bigl(\EGF{G}{\phi^*\VCyc};\bfL_{\phi^*\cala}^{\langle - \infty\rangle}\bigr) \to
L_n^{\langle - \infty\rangle}\left(\intgf{K}{\phi^*\cala}\right).
&
\end{eqnarray*}
is bijective for all $n \in \IZ$, where the family $\phi^*\VCyc$
of subgroups of $K$ consists of subgroups $L \subseteq K$ with
$\phi(L)$ virtually cyclic and $\phi^* \cala$ is the additive
$K$-category with involution obtained from $\cala$ by
restriction with $\phi$.
\end{definition}

Obviously the fibered version for the group $G$ of
Definition~\ref{def:fibered_L-theoretic_Farrell-Jones_Conjecture}
implies the version for the group $G$ of
Definition~\ref{def:L-theoretic_Farrell-Jones_Conjecture}, take $\phi
= \id$ in
Definition~\ref{def:fibered_L-theoretic_Farrell-Jones_Conjecture}.
The converse is also true.

\begin{theorem} \label{the:fibered_versus_unfibered}
  Let $G$ be a group. Then $G$ satisfies the fibered $L$-theoretic
  Farrell-Jones Conjecture  if
  and only if $G$ satisfies the $L$-theoretic Farrell-Jones Conjecture
  of Definition~\ref{def:L-theoretic_Farrell-Jones_Conjecture}.
\end{theorem}

A general statement of a \emph{Fibered Isomorphism Conjecture} and
the discussion of its inheritance properties under subgroups and
colimits of groups can be found
in~\cite[Section~4]{Bartels-Echterhoff-Lueck(2007colim)} (see
also~\cite[Appendix]{Farrell-Jones(1993a)},
\cite[Theorem~7.1]{Farrell-Linnell(2003a)}).  These very useful
inheritance properties do not hold for the unfibered version of
Definition~\ref{def:L-theoretic_Farrell-Jones_Conjecture}. The
next three corollaries are immediate consequences of
Theorem~\ref{the:fibered_versus_unfibered} and~\cite[Theorem~3.3,
Lemma~4.4, Lemma~4.5 and
Lemma~4.6]{Bartels-Echterhoff-Lueck(2007colim)}.

\begin{corollary}\label{cor:directed_colimits}
Let $\{G_i \mid i \in I\}$ be a directed system  (with not
necessarily injective) structure maps and let $G$ be its colimit
$\colim_{i \in I} G_i$. Suppose that $G_i$ satisfy the
Farrell-Jones Conjecture of
Definition~\ref{def:L-theoretic_Farrell-Jones_Conjecture} for
every $i \in I$.

Then $G$ satisfies the Farrell-Jones Conjecture
of Definition~\ref{def:L-theoretic_Farrell-Jones_Conjecture}.
\end{corollary}

\begin{corollary} \label{cor:extensions}
Let $1 \to K \to G \xrightarrow{p}  Q \to 1$ be an extension of groups.
Suppose that the group $Q$ and for any virtually cyclic subgroup $V \subseteq Q$ the
group $p^{-1}(V)$ satisfy the Farrell-Jones Conjecture
of Definition~\ref{def:L-theoretic_Farrell-Jones_Conjecture}.

Then the group $G$ satisfies the Farrell-Jones Conjecture
of Definition~\ref{def:L-theoretic_Farrell-Jones_Conjecture}.
\end{corollary}

\begin{corollary}\label{cor:subgroups}
If $G$ satisfies the Farrell-Jones Conjecture
of Definition~\ref{def:L-theoretic_Farrell-Jones_Conjecture},
then any subgroup $H \subseteq G$ satisfies the Farrell-Jones Conjecture
of Definition~\ref{def:L-theoretic_Farrell-Jones_Conjecture}.
\end{corollary}

Corollary~\ref{cor:extensions} and Corollary~\ref{cor:subgroups}
have also been proved in~\cite{Hambleton-Pedersen-Rosenthal(2007)}.

\begin{remark} \label{rem:products}
Suppose that the Farrell-Jones Conjecture
of Definition~\ref{def:L-theoretic_Farrell-Jones_Conjecture} has been proved for
the product of two virtually cyclic subgroups.

Then Corollary~\ref{cor:extensions} and Corollary~\ref{cor:subgroups} imply that
$G \times H$ satisfy the
Farrell-Jones Conjecture of Definition~\ref{def:L-theoretic_Farrell-Jones_Conjecture}
if and only if both $G$ and $H$ satisfy the
Farrell-Jones Conjecture of Definition~\ref{def:L-theoretic_Farrell-Jones_Conjecture}
\end{remark}

It is sometimes useful to have strict structures on $\cala$, e.g., the
involution is desired to be strict and there should be a (strictly
associative) functorial direct sum.  The functorial direct sum is
actually needed in some proofs in order to guarantee good
functoriality properties of certain categories arising from controlled
topology. We will show

\begin{theorem}\label{the:strict_coefficients}
  The group $G$ satisfies the $L$-theoretic Farrell-Jones Conjecture
  of Definition~\ref{def:L-theoretic_Farrell-Jones_Conjecture} if it
  satisfies the obvious version of it, where one only considers
  additive $G$-category with (strictly associative) functorial direct
  sum and strict involution (see
  Definition~\ref{def:additive_G-category_with_oplus_and_inv}).
\end{theorem}

The Farrell-Jones Conjecture with coefficients (in $K$- and $L$-theory)
has been introduced in~\cite{Bartels-Reich(2005)}.
Our treatment here is more general in that we allow involutions that are not
necessarily strict and also deal with twisted involutions on the crossed product ring.

All results mentioned here have obvious analogues for $K$-theory
whose proof is actually easier since one does not have to deal
with the involutions.

The work was financially supported by the
Sonderforschungsbereich 478 \--- Geometrische Strukturen in der
Mathematik \--- and the Max-Planck-Forschungspreis of the second
author.

The paper is organized as follows:\\[1mm]
\begin{tabular}{ll}%
\ref{sec:Additive_categories_with_involution}.
& Additive categories with involution
\\%
\ref{sec:Additive_categories_with_weak_(G,v)-action}
& Additive categories with weak $(G,v)$-action
\\%
\ref{sec:Making_an_additive_categories_with_weak_(G,v)-action_strict}
& Making an additive categories with weak $(G,v)$-action strict
\\%
\ref{sec:Crossed_product_rings_and_involutions}.
& Crossed product rings and involutions
\\%
\ref{sec:Connected_groupoids_and_additive_categories}.
& Connected groupoids and additive categories
\\%
\ref{sec:From_crossed_product_rings_to_additive_categories}.
& From crossed product rings to additive categories
\\%
\ref{sec:Connected_groupoids_and_additive_categories_with_involutions}.
& Connected groupoids and additive categories with involutions
\\%
\ref{sec:From_crossed_product_rings_with_involution_to_additive_categories_with_involution}.
& From crossed product rings with involution to additive categories with involution
\\%
\ref{sec:G-homology_theories}.
& $G$-homology theories
\\%
\ref{sec:Z-categories_and_additive_categories_with_involutionsG-homology_theories_and_restriction}.
& $\IZ$-categories and additive categories with involutions
\\%
\ref{sec:G-homology_theories_and_restriction}.
& $G$-homology theories and restriction
\\%
\ref{sec:Proof_of_the_main_theorems}.
& Proof of the main theorems
\\
& References
\end{tabular}

%%%%%%%%%%%%%%%%%%%%%%%%%%%%%%%%%%%%%%%%%%%%%%%%%%%%%%%%%%%%%%%%%%%%%%%%%%%%%%%%%%%%%%%%%
%%%%%%%%%%%%%%%%%%%% Additive categories with involution %%%%%%%%%%%%%%%%%%%%%%%%%%%%%%%%
%%%%%%%%%%%%%%%%%%%%%%%%%%%%%%%%%%%%%%%%%%%%%%%%%%%%%%%%%%%%%%%%%%%%%%%%%%%%%%%%%%%%%%%%%

\typeout{--------------  Section 1: Additive categories with involution -----------------}

\section{Additive categories with involution}
\label{sec:Additive_categories_with_involution}

In this section we will review the notion of an additive category with involution
as it appears and is used in the literature. This will be one of our main examples.

Let $\cala$ be an \emph{additive category}, i.e., a small category $\cala$ such that for
two objects $A$ and $B$ the morphism set $\mor_{\cala}(A,B)$ has the structure of an
abelian group and the direct sum $A \oplus B$ of two objects $A$ and $B$
exists and the obvious compatibility conditions hold.
A covariant \emph{functor of additive categories} $F \colon \cala_0 \to\cala_1$
is a covariant functor such that for two objects $A$ and $B$ in $\cala_0$ the map
$\mor_{\cala_0}(A,B) \to \mor_{\cala_1}(F(A), F(B))$ sending $f$ to $F(f)$ respects the
abelian group structures and $F(A \oplus B)$ is  a model for $F(A) \oplus F(B)$.
The notion of a contravariant functor of additive categories is defined
analogously.

An \emph{involution $(I,E)$ on an additive category} $\cala$ is contravariant functor
\begin{eqnarray}
&I \colon \cala  \to \cala &
\label{involution_I_on_an_additive_category}
\end{eqnarray}
of additive categories together with a natural equivalence of
such functors
\begin{eqnarray}
&E \colon \id_{\cala}  \to  I^2 := I \circ I&
\label{E_belonging_to_I}
\end{eqnarray}
such that for every object $A$ we have the
equality of morphisms
\begin{eqnarray}
E(I(A)) & = & I(E(A)^{-1}) \colon I(A) \to I^3(A).
\label{E(I(A))_is_I(E(A)-1)}
\end{eqnarray}
In the sequel we often write $I(A) =
A^*$ and $I(f) = f^*$ for a morphism $f \colon A \to B$ in $\cala$.
If $I^2 = \id_\cala$ and $E(A) = \id_A$ for all objects $A$,
then we call $I = (I,\id)$ a \emph{strict involution}.

\begin{definition}[Additive category with involution]
\label{def:additive_category_with_involution}
An \emph{additive category with involution} is an additive category
together with an involution $(I,E)$.

An \emph{additive category with strict involution} is an additive category
together with a strict involution $I$.
\end{definition}

The following example is a key example and illustrates why
one cannot expect in concrete situation that the involution is strict.

\begin{example} \label{exa:R-mod_as_add_cat}
Let $R$ be a ring. Let $\FGP{R}$ be the category of finitely
generated projective $R$-modules. This becomes an additive
category by the direct sum of $R$-modules and the elementwise
addition of $R$-homomorphisms.

A \emph{ring with involution} is a ring $R$ together with a map
$R \to R, \;r \mapsto \overline{r}$ satisfying
$\overline{1} = 1$, $\overline{r+s} = \overline{r} + \overline{s}$ and
$\overline{r\cdot s} = \overline{s} \cdot \overline{r}$ for $r,s \in R$.
Given a ring with involution $R$, define an involution $I$ on the additive category
$\FGP{R}$ as follows. Given a finitely generated projective $R$-module $P$,
let $I(P) = P^*$ be the finitely generated projective $\hom_R(P,R)$,where
for $r \in R$ and $f \in  \hom_R(P,R)$ the element $rf \in \hom_R(P,R)$ is defined by
$rf(x) = f(x) \cdot \overline{r}$ for $x \in P$. The desired natural transformation
$$E \colon \id_{\FGP{R}} \to I^2$$ assigns to a finitely generated projective $R$-module $P$
the $R$-isomorphism $P \xrightarrow{\cong} (P^*)^*$ sending $x \in P$ to
$\hom_R(P,R) \to R, \; f \mapsto \overline{f(x)}$.
\end{example}

A \emph{functor of additive categories with involution} $(F,T) \colon
\cala \to \calb$ consists of a covariant functor $F$ of the underlying additive
categories together with a natural  equivalence $T\colon
F \circ I_{\cala}\to  I_{\calb} \circ F$ such that for every
object $A$ in $\cala$ the following diagram commutes
\begin{eqnarray}
&
\xymatrix@!C=8em{F(A) \ar[d]^-{E_{\calb}(F(A))} \ar[r]^-{F(E_{\cala}(A))}
& F\left(A^{**}\right)  \ar[d]^-{T(A^*)}
\\
F(A)^{**} \ar[r]^-{T(A)^*}
& F(A^*)^*
}
&
\label{F(A)_F(A)aastast_F(Aastast)_F(Aast)ast}
\end{eqnarray}
If $T(A) = \id_A$ for all objects $A$, then we call $F$ a
\emph{strict} functor of additive categories with involution.

The \emph{composite} of functors of additive categories with
involution  $(F_1,T_1) \colon
\cala_1 \to \cala_2$ and $(F_2,T_2) \colon \cala_2 \to \cala_3$ is
defined to be $(F_2 \circ F_1, T_2 \circ T_1)$, where $F_2 \circ F_1$
is the composite of functors of additive categories
and the natural equivalence $T_2 \circ T_1$  assigns to an object $A \in
\cala_1$ the isomorphism in $\cala_3$
$$
F_2 \circ F_1  \circ I_{\cala_1}(A) \xrightarrow{F_2(T_1(A))} F_2 \circ I_{\cala_2} \circ F_1(A)
\xrightarrow{T_2(F_1(A))} I_{\cala_3} \circ F_2 \circ F_1(A).$$
A \emph{natural transformation $S \colon (F_1,T_1) \to (F_2,T_2)$
of functors $\cala_1 \to \cala_2$
of additive categories with involutions}
is a natural transformation $S \colon F_1 \to F_2$ of functors
of additive categories such that for every object $A$ in $\cala$
the following diagram  commutes
\begin{eqnarray}
&
\xymatrix@!C=8em{F_1(I_{\cala_1}(A))  \ar[r]^-{T_1(A)} \ar[d]^-{S(I_{\cala_1}(A))}
& I_{\cala_2}(F_1(A))
\\
F_2(I_{\cala_1}(A))  \ar[r]^-{T_2(A)}
& I_{\cala_2}F_2(A)) \ar[u]_-{I_{\cala_2}(S(A))}
}
&
\label{F_I_T}
\end{eqnarray}

%%%%%%%%%%%%%%%%%%%%%%%%%%%%%%%%%%%%%%%%%%%%%%%%%%%%%%%%%%%%%%%%%%%%%%%%%%%%%%%%%%%%%%%%%
%%%%%%%%%%%%%%%%%%%% Additive categories with weak (G,v)-action %%%%%%%%%%%%%%%%%%%%%%%%%%%%%%%%
%%%%%%%%%%%%%%%%%%%%%%%%%%%%%%%%%%%%%%%%%%%%%%%%%%%%%%%%%%%%%%%%%%%%%%%%%%%%%%%%%%%%%%%%%

\typeout{--------------  Section 2: Additive categories with weak $(G,v)$-action -----------------}

\section{Additive categories with weak $(G,v)$-action}
\label{sec:Additive_categories_with_weak_(G,v)-action}

In the sequel $G$ is a group and $v \colon G \to \{\pm 1\}$ is a group
homomorphism to the multiplicative group $\{\pm 1\}$.  In this section
we want to introduce the notion of an additive category with weak
$(G,v)$-action such that the notion of an additive category with
involution is the special case of an additive category with weak
$(\IZ/2,v)$-action for $v \colon \IZ/2 \to \{\pm 1\}$ the unique group
isomorphism and we can also treat $G$-actions up to natural
equivalence.  Notice that this will force us to deal with covariant
and contravariant functors simultaneously. The homomorphism $v$ will
take care of that.

We call a functor $+1$-variant if it is covariant and $-1$-variant if
it is contravariant.  If $F_1 \colon \calc_0 \to \calc_1$ is an
$\epsilon_1$-variant functor and $F_2 \colon \calc_1 \to \calc_2$ is
an $\epsilon_2$-variant functor, then the composite $F_2 \circ F_1
\colon \calc_0 \to \calc_2$ is $\epsilon_1\epsilon_2$-variant functor.
If $f \colon x_0 \to x_1$ is an isomorphism and $\epsilon \in \{\pm
1\}$, then define $f^{\epsilon} \colon x_0 \to x_1 $ to be $f$ if
$\epsilon = 1$ and $f^{\epsilon} \colon x_1 \to x_0 $ to be the
inverse of $f$ if $\epsilon = -1$.  If $F \colon \calc_0 \to \calc_1$
is $\epsilon$-variant and $f \colon x_0 \xrightarrow{\cong} y_0$ is an
isomorphism in $\calc_0$, then $F(f)^{\epsilon} \colon F(x_0) \to
F(x_1)$ is an isomorphism in $\calc_1$.

\begin{definition}[Additive category with weak $(G,v)$-action]
  \label{def:additive_category_with_weak_(G,v)-action}
  Let $G$ be a group together with a group homomorphism
  $v \colon G \to \{\pm 1\}$. An \emph{additive category with weak $(G,v)$-action} $\cala$ is an
  additive category together with the following data:

  \begin{itemize}

  \item For every $g \in G$ we have a $v(g)$-variant functor $R_g \colon
    \cala \to \cala$ of additive categories;

  \item For every two elements $g,h \in G$ there is a natural
    equivalence of $v(gh)$-variant functors of additive categories
    $$L_{g,h} \colon R_{gh} \to R_h \circ R_g.$$

\end{itemize}

We require:

\begin{enumerate}

\item $R_e = \id$ for $e \in G$ the unit element;

\item $L_{g,e} = L_{e,g} = \id$ for all $g \in G$;

\item \label{def:additive_category_with_weak_(G,v)-action:condition_for_(g,h,k)} The
  following diagram commutes for all $g,h,k \in G$ and objects $A$ in
  $\cala$
$$\xymatrix@!C=10em{R_{ghk}(A) \ar[r]^-{L_{gh,k}(A)} \ar[d]^-{L_{g,hk}(A)}
  & R_k(R_{gh}(A)) \ar[d]^-{R_k(L_{g,h}(A))^{v(k)}}
  \\
  R_{hk}(R_g(A)) \ar[r]^-{L_{h,k}(R_g(A))} & R_k(R_h(R_g(A)))}
$$
\end{enumerate}

If for every two elements $g,h \in G$ we have $L_{g,h} = \id$ and in particular
$R_{gh} = R_{h}R_{g}$, we call $\cala$ with these data an
\emph{additive category with strict $(G,v)$-action} or briefly a
\emph{additive $(G,v)$-category}.  If $v$ is trivial, we will omit it
from the notation.
\end{definition}

Let $\cala$ and $\calb$ be two additive categories with weak
$(G,v)$-action and let $\epsilon \in \{\pm 1\}$. An
\emph{$\epsilon$-variant functor $(F,T) \colon \cala \to \calb$ of
  additive categories with weak $(G,v)$-action} is a
$\epsilon$-variant functor $F \colon \cala \to \calb$ of additive
categories together with a collection $\{T_g \mid g \in G\}$ of
natural transformations of $\epsilon v(g)$-variant functors of
additive categories $T_g \colon F \circ R_g^{\cala} \to R_g^{\calb}
\circ F$. We require that for all $g,h \in G$ and all objects $A$ in
$\cala$ the following diagram commutes
\begin{eqnarray}
&&
\xymatrix@!C=10em{
F(R_{hg}^{\cala}(A))\ar[r]^-{F(L_{h,g}^{\cala}(A))^{\epsilon}} \ar[d]^-{T_{hg}(A)}
&
F(R_g^{\cala}(R_h^{\cala}(A))) \ar[r]^-{T_g(R_h(A))}
& R_g^{\calb}(F(R_h^{\cala}(A))) \ar[d]^-{R_g^{\calb}(T_{h}(A))^{v(g)}}
\\
R_{hg}^{\calb}(F(A)) \ar[rr]^-{L_{hg}^{\calb}(F(A))}
& &
R_g^{\calb}(R_h^{\calb}(F(A)))
}
\label{F_T_L_compatible}
\end{eqnarray}

The composite $(F_2,T_2) \circ (F_1,T_1) \colon \cala_1 \to \cala_3$
of an $\epsilon_1$-variant functor of additive categories with weak
$(G,v)$-action $(F_1,T_1) \colon \cala_1 \to \cala_2$ and an
$\epsilon_2$-variant functor of additive categories with weak
$(G,v)$-action $(F_2,T_2) \colon \cala_2 \to \cala_3$ is the
$\epsilon_1\epsilon_2$-variant functor of additive categories with
weak $(G,v)$-action whose underlying $\epsilon_1\epsilon_2$-variant
functor of additive categories is $F_2 \circ F_1 \colon \cala_1 \to
\cala_3$ and whose required natural transformations for $g \in G$ are
given for an object $A$ in $\cala_1$ by
$$F_2 \circ F_1 \circ R^{\cala_1}_g(A)
\xrightarrow{F_2((T_2)_g(A))^{\epsilon_2}} F_2 \circ R^{\cala_2}_g
\circ F_1(A) \xrightarrow{(T_2)_g(F_1(A))} R^{\cala_3}_g \circ F_2
\circ F_1(A).$$
A \emph{natural transformation $S \colon (F_1,T_1) \to
  (F_2,T_2)$ of functors $\cala_1 \to \cala_2$ of additive categories
  with weak $(G,v)$-action} is a natural transformation $S \colon F_1
\to F_2$ of functors of additive categories such that for all $g \in
G$ and objects $A$ in $\cala_1$ the following diagram commutes

\begin{eqnarray}
&
\xymatrix@!C=8em{F_1(R_g^{\cala_1}(A))  \ar[r]^-{(T_1)_g(A)} \ar[d]^-{S(R_g^{\cala_1}(A))}
& R_g^{\cala_2}(F_1(A)) \ar[d]^-{\left(R_g^{\cala_2}(S(A))\right)^{v(g)}}
\\
F_2(R_g^{\cala_1}(A))  \ar[r]^-{(T_2)_g(A)}
& R_g^{\cala_2}(F_2(A))
}
&
\label{F_I_T_add_G-cat}
\end{eqnarray}

An \emph{$\epsilon$-variant functor $F \colon \cala \to \calb$ of
additive categories with strict $(G,v)$-action} is an
$\epsilon$-variant functor $F \colon \cala \to \calb$ of additive
categories satisfying $F \circ R^{\cala_1}_g = R^{\cala_2}_g \circ F$
for all $g \in G$.  A \emph{natural transformation $S \colon F_1 \to F_2$
of $\epsilon$-variant functors $\cala_1 \to \cala_2$ of
additive categories with strict $(G,v)$-action} is a natural
transformation $S \colon F_1 \to F_2$ of additive categories
satisfying $S(R^{\cala_1}_g(A)) = R^{\cala_2}_g(S(A))^{v(g)}$ for all
$g \in G$ and objects $A$ in $\cals_1$.

\begin{example}[Additive categories with involution]
  \label{exa:additive_categories_with_involution}
  Given an additive category $A$, the structure of an additive
  category with weak $(\IZ/2,v)$-action for $v \colon \IZ/2 \to \{\pm 1\}$
  the unique group isomorphism is the same as an additive
  category with involution. Namely, let $t \in \IZ/2$ be the
  generator. Given an involution $(I,E)$ in the sense of
  Definition~\ref{def:additive_category_with_involution}, define the
  structure of an additive category with weak $(\IZ/2,v)$-action in
  the sense of Definition~\ref{def:additive_category_with_weak_(G,v)-action} by putting
  $R_e = \id$, $R_t = I$, $L_{e,e} = L_{t,e} = L_{t,e} = \id$ and
  $L_{t,t} = E$.
  Condition~\ref{def:additive_category_with_weak_(G,v)-action:condition_for_(g,h,k)} in
  Definition~\ref{def:additive_category_with_weak_(G,v)-action} follows
  from condition~\eqref{E(I(A))_is_I(E(A)-1)}. Given the
  structure of an additive category with weak $(\IZ/2,v)$-action,
  define the involution $(E,I)$ by $E = R_t$ and $I = L_{t,t}$.  The
  corresponding statement is true for functors of additive categories
  with weak $(\IZ/2,v)$-action and natural transformations between
  them, where diagram~\eqref{F(A)_F(A)aastast_F(Aastast)_F(Aast)ast}
  corresponds to diagram~\eqref{F_T_L_compatible}.

  Analogously we get that the structure of a additive category with
  strict $(\IZ/2,v)$-action is the same as an additive category with
  strict involution.
\end{example}

%%%%%%%%%%%%%%%%%%%%%%%%%%%%%%%%%%%%%%%%%%%%%%%%%%%%%%%%%%%%%%%%%%%%%%%%%%%%%%%%%%%%%%%%%
%%%%%%%%%% Making an additive categories with $(G,w)$-action strict %%%%%%%%%%%%%%%%%%%%%
%%%%%%%%%%%%%%%%%%%%%%%%%%%%%%%%%%%%%%%%%%%%%%%%%%%%%%%%%%%%%%%%%%%%%%%%%%%%%%%%%%%%%%%%%

\typeout{--------------  Section 3: Making an additive categories with $(G,w)$-action strict -----------------}

\section{Making an additive categories with weak $(G,v)$-action strict}
\label{sec:Making_an_additive_categories_with_weak_(G,v)-action_strict}

Many interesting examples occur as additive categories
with weak $(G,v)$-action which are
not necessarily strict.  On the other hand additive categories
with strict $(G,v)$-action are easier to handle. We explain how we can turn an additive
category with weak $(G,v)$-action $\cala$ to an additive category
with strict $(G,v)$-action which we will denote by $\cals(\cala)$.

\begin{definition}[$\cals(\cala)$] \label{def:cals(cala)}
An object in $\cals(\cala)$ is a pair $(A,g)$ consisting of an object
$A \in \cala$ and an element $g \in G$.  A morphism $(A,g)$ to $(B,h)$
is a morphism $\phi \colon R_g(A) \to R_h(B)$ in $\cala$.  The
composition of morphisms is given by the one in $\cala$. The category
$\cals(\cala)$ inherits the structure of an additive category from
$\cala$ in the obvious way.
\end{definition}

Next we define the structure of an additive category with strict $(G,v)$-action on
$\cals(\cala)$.  Define for $g\in G$ a functor $R^\cals_{g} \colon
\cals(\cala) \to \cals(\cala)$ of additive categories as follows.
Given an object $(A,h)$, define
$$R^{\cals}_{g}(A,h) = (A,hg).$$
Given a morphism $\phi \colon (A,h)
\to (B,k)$ define
$$R^{\cals}_{g}(\phi) \colon R^{\cals}_{g}(A,h) = (A,hg) \to
R^{\cals}_{g}(B,k) = (B,kg)$$
by the composite of morphisms in $\cala$
$$
R_{hg}(A) \xrightarrow{L_{h,g}(A)} R_{g}\left(R_{h}(A)\right)
\xrightarrow{R_{g}(\phi)} R_{g}\left(R_{k}(B)\right)
\xrightarrow{L_{k,g}(B)^{-1}} R_{kg}(B).
$$
if $v(g) = 1$
and
$$R^{\cals}_{g}(\phi) \colon
R^{\cals}_{g}(B,k) = (B,kg) \to R^{\cals}_{g}(A,h) = (A,hg)$$
by the composite of morphisms in $\cala$
$$R_{kg}(B)\xrightarrow{L_{k,g}(B)} R_{g}\left(R_{k}(B)\right)
\xrightarrow{R_{g}(\phi)} R_{g}\left(R_{h}(A)\right)
 \xrightarrow{L_{h,g}(A)^{-1}} R_{hg}(A)
$$
if $v(g) = -1$

A direct computation shows that $R^{\cals}_g$ is indeed a functor
of additive categories.  We conclude $R^{\cals}_e =
\id_{\cals(\cala)}$ from the conditions $R_e = \id$ and $L_{g,e} =
L_{e,g} = \id$.  We have to check $R^{\cals}_{g_2} \circ
R^{\cals}_{g_1} = R^{\cals}_{g_1g_2}$. We will do this for simplicity
only in the case $v(g_1)=v(g_2) = 1$, the other cases are analogous.
Given a morphism $\phi \colon
(A,h) \to (B,k)$, the morphism $R^{\cals}_{g_1g_2}(\phi)$ is given by
the composite in $\cala$
\begin{multline*}R_{hg_1g_2}(A)
  \xrightarrow{L_{h,g_1g_2}(A)} R_{g_1g_2}\left(R_hA)\right)
  \xrightarrow{R_{g_1g_2}(\phi)} R_{g_1g_2}\left(R_k(B)\right)
  \\
  \xrightarrow{L_{k,g_1g_2}(B)^{-1}} R_{kg_1g_2}(B).
\end{multline*}
The morphism $R^{\cals}_{g_2} \circ R^{\cals}_{g_1}(\phi)$ is given by
the composite in $\cala$
\begin{multline*}
  R_{hg_1g_2}(A) \xrightarrow{L_{hg_1,g_2}(A)}
  R_{g_2}\left(R_{hg_1}(A)\right)
  \\
  \xrightarrow{R_{g_2}(L_{h,g_1}(A))}
  R_{g_2}\left(R_{g_1}(R_h(A))\right)
  \xrightarrow{R_{g_2}\left(R_{g_1}(\phi)\right)}
  R_{g_2}\left(R_{g_1}(R_k(B))\right)
  \\
  \xrightarrow{R_{g_2}\left(L_{k,g_1}(B)^{-1}\right)}
  R_{g_2}\left(R_{kg_1}(B)\right) \xrightarrow{L_{kg_1,g_2}(B)^{-1}}
  R_{kg_1g_2}(B).
\end{multline*}
Next we compute that these two morphisms agree. Because of
condition~\ref{def:additive_category_with_weak_(G,v)-action:condition_for_(g,h,k)} in
Definition~\ref{def:additive_category_with_weak_(G,v)-action} have
\begin{eqnarray*}
R_{g_2}(L_{h,g_1}(A)) \circ L_{hg_1,g_2}(A) & = & L_{g_1,g_2}(R_h(A)) \circ L_{h,g_1g_2}(A);
\\
R_{g_2}(L_{k,g_1}(B)) \circ L_{kg_1,g_2}(B) & = & L_{g_1,g_2}(R_k(B)) \circ L_{k,g_1g_2}(B).
\end{eqnarray*}
Hence it suffices to show that the composite
\begin{multline*}
  R_{g_1g_2}\left(R_h(A)\right) \xrightarrow{L_{g_1,g_2}(R_h(A))}
  R_{g_2}\left(R_{g_1}(R_h(A))\right)
  \\
  \xrightarrow{R_{g_2}\left(R_{g_1}(\phi)\right)}
  R_{g_2}\left(R_{g_1}(R_k(B))\right)
  \xrightarrow{L_{g_1,g_2}(R_k(B))^{-1}} R_{g_1g_2}\left(R_k(B)\right)
\end{multline*}
agrees with
$$R_{g_1g_2}\left(R_k(B)\right) \xrightarrow{R_{g_1g_2}(\phi)}
R_{g_1g_2}\left(R_k(B)\right).$$
This follows from the fact that
$L_{g_1,g_2} \colon R_{g_1g_2} \to R_{g_2} \circ R_{g_2}$ is a natural
equivalence.

Let $(F,T) \colon \cala \to \calb$ be an $\epsilon$-variant functor of
additive categories with weak $(G,v)$-action.  It induces an
$\epsilon$-variant functor $\cals(F,T) \colon \cals(\cala) \to
\cals(\calb)$ of additive categories with strict $(G,v)$-action as
follows.  For simplicity we will only treat the case $\epsilon = 1$,
the other case $\epsilon = -1$ is analogous. The functor $\cals(F,T)$
sends an object $(A,h)$ in $\cals(\cala)$ to the object $(F(A),h)$ in
$\cals(\calb)$.  It sends a morphism $\phi \colon (A,h) \to (B,k)$ in
$\cals(\cala)$ which is given by a morphism $\phi \colon
R_h^{\cala}(A) \to R_k^{\cala}(B)$ in $\cala$ to the morphism
$\cals(F,T)(\phi) \colon (F(A),h) \to (F(B),k)$ in $\cals(\calb)$
which is given by the following composite of morphisms in $\calb$
$$
R^{\calb}_h(F(A)) \xrightarrow{T_h(A)^{-1}} F(R^{\cala}_h(A))
\xrightarrow{F(\phi)} F(R^{\cala}_k(B)) \xrightarrow{T_k(B)}
R^{\cala}_k(F(B)).
$$
We have to show $R^{\cals(\calb)}_g \circ \cals(F) = \cals(F) \circ
R^{\cals(\cala)}_g$ for every $g \in G$. We only treat the case $v(g) = 1$.
This is obvious on objects
since both composites send an object $(A,h)$ to $(F(A),hg)$. Let $\phi
\colon (A,h) \to (B,k)$ be a morphism in $\cals(\cala)$ which is given
by a morphism $\phi \colon R_h^{\cala}(A) \to R_k^{\cala}(B)$ in
$\cala$.  Then $R^{\cals(\calb)}_g \circ \cals(F)(\phi)$ is the
morphism $(F(A),hg) \to (F(B),kg)$ in $\cals(\calb)$ which is given by
the composite in $\calb$
\begin{multline*}
  R_{hg}^{\calb}(F(A)) \xrightarrow{L_{h,g}^{\calb}(F(A))}
  R_g^{\calb}(R_h^{\calb}(F(A)))
  \xrightarrow{R_g^{\calb}(T_h(A)^{-1})}
  R_g^{\calb}(F(R_h^{\cala}(A)))
  \\
  \xrightarrow{R_g^{\calb}(F(\phi))} R_g^{\calb}(F(R_k^{\cala}(B)))
  \xrightarrow{R_g^{\calb}(T_k(B))} R_g^{\calb}(R_k^{\calb}(F(B)))
  \xrightarrow{L_{k,g}^{\calb}(F(B))^{-1}} R_{kg}^{\calb}(F(B))
\end{multline*}
and $\cals(F) \circ R^{\cals(\cala)}_g(\phi)$ is the morphism
$(F(A),hg) \to (F(B),kg)$ in $\cals(\calb)$ which is given by the
composite in $\calb$
\begin{multline*}
  R_{hg}^{\calb}(F(A)) \xrightarrow{T_{hg}(A)^{-1}}
  F(R_{hg}^{\cala}(A)) \xrightarrow{F(L_{h,g}^{\cala}(A))}
  F(R_g^{\cala}(R_h^{\cala}(A)))
  \\
  \xrightarrow{F(R_g^{\cala}(\phi))} F(R_g^{\cala}(R_k^{\cala}(B)))
  \xrightarrow{F(L_{k,g}^{\cala}(B)^{-1})} F(R_{kg}^{\cala}(B))
  \xrightarrow{T_{kg}(A)} R_{kg}^{\calb}(F(B)).
\end{multline*}
Since $T_g\colon F\circ R_g^{\cala} \to R_g^{\calb} \circ F$ is a
natural transformation, the following diagram commutes
$$\xymatrix@!C=10em{ F(R_g^{\cala}(R_h^{\cala}(A)))
  \ar[r]^-{F(R_g^{\cala}(\phi))} \ar[d]^-{T_g(R_h^{\cala}(A))} &
  F(R_g^{\cala}(R_k^{\cala}(B))) \ar[d]^-{T_g(R_k^{\cala}(B))}
  \\
  R_g^{\calb}(F(R_h^{\cala}(A))) \ar[r]^-{R_g^{\calb}(F(\phi))} &
  R_g^{\calb}(F(R_k^{\cala}(B))) }
$$
Hence it suffices to show that the composite
\begin{multline*}R_{hg}^{\calb}(F(A))
  \xrightarrow{L_{h,g}^{\calb}(F(A))} R_g^{\calb}(R_h^{\calb}(F(A)))
  \xrightarrow{R_g^{\calb}(T_h(A)^{-1})}
  R_g^{\calb}(F(R_h^{\cala}(A)))
  \\
  \xrightarrow{T_g(R_h^{\cala}(A))^{-1}}
  F(R_g^{\cala}(R_h^{\cala}(A)))
\end{multline*}
agrees with the composite
$$
R_{hg}^{\calb}(F(A)) \xrightarrow{T_{hg}(A)^{-1}}
F(R_{hg}^{\cala}(A)) \xrightarrow{F(L_{h,g}^{\cala}(A))}
F(R_g^{\cala}(R_h^{\cala}(A)))
$$
and that the composite
\begin{multline*}
  F(R_g^{\cala}(R_k^{\cala}(B))) \xrightarrow{T_g(R_k^{\cala}(B))}
  R_g^{\calb}(F(R_k^{\cala}(B))) \xrightarrow{R_g^{\calb}(T_k(B))}
  R_g^{\calb}(R_k^{\calb}(F(B)))
  \\
  \xrightarrow{L_{k,g}^{\calb}(F(B))^{-1}} R_{kg}^{\calb}(F(B))
\end{multline*}
agrees with the composite
$$
F(R_g^{\cala}(R_k^{\cala}(B)))
\xrightarrow{F(L_{k,g}^{\cala}(B)^{-1})} F(R_{kg}^{\cala}(B))
\xrightarrow{T_{kg}(A)} R_{kg}^{\calb}(F(B)).
$$
This follows in both cases from the commutativity of the
diagram~\eqref{F_T_L_compatible}. This finishes the proof that
$\cals(F,T)$ is a functor of additive categories with strict $(G,v)$-action.

Let $S \colon (F_1,T_1) \to (F_2,T_2)$ be a natural transformation of
$\epsilon$-variant functors of additive categories with weak
$(G,v)$-action $(F_1,T_1) \colon \cala_1 \to \cala_2$ and $(F_2,T_2)
\colon \cala_1 \to \cala_2$.  It induces a natural transformation
$\cals(S) \colon \cals(F_1,T_1) \to \cals(F_2,T_2)$ of functors of
additive categories with strict $(G,v)$-action as follows.  Given an
object $(A,g)$ in $\cals(\cala)$, we have to specify a morphism in
$\cals(\cala)$
$$\cals(S)(A) \colon \cals(F_1,T_1)(A,g) = (F_1(A),g) \to
\cals(F_2,T_2)(A,g) = (F_2(A),g),$$
i.e., a morphism
$R^{\cala}_g(F_1(A)) \to R^{\cala}_g(F_2(A))$ in $\cala$. We take
$R^{\cala}_g(S(A))^{v(g)}$. We leave it to the reader to check that
this is indeed a natural transformation of $\epsilon$-variant functors
of additive categories with strict $(G,v)$-action using the
commutativity of the diagram~\eqref{F_I_T_add_G-cat}.

Let ${\Gaddcat}^{\epsilon}$ be the category of additive categories
with weak $(G,v)$-action with $\epsilon$-variant functors as morphisms
and let ${\sGaddcat}^{\epsilon} $ be the category of additive
categories with strict $(G,v)$-action with $\epsilon$-variant functors
as morphisms.  There is the forgetful functor
$$\forget \colon {\sGaddcat}^{\epsilon} \to {\Gaddcat}^{\epsilon}$$
and the functor constructed above
$$\cals \colon {\Gaddcat}^{\epsilon} \to {\sGaddcat}^{\epsilon}.$$

\begin{lemma}
\label{lem:adjoint_pair_(cals,forget)}
\begin{enumerate}

\item \label{lem:adjoint_pair_(cals,forget):adjoint}
We obtain an adjoint pair of functors $(\cals,\forget)$.

\item \label{lem:adjoint_pair_(cals,forget):equivalence}
We get for every additive category $\cala$ with weak $(G,v)$-action a functor
of additive categories with weak $(G,v)$-action
$$P_{\cala} \colon \cala \to \forget(\cals(\cala))$$
which is natural in $\cala$ and
whose underlying functor of additive categories is an equivalence
of additive categories.
\end{enumerate}
\end{lemma}
\begin{proof}
We will only treat the case, where $v$ is trivial and $\epsilon = 1$,
the other cases are analogous.
\\[2mm]\ref{lem:adjoint_pair_(cals,forget):adjoint}
We have to construct for any additive category $\cala$ with weak $G$-action and
any additive category $\calb$ with strict $G$-action
to one another inverse maps
$$\alpha \colon \func_{\sGaddcat}(\cals(\cala),\calb)
\to \func_{\Gaddcat}(\cala,\forget(\calb))$$
and
$$\beta \colon \func_{\Gaddcat}(\cala,\forget(\calb)) \to
\func_{\sGaddcat}(\cals(\cala),\calb).
$$
For a functor of additive categories with strict $G$-action
$F \colon \cals(\cala) \to \calb$, the functor of
additive categories with weak $G$-action, $\alpha(F) \colon \cala \to \forget(\calb)$
is given by a functor $\alpha(F) \colon \cala \to \forget(\calb)$
of additive categories and a collection of natural transformations
$T(F)_g \colon \alpha(F) \circ R_g^{\cala} \to R_g^{\calb} \circ \alpha(F)$
satisfying certain compatibility conditions. We first explain the
functor  $\alpha(F) \colon \cala \to \forget(\calb)$.
It sends a  morphism $f \colon A \to B$ in $\cala$ to the morphism in $\calb$
which is given by the value of $F$ on the morphism $(A,e) \to (B,e)$
in $\cals(\cala)$ defined by $f$.
For $g \in G$ the transformation $T(F)_g$ evaluated at an object $A$ in $\cala$
is the morphism
$$\alpha(F)(R_g^{\cala}(A)) = F(R_g^{\cala}(A),e) \to
R_g^{\calb}(\alpha(F)(A)) = R_g^{\calb}(F(A,e))$$
defined as follows. It is given by the composite of the image under $F$ of the
morphism  $(R^{\cala}_g(A),e) \to R_g^{\cals(\cala)}(A,e) = (A,g)$ in $\cals(\cala)$
which is defined by the identity morphism
$\id \colon R_g^{\cala}(A) \to R_g^{\cala}(A)$ in $\cala$ and the
identity $F(R_g^{\cals(\cala)}(A,e)) = R_g^{\calb}(F(A,e))$
which comes from the assumption that $F$ is a functor of strict
additive $G$-categories.
One easily checks that $\alpha(F)$ satisfies
condition~\eqref{F_T_L_compatible} since it is satisfied for $F$.

Given a functor of additive categories with weak $G$-action $(F,T)
\colon \cala \to \forget(\calb)$, the functor of additive categories
with strict $G$-action $\beta(F,T) \colon \cals(\cala) \to \calb$ is
defined as follows.  It sends an object $(A,h)$ to
$R_h^{\calb}(F(A))$. A morphism $\phi \colon (A,h) \to (B,k)$ in
$\cals(A)$ which is given by a morphism $\phi \colon R^{\cala}_hA \to
R^{\cala}_kB$ in $\cala$ is sent to morphism in $\calb$ given by the
composite
$$
R_h^{\calb}(F(A)) \xrightarrow{T_h(A)^{-1}} F(R^{\cala}_h(A))
\xrightarrow{F(\phi)} F(R^{\cala}_k(B)) \xrightarrow{T_k(B)}
R^{\calb}_k(F(B)).
$$
The following calculation shows that $\beta(F,T)$ is indeed a
functor of additive categories with strict $G$-action. Given an
element $g \in G$ the morphism $R_g^{\cals(\cala)}(\phi) \colon (A,hg)
\to (B,kg)$ in $\cals(\cala)$ is given by the morphism in $\cala$
$$R_{hg}^{\cala}(A) \xrightarrow{L_{hg}^{\cala}(A)}
R_g^{\cala}(R_h^{\cala}(A)) \xrightarrow{R_g^{\cala}(\phi)}
R_g^{\cala}(R_k^{\cala}(B)) \xrightarrow{L_{k,g}^{\cala}(B)^{-1}}
R_{kg}^{\cala}(B).
$$
Hence $\beta(F,T) \circ R_g^{\cals(\cala)}(\phi)$ is the morphism
in $\calb$ given by the composite
\begin{multline*}
  R_{hg}^{\calb}(F(A)) \xrightarrow{T_{hg}(A)^{-1}}
  F(R_{hg}^{\cala}(A)) \xrightarrow{F(L_{hg}^{\cala}(A))}
  F(R_g^{\cala}(R_h^{\cala}(A)))
  \\
  \xrightarrow{F(R_g^{\cala}(\phi))} F( R_g^{\cala}(R_k^{\cala}(B)))
  \xrightarrow{F(L_{k,g}^{\cala}(B)^{-1})} F(R_{kg}^{\cala}(B)))
  \\
  \xrightarrow{T_{kg}(B)} R_{kg}^{\cala}(F(B))).
\end{multline*}
The morphism $R_g^{\calb}(B) \circ \beta(F,T)(\phi)$ in $\calb$ is
given by the composite
\begin{multline*}
  R_g^{\calb}(R_h^{\calb}(F(A)))
  \xrightarrow{R_g^{\calb}(T_h(A)^{-1})}
  R_g^{\calb}(F(R^{\cala}_h(A))) \xrightarrow{R_g^{\calb}(F(\phi))}
  R_g^{\calb}(F(R^{\cala}_k(B)))
  \\
  \xrightarrow{R_g^{\calb}(T_k(B))} R_g^{\calb}(R^{\calb}_k(F(B))).
\end{multline*}
Since $\calb$ is a additive category with strict $G$-action by
assumption, we have the equalities $R_g^{\calb}(R_h^{\calb}(F(A))) =
R_{hg}^{\calb}(F(A))$ and $R_g^{\calb}(R^{\calb}_k(B)) =
R_{kg}^{\cala}(F(B)))$. We must show that under these identifications
the two morphisms in $\calb$ above agree.  Since $T_g$ is a natural
transformation $F \circ R^{\cala}_g \to R^{\calb}_g \circ F$, the
following diagram commutes
$$\xymatrix@!C=10em{ F(R_g^{\cala}(R_h^{\cala}(A)))
  \ar[r]^-{F(R_g^{\cala}(\phi))} \ar[d]^-{T_g(R_h^{\cala}(A))} &
  F(R_g^{\cala}(R_k^{\cala}(B))) \ar[d]^-{T_g(R_k^{\cala}(B))}
  \\
  R_g^{\calb}(F(R_h^{\cala}(A))) \ar[r]^-{R_g^{\calb}(F(\phi))} &
  R_g^{\calb}(F(R_k^{\cala}(B))) }
$$
Hence it suffices to show that the composites
$$R_g^{\calb}(R_h^{\calb}(F(A))) = R_{hg}^{\calb}(F(A))
\xrightarrow{T_{hg}(A)^{-1}} F(R_{hg}^{\cala}(A))
\xrightarrow{F(L_{hg}^{\cala}(A))} F(R_g^{\cala}(R_h^{\cala}(A)))
$$
and
$$R_g^{\calb}(R_h^{\calb}(F(A)))
\xrightarrow{R_g^{\calb}(T_h(A)^{-1})} R_g^{\calb}(F(R^{\cala}_h(A)))
\xrightarrow{T_g(R_h^{\cala}(A))^{-1}} F(R_g^{\cala}(R_h^{\cala}(A)))
$$
agree and that the composites
$$
F(R_g^{\cala}(R_k^{\cala}(B)))
\xrightarrow{F(L_{k,g}^{\cala}(B)^{-1})} F(R_{kg}^{\cala}(B)))
\xrightarrow{T_{kg}(B)} R_{kg}^{\cala}(F(B))) =
R_g^{\calb}(R^{\calb}_k(F(B)))$$
and
$$F( R_g^{\cala}(R_k^{\cala}(B))) \xrightarrow{T_g(R_k^{\cala}(B))}
R_g^{\calb}(F(R_k^{\cala}(B))) \xrightarrow{R_g^{\calb}(T_k(B))}
R_g^{\calb}(R^{\calb}_k(F(B))$$
agree. This follows in both cases from
the commutativity of the diagram~\eqref{F_T_L_compatible}. This
finishes the proof that $\beta(F)$ is a functor of additive categories
with strict $G$-action.  We leave it to the reader to check that both
composites $\beta \circ \alpha$ and $\alpha \circ \beta$ are the
identity.  \\[1mm]\ref{lem:adjoint_pair_(cals,forget):equivalence} The
in $\cala$ natural functor of additive categories with weak
$(G,v)$-action
$$P_{\cala} \colon \cala \to \forget(\cals(\cala))$$
is defined to be
the adjoint of the identity functor $\id \colon \cals(\cala) \to
\cals(\cala)$. Explicitly it sends an object $A$ to the object $(A,e)$
and a morphism $\phi \colon A \to B$ to the morphism $(A,e) \to (B,e)$
given by $\phi$. Obviously $P_{\cala}$ induces a bijection
$\mor_{\cala}(A,B) \to \mor_{\cals(\cala)}(P_{\cala}(A),P_{\cala}(B))$
and for every object $(A,g)$ in $\cals(\cala)$ there is an object in
the image of $P_{\cala}$ which is isomorphic to $(A,g)$, namely,
$P_{\cala}(R^{\cala}_{g}(A)) = (R_{g}^{\cala}(A),e)$.  Hence the
underlying functor $R_{\cala}$ is an equivalence of additive
categories.
\end{proof}

%%%%%%%%%%%%%%%%%%%%%%%%%%%%%%%%%%%%%%%%%%%%%%%%%%%%%%%%%%%%%%%%%%%%%%%%%%%%%%%%%%%%%%%%%
%%%%%%%%%%%%%%%%%%%%%%% Crossed product rings and involutions %%%%%%%%%%%%%%%%%%%%%%%%%%%%%%%%
%%%%%%%%%%%%%%%%%%%%%%%%%%%%%%%%%%%%%%%%%%%%%%%%%%%%%%%%%%%%%%%%%%%%%%%%%%%%%%%%%%%%%%%%%

\typeout{-------------------------  Section 4:  Crossed product rings and involutions
  ----------------------}

\section{Crossed product rings and involutions}
\label{sec:Crossed_product_rings_and_involutions}

In this subsection we will introduce the concept of  a crossed product ring.
Let $R$ be a ring and let $G$ be a group. Let $e\in G$ be the unit in $G$ and
denote by $1$ the multiplicative unit in $R$. Suppose that we are given
maps of sets
\begin{eqnarray}
c\colon G & \to & \aut(R), \quad g \mapsto c_g
\label{c:G_toaut(G)};
\\
\tau\colon G \times G & \to & R^{\times}.
\label{tau:G_times_G_to_R_times}
\end{eqnarray}
We require
\begin{eqnarray}
c_{\tau(g,g')} \circ c_{gg'} & = & c_{g} \circ c_{g'};
\label{tau_and_c_group_homo}
\\
\tau (g,g') \cdot \tau (gg',g'') & = & c_{g}(\tau (g',g''))\cdot \tau (g,g'g'');
\label{tau_and_c_cocykel}
\\ c_e & = & \id_R; \label{c_e-Is_id}
\\
\tau(e,g) & = & 1; \label{tau(g,e)_is_1}
\\
\tau(g,e) & = & 1, \label{tau(e,g)_is_1}
\end{eqnarray}
for $g,g',g'' \in G$, where $c_{\tau(g,g')}\colon R \to R$ is
conjugation with $\tau(g,g')$, i.e., it sends $r$ to
$\tau(g,g')r\tau(g,g')^{-1}$.  Let $R \ast G = R \ast_{c,\tau} G$ be
the free $R$-module with the set $G$ as basis. It becomes a ring with
the following multiplication
$$\left(\sum_{g \in G} \lambda_g g\right) \cdot \left(\sum_{h \in G}
  \mu_h h\right) = \sum_{g \in G} \left(\sum_{\substack{g',g'' \in
      G,\\g'g'' = g}} \lambda_{g'} c_{g'}(\mu_{g''}) \tau(g',g'')
\right) g.$$
This multiplication is uniquely determined by the
properties $g\cdot r = c_g(r)\cdot g$ and $g \cdot g' = \tau(g,g')
\cdot (gg')$. The conditions~\eqref{tau_and_c_group_homo}
and~\eqref{tau_and_c_cocykel} relating $c$ and $\tau$ are equivalent
to the condition that this multiplication is associative. The other
conditions~\eqref{c_e-Is_id}, \eqref{tau(g,e)_is_1}
and~\eqref{tau(e,g)_is_1} are equivalent to the condition that the
element $1 \cdot e$ is a multiplicative unit in $R\ast G$. We call
\begin{eqnarray}
R \ast G & = & R \ast_{c,\tau} G
\label{R_ast_G}
\end{eqnarray}
the \emph{crossed product} of $R$ and $G$ with respect to $c$ and
$\tau$.

\begin{example} \label{exa:RG_as_cros._prod.}
Let $1 \to H \xrightarrow{i} G \xrightarrow{p} Q \to 1$ be an extension of groups.
Let $s\colon Q \to G$ be a map satisfying $p \circ s = \id$ and $s(e) = e$. We do not
require $s$ to be a group homomorphism.
Define $c\colon Q \to \aut(RH)$ by $c_q(\sum_{h \in H} \lambda_h h) =
\sum_{h \in H} \lambda_h s(q)hs(q)^{-1}$. Define $\tau\colon Q \times Q \to
(RH)^{\times}$ by $\tau (q,q') = s(q)s(q')s(qq')^{-1}$. Then we obtain a ring
isomorphism $RH \ast Q \to RG$ by sending $\sum_{q \in Q} \lambda_q q$
to $\sum_{q \in Q} i(\lambda_q) s(q)$, where $i\colon RH \to RG$ is the
ring homomorphism induced by $i\colon H \to G$. Notice that $s$ is a group
homomorphism if and only if $\tau$ is constant with value $1 \in R$.
\end{example}

Next we consider the additive category with involution $\FGP{R}$
of finitely generated projective $R$-modules. For $g \in G$ we obtain a functor
$\res_{c_g} \colon \FGP{R} \to \FGP{R}$ by restriction with the ring automorphism
$c_g\colon R \to R$. Define natural transformation of functors $\FGP{R} \to \FGP{R}$
$$L_{\tau(g,h)} \colon \res_{c_{gh}} \to \res_{c_h} \circ \res_{c_g}$$
by assigning to a finitely generated projective $R$-module the $R$-homomorphism
$$\res_{c_{gh}}P \to \res_{c_h}\res_{c_g}P,\quad p \mapsto \tau(g,h) p.$$
This is indeed a $R$-linear map because of the following computation
for $r \in R$ and $p \in P$
\begin{multline*}
\tau(g,h)c_{gh}(r) = \tau(g,h)c_{gh}(r)\tau(g,h)^{-1}\tau(g,h) =
c_{\tau(g,h)} \circ c_{gh}(r)\tau(g,h)
\\
= c_g \circ c_h(r)\tau(g,h).
\end{multline*}

\begin{lemma} \label{lem:weak_G-structure_on_FGP(R)}
We get from the collections $\{\res_{c_g} \mid g \in G\}$ and
$\{L_{\tau(g,h)} \mid g,h \in G\}$ the structure of an additive category with
weak $G$-action on $\FGP{R}$.
\end{lemma}
\begin{proof}
Condition~\eqref{tau_and_c_cocykel} implies
that for every finitely generated projective $R$-module the composites
$$\res_{c_{gg'g''}} P \xrightarrow{L_{\tau(g,g'g'')}}
  \res_{c_{g'g''}}\res_{c_g} P
  \xrightarrow{L_{c_g(\tau(g',g''))}}
  \res_{c_{g''}} \res_{C_{g'}} \res_{c_g} P$$
and
$$\res_{c_{gg'g''}} P \xrightarrow{L_{\tau(gg',g'')}}
  \res_{c_{g''}}\res_{c_{gg'}} P
  \xrightarrow{L_{\tau(g,g')}}
  \res_{c_{g''}} \res_{C_{g'}} \res_{c_g} P$$
agree. This takes care of
condition~\ref{def:additive_category_with_weak_(G,v)-action:condition_for_(g,h,k)}
in Definition~\ref{def:additive_category_with_weak_(G,v)-action}.
We conclude $(\res_{c(e)} = \id$, $L_{\tau(g,e)} = \id$ and
$L_{\tau(e,g)} = \id$ for all $g \in G$ from~\eqref{c_e-Is_id},
\eqref{tau(g,e)_is_1} and~\eqref{tau(e,g)_is_1}.
\end{proof}

 Because of Lemma~\ref{lem:weak_G-structure_on_FGP(R)} we
obtain two additive categories with strict $G$-action from the constructions of
Section~\ref{sec:Making_an_additive_categories_with_weak_(G,v)-action_strict}

\begin{eqnarray}
\FGP{R}_{c,\tau} & := & \cals(\FGP{R});
\label{FGP(R)_(c,tau)}
\end{eqnarray}

>From now on assume that $R$ comes with an involution of rings $r
\mapsto \overline{r}$. We want to consider extensions of it to an
involution on $R \ast G$. Suppose that additionally we are given a map
\begin{eqnarray}
&w \colon G  \to  R. &
\label{map_w}
\end{eqnarray}
We require the following conditions for $g,h \in G$ and $r \in R$
\begin{eqnarray}
w(e) & = & 1;
\label{w(e)_is_1}
\\
w(gh)
& = &
w(h)c_{h^{-1}}(w(g))\tau (h^{-1},g^{-1})
      c_{(gh)^{-1}}\left(\overline{\tau(g,h)}\right)^{-1};
\label{w(gh)}
\\
\overline{w(g)} & = & w(g)c_{g}^{-1}\left(\tau (g,g^{-1})\overline{\tau (g,g^{-1})}^{-1}\right);
\label{overlinew(g)}
\\
\overline{c_{g}(r)}
& = &
c_g\left(\left(w(g)\tau(g^{-1},g)\right)^{-1}\overline{r}\left(w(g)\tau(g^{-1},g)\right)\right)
\label{overlinec_g(r)}.
\end{eqnarray}
We claim that there is precisely one involution on $R \ast G$ with the properties
that it extends the involution on $R$ and sends $g$ to $w(g) \cdot g^{-1}$.
The candidate for the involution is
\begin{eqnarray}
\overline{\sum_{g \in G} r_g \cdot g}
& := &
\sum_{g \in G} w(g)c_{g^{-1}}(\overline{r_g}) \cdot g^{-1}.
\label{inv_onR_ast_G}
\end{eqnarray}
One easily concludes from the requirements and the axioms
of an involution that this is the only possible formula for such an involution. Namely,
\begin{multline*}
\overline{\sum_{g \in G} r_g \cdot g}
= \sum_{g \in G} \overline{r_g \cdot g}
= \sum_{g \in G} \overline{g} \cdot\overline{r_g}
= \sum_{g \in G} w(g) \cdot g^{-1} \cdot \overline{r_g}
\\
= \sum_{g \in G} w(g) \cdot \left(g^{-1} \cdot \overline{r_g} \cdot g\right) \cdot g^{-1}
= \sum_{g \in G} (w(g)c_{g^{-1}}(\overline{r_g}) \cdot g^{-1}.
\end{multline*}
Before we explain that this definition indeed satisfies the axioms for
an involution, we show that the conditions about $w$ above are
necessary for this map to be an involution on $R \ast G$. So assume
that we have an involution on $R \ast G$ that extends the involution
on $R$ and sends $g$ to $w(g) \cdot g^{-1}$ for a given map $w \colon
G \to R$.  Denote by $1$ the multiplicative unit in both $R$ and $R
\ast G$.  From
$$1 \cdot e = 1 =  \overline{1} = \overline{1  \cdot e}= w(e) \cdot e$$
we conclude~\eqref{w(e)_is_1}.
The equality
\begin{multline*}
w(gh) c_{(gh)^{-1}}\left(\overline{\tau(g,h)}\right) \cdot (gh)^{-1}
= \overline{\tau(g,h) \cdot gh}
= \overline{g\cdot h}
\\
= \overline{h} \cdot \overline{g}
= w(h) \cdot h^{-1} \cdot w(g) \cdot g^{-1}
= w(h)\left(h^{-1} \cdot w(g) \cdot h\right) \cdot h^{-1} \cdot g^{-1}
\\
= w(h)c_{h^{-1}}(w(g))\tau (h^{-1},g^{-1}) \cdot (gh)^{-1}
\end{multline*}
implies~\eqref{w(gh)}. If we take $h = g^{-1}$ in~\eqref{w(gh)} and use~\eqref{w(e)_is_1}, we get
\begin{eqnarray}
& 1 = w(e) = w(gg^{-1}) =  w(g^{-1})c_{g}(w(g))\tau (g,g^{-1}) \overline{\tau(g,g^{-1})}^{-1}.&
\label{eqn_about_1_is_w(e)_is}
\end{eqnarray}
This implies that for all $g \in G$ the element $w(g)$ is a unit in $R$ with inverse
$$w(g)^{-1} = c_{g^{-1}}(w(g^{-1}))\tau (g^{-1},g)
\overline{\tau(g^{-1},g)}^{-1}.$$
The equality
\begin{multline*}
g =  \overline{\overline{g}}
= \overline{w(g) \cdot g^{-1}}
= \overline{g^{-1}} \cdot \overline{w(g)}
= w(g^{-1}) \cdot g \cdot \overline{w(g)}
\\
= w(g^{-1})\cdot \left(g\cdot \overline{w(g)}\cdot g^{-1}\right) \cdot g
= w(g^{-1})c_{g}\left(\overline{w(g)}\right) \cdot g
\end{multline*}
together with~\eqref{eqn_about_1_is_w(e)_is} implies
$$w(g^{-1})c_{g}\left(\overline{w(g)}\right) = 1 =
w(g^{-1})c_{g}(w(g))\tau (g,g^{-1})
\overline{\tau(g,g^{-1})}^{-1}.$$ If we multiply this equation
with $w(g^{-1})^{-1}$ and apply the inverse $c_g^{-1}$ of $c_g$,
we derive condition~\eqref{overlinew(g)}. The equality
\begin{multline*}
\overline{r} \cdot w(g) \cdot g^{-1}
= \overline{r} \cdot \overline{g}
= \overline{g\cdot r}
= \overline{(g\cdot r \cdot g^{-1}) \cdot g}
= \overline{c_{g}(r)\cdot g}
=
\overline{g} \cdot \overline{c_{g}(r)}
\\
=
w(g) \cdot g^{-1} \cdot \overline{c_{g}(r)}
=
w(g) \cdot \left(g^{-1} \cdot \overline{c_{g}(r)} \cdot g\right) \cdot g^{-1}
=
w(g) \cdot c_{g^{-1}}\left(\overline{c_{g}(r)}\right) \cdot g^{-1}
\end{multline*}
implies that for all $g \in G$ and $r \in R$ we have
$\overline{r} \cdot w(g) = w(g) \cdot c_{g^{-1}}\left(\overline{c_{g}(r)}\right)$
and hence
$$\overline{c_{g}(r)} = c_{g^{-1}}^{-1}\left(w(g)^{-1}\overline{r}w(g)\right).$$
>From the relation~\eqref{tau_and_c_group_homo}
we conclude $c_{\tau(g^{-1},g)} = c_{g^{-1}} \circ c_{g}$
and hence $ c_{g^{-1}}^{-1} = c_g \circ c_{\tau(g^{-1},g)}^{-1}$.
Now condition~\eqref{overlinec_g(r)} follows.

Finally we show that the conditions~\eqref{w(e)_is_1}, \eqref{w(gh)},
\eqref{overlinew(g)} and~\eqref{overlinec_g(r)} on $w$ do imply that
we get an involution of rings on $R \ast G$ by the
formula~\eqref{inv_onR_ast_G}.  Obviously this formula is compatible
with the additive structure on $R \ast G$ and sends $1$ to $1$. In
order to show that it is an involution and compatible with the
multiplicative structure we have to show $\overline{g\cdot h} =
\overline{h}\cdot \overline{g}$, $\overline{rs} = \overline{s} \cdot
\overline{r}$, $\overline{r \cdot g} = \overline{g} \cdot
\overline{r}$, $\overline{g \cdot r} = \overline{r} \cdot
\overline{g}$, $\overline{\overline{r}} = r$ and
$\overline{\overline{g}} = g$ for $r,s \in R$ and $g,h \in G$. We get
$\overline{rs} = \overline{s} \cdot \overline{r}$ and
$\overline{\overline{r}} = r$ from the fact that we start with an
involution on $R$. The other equations follow from the proofs above
that~\eqref{inv_onR_ast_G} is the only possible candidate and that the
conditions about $w$ are necessary for the existence of the desired
involution on $R \ast G$, just read the various equations and
implications backwards.  We will denote the resulting ring with
involution by
\begin{eqnarray}
&R \ast_{c,\tau,w} G. &
\label{r_ast_c,tau,w_G}
\end{eqnarray}

\begin{example} \label{exa:RG_as_cros._prod.with_inv}
Suppose that we are in the situation of Example~\ref{exa:RG_as_cros._prod.}.
Suppose that we are additionally given a group homomorphism
$w_1 \colon G \to \cent(R)^{\times}$ to the abelian group of
invertible central elements  in $R$ satisfying
$\overline{w_1(g)} = w_1(g)$ for all $g \in G$. The $w_1$-twisted
involution on $RG$ is defined by $\overline{\sum_{g \in G} r_g \cdot g} =
\sum_{g \in G} \overline{r_g}w_1(g) \cdot g^{-1}$. It extends the
$w_1|_H$-involution on $RH$.
We obtain an involution on $RH \ast Q$ if we conjugate the $w_1$-twisted
involution with the isomorphism $RH \ast Q \xrightarrow{\cong} RG$
which we have introduced in Example~\ref{exa:RG_as_cros._prod.}.
This involution on $RH \ast Q$  sends $q \in Q$ to the element
$w_1(s(q))\tau(q^{-1},q)^{-1} \cdot q^{-1}$ because of the
following calculation in $RG$ for $q \in Q$
\begin{multline*}
\overline{s(q)}
= w_1(s(q)) \cdot s(q)^{-1}
= w_1(s(q)) \cdot s(q)^{-1}\cdot s(q^{-1})^{-1} \cdot s(q^{-1})
\\
= w_1(s(q)) \cdot \left(s(q^{-1}) \cdot s(q)\right)^{-1} \cdot s(q^{-1})
= w_1(s(q)) \cdot \left(\tau(q^{-1},q)s(q^{-1}q)\right)^{-1} \cdot s(q^{-1})
\\
= w_1(s(q))\tau(q^{-1},q)^{-1} \cdot s(q^{-1}).
\end{multline*}
Define
$$w \colon Q \to RH, \quad q \mapsto w_1(s(q))\tau(q^{-1},q)^{-1}.$$
Then $w$ satisfies the conditions \eqref{w(e)_is_1}, \eqref{w(gh)},
\eqref{overlinew(g)} and~\eqref{overlinec_g(r)} and the involution
on $RH \ast Q$ determined by $w$ corresponds under the
isomorphism $RH \ast Q \xrightarrow{\cong} RG$ to the $w_1$-twisted
involution on $RG$.
\end{example}

Let
\begin{eqnarray}
& t_g \colon \res_{c_g} \circ I_{\FGP{R}} \to I_{\FGP{R}} \circ \res_{c_g} &
\label{t_g}
\end{eqnarray}

be the natural transformation which assigns to a finitely generated projective $R$-module
$P$ the $R$-isomorphism
$t_g(P) \colon  \res_{c_g} P^* \to (\res_{c_g} P)^*$
which sends the $R$-linear map $f \colon P \to R$ to
the $R$-linear map
$$t_g(P)(f) \colon \res_{c_g} P \to R,
\quad p \mapsto c_g^{-1}(f(p))\left(w(g)\tau(g^{-1},g)\right)^{-1}.$$
We firstly check that $t_g(P)(f) \colon \res_{c_g} P \to R$
is $R$-linear by the following computation
\begin{eqnarray*}
t_g(P)(f)(c_g(r)p)
& = &
c_g^{-1}\left(f(c_g(r)p)\right)\left(w(g)\tau(g^{-1},g)\right)^{-1}
\\
& = &
c_g^{-1}\left(c_g(r)f(p)\right)\left(w(g)\tau(g^{-1},g)\right)^{-1}
\\
& = &
c_g^{-1}\left(c_g(r)\right)c_g^{-1}(f(p))\left(w(g)\tau(g^{-1},g)\right)^{-1}
\\
& = &
rc_g^{-1}(f(p))\left(w(g)\tau(g^{-1},g)\right)^{-1}
\\
& = &
rt_g(P)(f)(p).
\end{eqnarray*}
Finally we check that $t_g(P) \colon  \res_{c_g} P^* \to (\res_{c_g} P)^*$
is $R$-linear by the following calculation for $f \in P^*$ and $p \in P$
\begin{eqnarray*}
\lefteqn{t_g(P)\left((c_g(r) f)\right)(p)}
& &
\\
& = &
c_g^{-1}\left((c_g(r) f)(p)\right)\left(w(g)\tau(g^{-1},g)\right)^{-1}
\\
& = &
c_g^{-1}\left(f(p)\overline{c_g(r)}\right)\left(w(g)\tau(g^{-1},g)\right)^{-1}
\\
& = &
c_g^{-1}(f(p))c_g^{-1}(\overline{c_g(r)})\left(w(g)\tau(g^{-1},g)\right)^{-1}
\\
& = &
c_g^{-1}(f(p))c_g^{-1}\left(c_g\left(\left(w(g)\tau(g^{-1},g)\right)^{-1}
\overline{r}\left(w(g)\tau(g^{-1},g)\right)\right)\right)
\left(w(g)\tau(g^{-1},g)\right)^{-1}
\\
& = &
c_g^{-1}(f(p)\left(w(g)\tau(g^{-1},g)\right)^{-1}\overline{r}\left(w(g)\tau(g^{-1},g)\right)
\left(w(g)\tau(g^{-1},g)\right)^{-1}
\\
& = &
c_g^{-1}(f(p)\left(w(g)\tau(g^{-1},g)\right)^{-1}\overline{r}
\\
& = &
t_g(P)(f)(p)\overline{r}
\\
& = &
\left(rt_g(P)\right)(f)(p).
\end{eqnarray*}

\begin{definition}\label{def:additive_G-category_with_involution}
  An \emph{additive $G$-category with involution} $\cala$ is an
  additive $G$-category, which is the same as an additive category with
  strict $G$-action (see
  Definition~\ref{def:additive_category_with_weak_(G,v)-action}),
  together with an involution $(I,E)$ of additive categories
  (see~\eqref{involution_I_on_an_additive_category} and
  \eqref{E_belonging_to_I}) with the following properties:
  $I  \colon \cala \to \cala$ is a contravariant functor of additive
  $G$-categories, i.e., $R_g \circ I = I \circ R_g$ for all $g \in G$,
  and $E \colon \id_{\cala} \to I \circ I$ is a natural transformation
  of functors of additive $G$-categories, i.e., for every $g \in G$
  and every object $A$ in $\cala$ the morphisms $E(R_g(A))$ and
  $R_g(E(A))$ from $R_g(A)$ to $I^2 \circ (R_g(A) = R_g \circ I^2(A)$
  agree.
\end{definition}

\begin{lemma} \label{lem:involution_on_FGP(R_c,w)}
  The additive category with strict $G$-action $\FGP{R}_{c,\tau}$ of~\eqref{FGP(R)_(c,tau)}
  inherits the structure of an additive $G$-category with involution
  in the sense of Definition~\ref{def:additive_G-category_with_involution}.
\end{lemma}
\begin{proof}
  We firstly show that
  $$I_{\FGP{R}} \colon \FGP{R} \to \FGP{R}$$
  together with the
  collection of the $\{t_g^{-1} \colon I_{\FGP{R}} \circ \res_{c_g}
  \to \res_{c_g} \circ I_{\FGP{R}} \mid g \in G\}$ (see~\eqref{t_g})
  is a contravariant
  functor of additive categories with weak $G$-action. We have to
  verify that the diagram~\eqref{F_T_L_compatible} commutes. This is
  equivalent to show for every finitely generated projective
  $R$-module $P$ and $g,h \in G$ that the following diagram commutes
$$
\xymatrix@!C=10em{\res_{c_{gh}} P^* \ar[rr]^-{t_{gh}(P)} \ar[d]^-{L_{\tau(g,h)}(P^*)}
& &
(\res_{c_{gh}} P)^*
\\
\res_{c_{h}} \res_{c_{g}} P^*  \ar[r]^-{\res_{c_h} t_g(P)}
&
\res_{c_{h}} (\res_{c_{g}} P)^* \ar[r]^-{t_h(\res_gP)}
& (\res_{c_{h}} \res_{c_{g}} P)^* \ar[u]^-{L_{\tau(g,h)}(P)^*}
}
$$
We start with an element $f \colon P \to R$ in the left upper corner.
Its image under the upper horizontal arrow is
$p \mapsto c_{gh}^{-1}(f(p))\left(w(gh)\tau((gh)^{-1},gh)\right)^{-1}$.
Next we list successively how its image looks like if we go in the anticlockwise direction
from the left upper corner to the right upper corner. We first get $p
\mapsto f(p)\overline{\tau (g,h)}$.
After the second map we get
$p \mapsto c_g^{-1}\left(f(p)\overline{\tau(g,h)}\right)\left(w(g)\tau(g^{-1},g)\right)^{-1}$.
After applying the third map we obtain
$p \mapsto c_h^{-1}\left(c_g^{-1}\left(f(p)\overline{\tau(g,h)}\right)
\left(w(g)\tau(g^{-1},g)\right)^{-1}\right)
\left(w(h)\tau(h^{-1},h)\right)^{-1}$.
Finally we get
$p \mapsto c_h^{-1}\left(c_g^{-1}\left(f(\tau(g,h)p)\overline{\tau(g,h)}\right)
\left(w(g)\tau(g^{-1},g)\right)^{-1}\right) \left(w(h)\tau(h^{-1},h)\right)^{-1}$.
Since $f$ lies in $P^*$, we have $f(\tau(g,h)p) = \tau(g,h)f(p)$.
Hence it suffices to show for all $r \in R$
\begin{multline*}
c_h^{-1}\left(c_g^{-1}\left(\tau(g,h)r\overline{\tau(g,h)}\right)
\left(w(g)\tau(g^{-1},g)\right)^{-1}\right)
\left(w(h)\tau(h^{-1},h)\right)^{-1}
\\
= c_{gh}^{-1}(r)\left(w(gh)\tau((gh)^{-1},gh)\right)^{-1}.
\end{multline*}
(Notice that now $f$ has been eliminated.)
By applying $c_{gh}$ we see that this is equivalent to showing
\begin{multline*}
c_{gh}\left(c_h^{-1}\left(c_g^{-1}\left(\tau(g,h)r\overline{\tau(g,h)}\right)\right)\right)
\\
= rc_{gh}\left(\left(w(gh)\tau((gh)^{-1},gh)\right)^{-1}\left(w(h)\tau(h^{-1},h)\right)
c_h^{-1}\left(w(g)\tau(g^{-1},g)\right)\right).
\end{multline*}
>From the relation~\eqref{tau_and_c_group_homo}
we conclude that
$c_{gh} \circ c_{h^{-1}} \circ c_{g^{-1}}(s) = \tau(g,h)^{-1}s\tau(g,h)$ holds for
all $s \in R$. Hence it remains to show
\begin{multline*}
 \tau(g,h)^{-1}\left(\tau(g,h)r\overline{\tau(g,h)}\right) \tau(g,h)
\\
= rc_{gh}\left(\left(w(gh)\tau((gh)^{-1},gh)\right)^{-1}w(h)\tau(h^{-1},h)
c_h^{-1}\left(w(g)\tau(g^{-1},g)\right)\right).
\end{multline*}
This reduces to proving for $g,h \in G$
\begin{multline*}
\overline{\tau(g,h)} \tau(g,h)
\\
= c_{gh}\left(\tau((gh)^{-1},gh)^{-1}w(gh)^{-1}
w(h)\tau(h^{-1},h) c_h^{-1}\left(w(g)\tau(g^{-1},g)\right)\right).
\end{multline*}
(Notice that now $r$ has been eliminated.)
By inserting condition~\eqref{w(gh)} and the conclusions
$c_{\tau(h^{-1},h)} \circ c_h^{-1} = c_{h^{-1}}$ and
$c_{\tau((gh)^{-1},gh)} \circ c_{gh}^{-1} = c_{(gh)^{-1}}$ from
conditions~\eqref{tau_and_c_group_homo} and~\eqref{c_e-Is_id} we get
\begin{eqnarray*}
\lefteqn{w(gh)^{-1}
w(h)\tau(h^{-1},h) c_h^{-1}\left(w(g)\tau(g^{-1},g)\right)}
& &
\\
& = &
\left(w(h)c_{h^{-1}}(w(g))\tau (h^{-1},g^{-1})
      c_{(gh)^{-1}}\left(\overline{\tau(g,h)}\right)^{-1}\right)^{-1}w(h)
\\
& & \hspace{30mm}  \tau(h^{-1},h) c_h^{-1}\left(w(g)\tau(g^{-1},g)\right)
\tau(h^{-1},h)^{-1} \tau(h^{-1},h)
\\
& = &
c_{(gh)^{-1}}\left(\overline{\tau(g,h)}\right)\tau (h^{-1},g^{-1})^{-1}
c_{h^{-1}}(w(g))^{-1}w(h)^{-1}w(h)
\\
& & \hspace{30mm}
c_{h^{-1}}\left(w(g)\tau(g^{-1},g)\right)\tau(h^{-1},h)
\\
& = &
c_{(gh)^{-1}}\left(\overline{\tau(g,h)}\right)\tau (h^{-1},g^{-1})^{-1}
c_{h^{-1}}(w(g))^{-1}c_{h^{-1}}(w(g))
\\
& & \hspace{30mm}
c_{h^{-1}}\left(\tau(g^{-1},g)\right) \tau(h^{-1},h)
\\
& = &
c_{(gh)^{-1}}\left(\overline{\tau(g,h)}\right)\tau (h^{-1},g^{-1})^{-1}
c_{h^{-1}}\left(\tau(g^{-1},g)\right) \tau(h^{-1},h)
\\
& = &
\tau((gh)^{-1},gh) c_{gh}^{-1}\left(\overline{\tau(g,h)}\right)
\tau((gh)^{-1},gh)^{-1}\tau (h^{-1},g^{-1})^{-1}
\\
& & \hspace{30mm}
c_{h^{-1}}\left(\tau(g^{-1},g)\right) \tau(h^{-1},h).
\end{eqnarray*}
This implies
\begin{eqnarray*}
\lefteqn{c_{gh}\left(\tau((gh)^{-1},gh)^{-1}w(gh)^{-1}
w(h)\tau(h^{-1},h) c_h^{-1}\left(w(g)\tau(g^{-1},g)\right)\right)}
& &
\\
& = &
c_{gh}\left(\tau((gh)^{-1},gh)^{-1}\tau((gh)^{-1},gh) c_{gh}^{-1}
\left(\overline{\tau(g,h)}\right)\tau((gh)^{-1},gh)^{-1}\right.
\\
& & \hspace{30mm}
\left.\tau (h^{-1},g^{-1})^{-1}c_{h^{-1}}\left(\tau(g^{-1},g)\right) \tau(h^{-1},h) \right)
\\
& = &
\overline{\tau(g,h)}c_{gh}\left(\tau((gh)^{-1},gh)^{-1}\tau (h^{-1},g^{-1})^{-1}
c_{h^{-1}}\left(\tau(g^{-1},g)\right) \tau(h^{-1},h) \right).
\end{eqnarray*}
Hence it remains to show
$$\tau(g,h) = c_{gh}\left(\tau((gh)^{-1},gh)^{-1}\tau (h^{-1},g^{-1})^{-1}
c_{h^{-1}}\left(\tau(g^{-1},g)\right) \tau(h^{-1},h) \right).$$
(Notice that we have eliminated any expression involving the involution.)
>From condition~\eqref{tau_and_c_group_homo}, \eqref{tau_and_c_cocykel}
and~\eqref{c_e-Is_id}  we conclude
\begin{eqnarray*}
\tau(h^{-1},g^{-1}) \tau((gh)^{-1},g) & = & c_{h^{-1}}(\tau(g^{-1},g));
\\
\tau(gh)^{-1},g)\tau(h^{-1},h) & = &c_{(gh)^{-1}}(\tau(g,h))\tau((gh)^{-1},gh);
\\
c_{gh}^{-1} & = & c_{\tau((gh)^{-1},gh)^{-1}} \circ c_{(gh)^{-1}}.
\end{eqnarray*}
Hence
\begin{eqnarray*}
\lefteqn{\tau((gh)^{-1},gh)^{-1}\tau (h^{-1},g^{-1})^{-1}
c_{h^{-1}}\left(\tau(g^{-1},g)\right) \tau(h^{-1},h)}
& &
\\
& = &
\tau((gh)^{-1},gh)^{-1}\tau (h^{-1},g^{-1})^{-1}
\tau(h^{-1},g^{-1}) \tau((gh)^{-1},g)  \tau(h^{-1},h)
\\ & = &
\tau((gh)^{-1},gh)^{-1}\tau((gh)^{-1},g)  \tau(h^{-1},h)
\\ & = &
\tau((gh)^{-1},gh)^{-1}c_{(gh)^{-1}}(\tau(g,h))\tau((gh)^{-1},gh)
\\ & = &
c_{gh}^{-1}(\tau(g,h)).
\end{eqnarray*}
This finishes the proof of the commutativity of the
diagram~\eqref{F_T_L_compatible}.

Next we show that $E_{\FGP{R}} \colon \id_{\FGP{R}} \to I_{\FGP{R}}
\circ I_{\FGP{R}}$ is a natural transformation of contravariant
functors of additive categories with weak $G$-action.  We have to show
that the diagram~\eqref{F_I_T_add_G-cat} commutes.  This is equivalent
to show for every finitely generated projective $R$-module $P$ the
following diagram commutes
$$
\xymatrix@!C=10em{\res_{c_g} P \ar[r]^-{E_{\FGP{R}}(\res_{c_g} P)}
  \ar[d]^-{\res_{c_g} E_{\FGP{R}}(P)} & (\res_{c_g} P)^{**}
  \ar[d]^-{t_g(P)^*}
  \\
  \res_{c_g} \left(P^{**}\right) \ar[r]^-{t_g(P^*)} & (\res_{c_g}
  P^*)^* }
$$
We start with an element $p \in P$ in the left upper corner. It is
sent under the left vertical arrow to the element given by $f \mapsto
f(p)$. The image of this element under the lower horizontal is given
by $f \mapsto c_g^{-1}(f(p))\left(w(g)\tau(g^{-1},g)\right)^{-1}$. The
image of $p \in P$ under the upper horizontal arrow is $f \mapsto
f(p)$.  The image of this element under the right vertical arrow sends
$f$ to $f \circ t_g(P)(p) =
c_g^{-1}(f(p))\left(w(g)\tau(g^{-1},g)\right)^{-1}$.

>From the naturality of the construction of the additive category with
strict $G$-action $\FGP{R}_{c,\tau} := \cals(\FGP{R})$ (see
Section~\ref{sec:Making_an_additive_categories_with_weak_(G,v)-action_strict})
we conclude that $(I_{\FGP{R}},\{t_g \mid g \in G\})$ induces a functor
of additive categories with strict $G$-action
$$I_{\FGP{R}_{c,\tau}}  \colon \FGP{R}_{c,\tau} \to \FGP{R}_{c,\tau}$$
and $E_{\FGP{R}}$ induces a  natural transformation of functors of
additive categories with strict $G$-action
$$E_{\FGP{R}_{c,\tau}} \colon   \id_{\FGP{R}} \to
I_{\FGP{R}_{c,\tau}} \circ I_{\FGP{R}_{c,\tau}}.$$
It remains to prove that condition~\eqref{E(I(A))_is_I(E(A)-1)} holds
for $(I_{\FGP{R}_{c,\tau}},E_{\FGP{R}_{c,\tau}})$.
But this follows easily from the fact that
condition~\eqref{E(I(A))_is_I(E(A)-1)} holds
for $(I_{\FGP{R}},E_{\FGP{R}})$.
\end{proof}

The additive $G$-category with involution constructed in
Lemma~\ref{lem:involution_on_FGP(R_c,w)} will be denoted in the
sequel by
\begin{eqnarray}
& \FGP{R}_{c,\tau,w}.&
\label{FGP(R)_(c,tau,w)}
\end{eqnarray}

%%%%%%%%%%%%%%%%%%%%%%%%%%%%%%%%%%%%%%%%%%%%%%%%%%%%%%%%%%%%%%%%%%%%%%%%%%%%%%%%%%%%%%%%%
%%%%%%%%%%%%% Connected groupoids an additive categories %%%%%%%%%%%%%%%%%%%%%%%%%%%%%%%%
%%%%%%%%%%%%%%%%%%%%%%%%%%%%%%%%%%%%%%%%%%%%%%%%%%%%%%%%%%%%%%%%%%%%%%%%%%%%%%%%%%%%%%%%%

\typeout{-------------------------  Section 5:  Connected groupoids an additive categories
  ----------------------}

\section{Connected groupoids and additive categories}
\label{sec:Connected_groupoids_and_additive_categories}

Groupoids are always to be understood to be small.
A groupoid is called \emph{connected} if for two objects $x$ and
$y$ there exists a morphism $f \colon x \to y$.
Let $\calg$ be a connected groupoid. Let $\addcat$ be the category of small additive categories.
Given a contravariant functor $F \colon \calg \to \addcat$, we
define a new small additive category, which we call its \emph{homotopy colimit}
(see for instance~\cite{Thomason(1979)})
\begin{eqnarray}
& \intgf{\calg}{F} &
\label{int_calg_F}
\end{eqnarray}
as follows.  An object is a pair $(x,A)$ consisting of an object $x$ in $\calg$ and
an object $A$ in $F(x)$.  A morphism in $\intgf{\calg}{F}$ from $(x,A)$ to $(y,B)$
is a formal sum
$$\sum_{f \in \mor_{\calg}(x,y)} f \cdot \phi_f$$
where $\phi_f \colon A \to F(f)(B)$ is a morphism in $F(x)$
and only finitely many coefficients $\phi_f$ are
different from zero.  The composition of a morphism $\sum_{f \in \mor_{\calg}(x,y)} f \cdot \phi_f
\colon (x,A) \to (y,B)$ and a morphism
$\sum_{g \in \mor_{\calg}(y,z)} g \cdot \phi_g
\colon (y,B) \to (z,C)$ is given by the formula
$$\sum_{h \in \mor_{\calg}(x,z)} h \cdot
\Biggl(\sum_{\substack{f \in \mor_{\calg}(x,y)\\g \in \mor_{\calg}(y,z)\\ h = g\circ f}}
F(f)(\psi_g) \circ \phi_f)\Biggr).$$
The decisive special case is
$$(g \cdot \psi ) \circ (f \cdot \phi) = (g \circ f) \cdot (F(f)(\psi) \circ \phi).$$
The $\IZ$-module structure on $\mor_{\intgf{\calg}{F}}(x,y)$ is given by
\begin{eqnarray*}
\left(\sum_{f \in \mor_{\calg}(x,y)} f \cdot \phi_f\right) + \left( \sum_{f \in \mor(\calg)} f \cdot \psi_f\right)
& = &
\sum_{f \in \mor_{\calg}(x,y)} f \cdot (\phi_f + \psi_f).
\end{eqnarray*}
A model for the sum of two objects $(x,A)$ and $(x,B)$ is $(x,A \oplus B)$
if $A \oplus B$ is a model for the sum of $A$ and $B$ in $F(x)$. Since $\calg$ is by
assumption connected, we can choose for any object
$(y,B)$ in $\intgf{\calg}{F}$ and any object $x$ in $\calg$
an isomorphism $f \colon x \to y$ and the objects $(x,F(f)(B))$ and $(y,B)$
in $\intgf{\calg}{F}$ are isomorphic. Namely $f \cdot \id_{F(f)(B)}$ is an isomorphism
$(x,F(f)(B) \xrightarrow{\cong} (y,B)$ whose inverse is $f^{-1} \cdot \id_{B} \cdot $
Hence the direct sum of two arbitrary objects $(x,A)$ and $(y,B)$ exists in
$\intgf{\calg}{F}$.

Notice that we need the connectedness of $\calg$ only to show the
existence of a direct sum. This will become important later when
we deal with non-connected groupoids.

This construction is functorial in $F$. Namely, if $S \colon F_0 \to F_1$ is a natural
transformation of contravariant functors $\calg \to \addcat$, then it induces a functor
\begin{eqnarray}
& \intgf{\calg}{S}  \colon \intgf{\calg}{F_0} \to \intgf{\calg}{F_1} &
\label{intgf(calg)(S)}
\end{eqnarray}
of additive categories as follows.
It sends an object $(x,A)$ in $\intgf{\calg}{F_0}$ to the object
$(x,S(x)(A))$ in $\intgf{\calg}{F_1}$.
A morphism $\sum_{f \in \mor_{\calg}(x,y)} f \cdot \phi_f
\colon (x,A) \to (y,B)$ is sent to the morphism
$$\sum_{f \in \mor_{\calg}(x,y)} f \cdot S(x)(\phi_f)
\colon (x,s(x)(A)) \to (y,s(y)(B)).$$ This makes sense since
$S(x)(\phi_f)$ is a morphism in $F_1(x)$ from $S(x)(A)$ to
$S(x)(F_0(f)(B)) = F_1(f)(S(y)(B))$. The decisive special case is
that $\intgf{\calg}{S}$ sends $(f \colon x \to y) \cdot \phi$ to
$(f \colon x \to y) \cdot S(x)(\phi)$. One easily checks that
$\intgf{\calg}{S}$ is compatible with the structures of additive
categories and we have
\begin{eqnarray}
\left(\intgf{\calg}{S_2}\right) \circ \left(\intgf{\calg}{S_1}\right)
& = &
\intgf{\calg}{(S_2 \circ S_1)};
\label{int_S_2_circint_S_1_is_int_S_2_circS_1}
\\
\intgf{\calg}{\id_F} & = & \id_{\intgf{\calg}{F}}.
\label{int_id_is_id}
\end{eqnarray}

The construction is also functorial in $\calg$. Namely, let $W \colon \calg_1 \to \calg_2$ be
a covariant functor of groupoids.  Then we obtain a covariant functor
\begin{eqnarray}
& W_* \colon \intgf{\calg_1}{F} \circ W  \to  \intgf{\calg_2}{F} &
\label{W_ast}
\end{eqnarray}
of additive categories as follows.
An object $(x_1,A)$ in $\intgf{\calg_1}{F \circ W}$ is sent to the object
$(W(x_1),A)$  in $\intgf{\calg_2}{F}$. A morphism
$\sum_{f \in \mor_{\calg_1}(x_1,y_1)} f \cdot \phi_f
\colon (x_1,A) \to (y_1,B)$ in  $\intgf{\calg_1}{F \circ W}$ is sent to the morphism
$$\sum_{f \in \mor_{\calg_2}(W(x_1),W(y_1))} f \cdot
\Biggl(\sum_{\substack{f_1 \in \mor_{\calg_1}(x_1,y_1)\\ W(f_1) = f}}
\phi_{f_1}\Biggr) \colon(W(x_1),A) \to (W(y_1),B)$$
in $\intgf{\calg_2}{F}$. Here the decisive special case is that
$W_*$ sends the morphism $f \cdot \phi$ to $W(f) \cdot \phi$.
One easily checks that $W_*$ is compatible with the structures of additive categories
and we have for covariant functors
$W_1 \colon \calg_1 \to \calg_2$, $W_2 \colon \calg_2 \to \calg_3$
and a contravariant functor $F \colon \calg \to \addcat$
\begin{eqnarray}
(W_2)_* \circ (W_1)_* & = & (W_2 \circ W_1)_*;
\label{(W_2)_ast_circ_(W_1)_ast_is_(W_2_circ_W_1)_ast}
\\
(\id_{\calg})_* & = & \id_{\intgf{\calg}{F}}.
\label{(id_calg)_ast_is_id_int_calg_F}
\end{eqnarray}
These two constructions are compatible. Namely, given a natural transformation
$S_1 \colon F_1 \to F_2$ of contravariant functors $\calg \to \addcat$ and a covariant
functor $W \colon \calg_1 \to \calg$, we get
\begin{eqnarray}
\left(\intgf{\calg}{S} \right) \circ W_*
& = &
W_* \circ \left(\intgf{\calg_1}{(S \circ W)} \right).
\label{compatibility_of_W_ast_and_int_S}
\end{eqnarray}

A functor $F \colon \calc_0 \to \calc_1$ of categories is called
an \emph{equivalence} if there exists a functor $F' \colon \calc_1
\to \calc_0$ with the property that $F' \circ F$ is naturally
equivalent to the identity functor $\id_{\calc_0}$ and $F \circ
F'$ is naturally equivalent to the identity functor
$\id_{\calc_1}$. A functor $F$ is a natural equivalence if and
only if it is \emph{full} and \emph{faithful}, i.e., it induces a
bijection on the isomorphism classes of objects and for any two
objects $c,d$ in  $\calc_0$ the induced map $\mor_{\calc_0}(c,d)
\to \mor_{\calc_1}(F(c),F(d))$ is bijective. If $\calc_0$ and
$\calc_1$ come with an additional structure such as of an additive
category (with involution) and $F$ is compatible with this
structure, we require that $F'$ and the two natural equivalences
$F' \circ F \simeq \id_{\calc_0}$ and $F \circ F' \simeq
\id_{\calc_1}$ are compatible with these. In this case it still
true that $F$ is an equivalence of categories with this additional
structure if and only if $F$ is full and faithful.

One easily checks
\begin{lemma}
\label{lem:(F_1)_ast_and_int_calg_S_and_equivalences}
\begin{enumerate}
\item \label{lem:(F_1)_ast_and_int_calg_S_and_equivalences:(F_1)_ast}
Let $W \colon \calg_1 \to \calg$ be an equivalence of connected groupoids.
Let $F \colon \calg \to \addcat$ be a contravariant functor.  Then
$$W_* \colon \intgf{\calg_1}{F \circ W} \to \intgf{\calg}{F}$$
is an equivalence of additive categories.

\item \label{lem:(F_1)_ast_and_int_calg_S_and_equivalences:int_calg_S}
Let $\calg$ be a connected groupoid. Let $S \colon F_1 \to F_2$ be a transformation
of contravariant functors $\calg \to \addcat$ such that for every
object $x$ in $\calg$ the functor $S(x) \colon F_0(x) \to F_1(x)$ is
an equivalence of additive categories. Then
$$\intgf{\calg}{S} \colon \intgf{\calg}{F_1} \to \intgf{\calg}{F_2}$$
is an equivalence of additive categories.
\end{enumerate}
\end{lemma}

%%%%%%%%%%%%%%%%%%%%%%%%%%%%%%%%%%%%%%%%%%%%%%%%%%%%%%%%%%%%%%%%%%%%%%%%%%%%%%%%%%%%%%%%%
%%%%%%%%%%%%%%%%%%%   From crossed product rings to additive categories %%%%%%%%%%%%%%%%%
%%%%%%%%%%%%%%%%%%%%%%%%%%%%%%%%%%%%%%%%%%%%%%%%%%%%%%%%%%%%%%%%%%%%%%%%%%%%%%%%%%%%%%%%%

\typeout{-------------------------  Section 6:  From crossed product rings to additive categories
  ----------------------}

\section{From crossed product rings to additive categories}
\label{sec:From_crossed_product_rings_to_additive_categories}

\begin{example}\label{exa:R_ast_c,tau_as_add_cat_with_G-action}
Here is our main example of a contravariant functor $\calg \to \addcat$.
Notice that a group $G$ is the same as a groupoid with one object
and hence a contravariant functor from a group $G$ to $\addcat$
is the same as an additive $G$-category what is the same as an
additive category  with strict $G$-action
(see Definition~\ref{def:additive_category_with_weak_(G,v)-action}).
Let $R$ be a ring together with maps of sets
\begin{eqnarray*}
c\colon G & \to & \aut(R), \quad g \mapsto c_g;
\\
\tau\colon G \times G & \to & R^{\times}.
\end{eqnarray*}
satisfying~\eqref{tau_and_c_group_homo}, \eqref{tau_and_c_cocykel}, \eqref{c_e-Is_id},
\eqref{tau(g,e)_is_1} and \eqref{tau(e,g)_is_1}.
We have introduced the  additive $G$-category
$\FGP{R}_{c,\tau}$ in~\eqref{FGP(R)_(c,tau)}. All the construction restrict
to the subcategory $\FGF{R} \subseteq \FGP{R}$ of finitely generated free $R$-modules
and lead to the  additive $G$-category
\begin{eqnarray}
\FGF{R}_{c,\tau} & := & \cals(\FGP{R});
\label{FGF(R)_(c,tau)}
\end{eqnarray}
\end{example}

\begin{lemma} \label{lem:int_G_R_FGF(R)_c,tau_is_R_ast_c,tauG-FGF}
Consider the data $(R,c,\tau)$ and the additive category
$\FGF{R}_{c,\tau}$ appearing in
Example~\ref{exa:R_ast_c,tau_as_add_cat_with_G-action}.
Let $\intgf{G}{\FGF{R}_{c,\tau}}$ be the additive category defined
in~\eqref{int_calg_F}. Since $G$ regarded as a groupoid has
precisely one object, we can (and will) identify the set of objects in
$\intgf{G}{\FGF{R}_{c,\tau}}$ with the set of objects in $\FGF{R}_{c,\tau}$
which consists of pairs $(M,g)$ for $M$ a finitely generated free $R$-module and $g \in G$.
Denote by
$\left(\intgf{G}{\FGF{R}_{c,\tau}}\right)_e$ the full subcategory of
$\intgf{G}{\FGF{R}_{c,\tau}}$ consisting of objects of the shape
$(M,e)$ for $e \in G$ the unit element. Denote by $R \ast G = R
\ast_{c,\tau} G$ the crossed product ring (see~\eqref{R_ast_G}).
Then

\begin{enumerate}
\item
\label{lem:int_G_R_FGF(R)_c,tau_is_R_ast_c,tauG-FGF:R_ast_G} There
is an equivalence of additive categories
$$
\alpha \colon \left(\intgf{G}{\FGF{R}_{c,\tau}}\right)_e  \to  \FGF{R \ast_{c,\tau}G};
$$

\item  \label{lem:int_G_R_FGF(R)_c,tau_is_R_ast_c,tauG-FGF:e}
The inclusion
$$\left(\intgf{G}{\FGF{R}_{c,\tau}}\right)_e \to \intgf{G}{\FGF{R}_{c,\tau}}$$
is an equivalence of additive categories.
\end{enumerate}

\end{lemma}
\begin{proof}\ref{lem:int_G_R_FGF(R)_c,tau_is_R_ast_c,tauG-FGF:R_ast_G}
An object $(M,e)$ in $\left(\intgf{G}{\FGF{R}_{c,\tau}}\right)_e$ is
sent under $\alpha$ to the finitely generated free $R \ast_{c,\tau}G$-module $R
\ast_{c,\tau}G \otimes_R M$. A morphism $\phi = \sum_{g \in G}g
\cdot \left(\phi_g \colon M \to \res_{c_g}(N)\right)$ from $(M,e)$
to $(N,e)$ is sent to the $R \ast_{c,\tau}G$-homomorphism
$$\alpha(\phi) \colon R \ast_{c,\tau}G \otimes_R M \to R \ast_{c,\tau}G \otimes_R N,
\quad u \otimes x \mapsto \sum_{g \in G} u\cdot
\tau(g^{-1},g)^{-1} \cdot g^{-1} \otimes \phi_g(x)$$ for $u \in  R
\ast_{c,\tau}G$ and $x \in M$. This is well-defined, i.e.,
compatible with the tensor relation, by the following calculation
for $r \in R$ using~\eqref{tau_and_c_group_homo} and
\eqref{c_e-Is_id}.
\begin{eqnarray*}
\lefteqn{u\cdot \tau(g^{-1},g)^{-1} \cdot g^{-1} \otimes \phi_g(rx)}
& &
\\
& = &
u\cdot \tau(g^{-1},g)^{-1} \cdot g^{-1} \otimes c_g(r) \phi_g(x)
\\
& = &
u\cdot \tau(g^{-1},g)^{-1} \cdot g^{-1} \cdot c_g(r) \otimes \phi_g(x)
\\
& = &
u\cdot \tau(g^{-1},g)^{-1} c_{g^{-1}}(c_g(r)) \cdot g^{-1}\otimes \phi_g(x)
\\
& = &
u\cdot \tau(g^{-1},g)^{-1} c_{g^{-1}}(c_g(r)) \tau(g^{-1},g)\tau(g^{-1},g)^{-1}
\cdot g^{-1}\otimes \phi_g(x)
\\
& = &
u \cdot c_{g^{-1}g}(r)\tau(g^{-1},g)^{-1} \cdot g^{-1}\otimes \phi_g(x)
\\
& = &
u \cdot c_{e}(r)\tau(g^{-1},g)^{-1} \cdot g^{-1} \otimes \phi_g(x)
\\
& = &
(u \cdot r) \tau(g^{-1},g)^{-1} \cdot g^{-1}\otimes \phi_g(x).
\end{eqnarray*}
Next we show that $\alpha$ is a covariant functor.
Obviously $\alpha(\id_{(M,e)}) = \id_{\alpha(M,e)}$.
Consider morphisms $\phi = \sum_{g \in G} g \cdot \phi_g  \colon (M,e) \to (N,e)$  and
$\psi = \sum_{g \in G} g \cdot \psi_g\colon (N,e) \to (P,e)$ in
$\left(\intgf{G}{\FGF{R}_{c,\tau}}\right)_e$.
A direct computation shows for $u \in  R \ast_{c,\tau}G$ and $x \in M$
\begin{eqnarray*}
\lefteqn{\alpha(\psi)\left(\alpha(\phi)(u \otimes x)\right)}
& &
\\
& = &
\alpha(\psi)\left(\sum_{k \in G} u\cdot \tau(k^{-1},k)^{-1} \cdot k^{-1} \otimes \phi_k(x)\right)
\\
& = &
\sum_{h \in G}\sum_{k \in G} u\cdot \tau(k^{-1},k)^{-1} \cdot k^{-1}\cdot
\tau(h^{-1},h)^{-1} \cdot h^{-1}
\otimes \psi_h\circ \phi_k(x)
\\
& = &
\sum_{h,k \in G} u\cdot \tau(k^{-1},k)^{-1} c_{k^{-1}}(\tau(h^{-1},h)^{-1}) \cdot k^{-1}\cdot
h^{-1}  \otimes \psi_h\circ \phi_k(x)
\\
& = &
\sum_{h,k \in G} u\cdot \tau(k^{-1},k)^{-1} c_{k^{-1}}(\tau(h^{-1},h)^{-1}) \tau(k^{-1},h^{-1})
\cdot  (hk)^{-1}
\otimes \psi_h\circ \phi_k(x)
\end{eqnarray*}
and
\begin{eqnarray*}
\lefteqn{\alpha(\psi \circ \phi)(u \otimes x)}
& &
\\
& = &
\sum_{g \in G} u\cdot \tau(g^{-1},g)^{-1} \cdot g^{-1} \otimes (\psi \circ \phi)_g(x)
\\
& = &
\sum_{g \in G} u\cdot \tau(g^{-1},g)^{-1} \cdot g^{-1} \otimes
\left(\sum_{\substack{hk \in G,\\hk = g}} r_k(\psi_{h}) \circ \phi_{k}(x)\right)
\\
& = &
\sum_{g \in G} u\cdot \tau(g^{-1},g)^{-1} \cdot g^{-1} \otimes
\left(\sum_{\substack{hk \in G,\\hk = g}} \tau(h,k)^{-1} \psi_{h} \circ \phi_{k}(x)\right)
\\
& = &
\sum_{g \in G} \sum_{\substack{hk \in G,\\hk = g}} u\cdot \tau(g^{-1},g)^{-1} \cdot
g^{-1} \cdot \tau(h,k)^{-1}
\otimes \psi_{h} \circ \phi_{k}(x)
\\
& = &
\sum_{g \in G} \sum_{\substack{hk \in G,\\hk = g}} u\cdot \tau(g^{-1},g)^{-1}c_{g^{-1}}(\tau(h,k)^{-1})
\cdot g^{-1} \otimes \psi_{h} \circ \phi_{k}(x)
\\
& = &
\sum_{hk \in G} u\cdot \tau((hk)^{-1},hk)^{-1}c_{(hk)^{-1}}(\tau(h,k)^{-1})
\cdot (hk)^{-1} \otimes \psi_{h} \circ \phi_{k}(x).
\end{eqnarray*}
Hence it remains to show for $h,k \in G$
\begin{eqnarray*}
\tau(k^{-1},k)^{-1} c_{k^{-1}}(\tau(h^{-1},h)^{-1}) \tau(k^{-1},h^{-1})
& = &
\tau((hk)^{-1},hk)^{-1}c_{(hk)^{-1}}(\tau(h,k)^{-1}),
\end{eqnarray*}
or, equivalently,
\begin{eqnarray*}
\tau(k^{-1},h^{-1})c_{(hk)^{-1}}(\tau(h,k))\tau((hk)^{-1},hk)
& = &
c_{k^{-1}}(\tau(h^{-1},h)) \tau(k^{-1},k).
\end{eqnarray*}
Since~\eqref{tau_and_c_cocykel} yields
$$\tau((hk)^{-1},h)\tau(k^{-1},k) = c_{(hk)^{-1}}(\tau(h,k)\tau((hk)^{-1},hk),$$
it suffices to show
\begin{eqnarray*}
\tau(k^{-1},h^{-1})\tau((hk)^{-1},h)
& = &
c_{k^{-1}}(\tau(h^{-1},h)).
\end{eqnarray*}
But this follows from~\eqref{tau_and_c_cocykel}
and~\eqref{tau(e,g)_is_1}. This finishes the proof that $\alpha$
is a covariant functor. Obviously it is compatible with the
structures of an additive category. One easily checks that
$\alpha$ induces a bijection between the isomorphism classes of
objects. In order to show that $\alpha$ is a weak equivalence, we
have to show for two objects $(M,e)$ and $(N,e)$ that $\alpha$
induces a bijection
$$\mor_{\left(\intgf{G}{\FGF{R}_{c,\tau}}\right)_e}((M,e),(N,e)) \xrightarrow{\cong}
\hom_{R \ast_{c,\tau}G}(R \ast_{c,\tau}\otimes_R M, R \ast_{c,\tau}\otimes_R N).$$
Since $\alpha$ is compatible with the structures of an additive category, it suffices
to check this in the special case $M = N = R$, where it is obvious.
\\[1mm]\ref{lem:int_G_R_FGF(R)_c,tau_is_R_ast_c,tauG-FGF:e}
An object
of the shape $(M,g)$ in $\intgf{G}{\FGF{R}_{c,\tau}}$ is isomorphic to the object
$(M,e)$, namely an isomorphism $(M,g) \xrightarrow{\cong} (M,e)$ in
 $\intgf{G}{\FGF{R}_{c,\tau}}$ is given by $g \cdot \id_{\res_ {c_g}(M)}$.
\end{proof}

%%%%%%%%%%%%% Connected groupoids an additive categories with involution %%%%%%%%%%%%%%%%
%%%%%%%%%%%%%%%%%%%%%%%%%%%%%%%%%%%%%%%%%%%%%%%%%%%%%%%%%%%%%%%%%%%%%%%%%%%%%%%%%%%%%%%%%

\typeout{-  Section 7:  Connected groupoids and additive categories with involutions ---}

\section{Connected groupoids and additive categories with involutions}
\label{sec:Connected_groupoids_and_additive_categories_with_involutions}

Next we want to enrich the constructions of
Section~\ref{sec:Connected_groupoids_and_additive_categories}
to additive categories with involutions.
Let $\addcatinv$ be the category of  additive categories with involution.
Given a contravariant functor $(F,T) \colon \calg \to \addcatinv$, we want to
define on the additive category $\intgf{\calg}{F}$ the structure
of an additive category with involution.
Here the pair $(F,T)$ means that we assign to every object $x$ in $\calg$
an additive category with involution $F(x)$ and for every morphism $f \colon x \to y$ in $\calg$
we have a functor of additive categories with involution
$(F(f),T(f)) \colon F(y) \to F(x)$.

Next we construct for a  functor $\calg \to \addcatinv$
an involution of additive categories
\begin{eqnarray}
& (I_{\intgf{\calg}{F}},E_{\intgf{\calg}{F}}) &
\label{int_calg_F,E_int_calg_F}
\end{eqnarray}
on the additive category
$\intgf{\calg}{F}$ which we have introduced in~\eqref{int_calg_F}.
On objects we put
$$I_{\intgf{\calg}{F}}(x,A) := (x,I_{\calg}(A)) = (x,A^*).$$
Let $\phi = \sum_{f \in \mor_{\calg}(x,y)}f \cdot \phi_f \colon (x,A) \to (y,B)$
be a morphism  in $\intgf{\calg}{F}$. Define $I_{\intgf{\calg}{F}}(\phi)\colon B^* \to A^*$
to be the morphism $\phi^* =
\sum_{f \in \mor_{\calg}(y,x)} f \cdot (\phi^*)_f \colon (y,B^*) \to (x,A^*)$
in $\intgf{\calg}{F}$ whose component for $f \in \mor_{\calg}(y,x)$ is given by the
composite
\begin{multline*}(\phi^*)_f \colon B^* = F(f)\left(F(f^{-1})(B^*)\right)
\xrightarrow{F(f)(T(f^{-1})(B))} F(f)\left(F(f^{-1})(B)^*\right)
\\
\xrightarrow{F(f)\left((\phi_{f^{-1}})^*\right)} F(f)(A^*).
\end{multline*}
Next we show that $I_{\intgf{\calg}{F}}$ is a contravariant functor.
Obviously $I_{\intgf{\calg}{F}}$ sends the identity $\id_{A}$ to $\id_{I_{\intgf{\calg}{F}}(A)}$.
We have to show
$$I_{\intgf{\calg}{F}}(\psi \circ \phi) = I_{\intgf{\calg}{F}}(\phi) \circ I_{\intgf{\calg}{F}}(\psi)$$
for morphisms $\phi = \sum_{h \in \mor_{\calg}(x,y)} h \cdot \phi_h \colon (x,A) \to (y,B)$
and $\psi \colon \sum_{k \in \mor_{\calg}(y,z)} k \cdot \psi_k \colon (y,B \to (z,C)$,
or in short notation
$(\psi \circ \phi)^* = \phi^* \circ \psi^*$.

By definition
$(\phi^* \circ \psi^*) = \sum_{g \in \mor_{\calg}(z,x)} g \cdot (\phi^* \circ \psi^*)_g$
for
$$(\phi^* \circ \psi^*)_{g} :=
\sum_{\substack{k \in \mor_{\calg}(z,y),\\h \in \mor_{\calg}(y,x),\\hk = g}}
F(k)((\phi^*)_{h}) \circ (\psi^*)_{k}.$$
By definition
\begin{multline*}
(\psi^*)_{k} \colon C^* = F(k)(F(k^{-1})(C^*)) \xrightarrow{F(k)(T(k^{-1})(C))} F(k)(F(k^{-1})(C)^*)
\\
\xrightarrow{F(k)((\psi_{k^{-1}})^*)} F(k)(B^*)
\end{multline*}
and
\begin{multline*}
(\phi^*)_{h} \colon B^* = F(h)(F(h^{-1})(B^*)) \xrightarrow{F(h)(T(h^{-1})(B))}
F(h)(F(h^{-1})(B)^*)
\\
\xrightarrow{F(h)((\phi_{h^{-1}})^*)} F(h)(A^*).
\end{multline*}
Hence the component $(\phi^* \circ \psi^*)_{g}$ of $(\phi^* \circ \psi^*)$ at $g \colon z \to x$
is given by the sum
of morphisms from $C^*$ to $F(g)(A^*)$
\begin{multline*}
\sum_{\substack{k \in \mor_{\calg}(z,y),\\h \in \mor_{\calg}(y,x),\\hk = g}}
F(k)\left(F(h)((\phi_{h^{-1}})^*)\right) \circ F(k)\left(F(h)(T(h^{-1})(B))\right)
\\
\circ F(k)((\psi_{k^{-1}})^*) \circ F(k)(T(k^{-1})(C)).
\end{multline*}
The component of $(\psi \circ \phi)^*_{g}$ of $(\psi \circ \phi)^*$ at
$g\colon z \to x$ is given by
\begin{multline*}
C^* = F(g)\left(F(g^{-1})(C^*)\right)
\xrightarrow{F(g)(T(g^{-1})(C))} F(g)\left(F(g^{-1})(C)^*\right)
\\
\xrightarrow{F(g)\left(((\psi \circ \phi)_{g^{-1}})^*\right)} F(g)(A^*).
\end{multline*}
Since for $g \colon z \to x$ we have
$$(\psi \circ \phi)_{g^{-1}} =
\sum_{\substack{h \in \mor_{\calg}(y,z),\\k \in \mor_{\calg}(x,y),\\hk = g^{-1}}}
F(k)(\psi_{h}) \circ \phi_{k},$$
the component of $(\psi \circ \phi)^*_{g}$ of $(\psi \circ \phi)^*$ at $g\colon z \to x$
is given by the sum of morphisms $C^*$ to $F(g)(A^*)$
$$\sum_{\substack{h \in \mor_{\calg}(y,z),\\k \in \mor_{\calg}(x,y),\\hk = g^{-1}}}
 F(g)\left((\phi_k)^*\right) \circ F(g)\left(F(k)(\psi_h)^*\right)
\circ F(g)(T(g^{-1})(C)).$$
By changing the indexing by replacing $h$ with $k^{-1}$ and $k$ by $h^{-1}$, this transforms to
$$\sum_{\substack{k \in \mor_{\calg}(z,y),\\h \in \mor_{\calg}(y,x),\\hk = g}}
F(g)\left((\phi_{h^{-1}})^*\right) \circ F(g)\left(F(h^{-1})(\psi_{k^{-1}})^*\right)
\circ F(g)(T(g^{-1})(C)).$$

Hence we have to show for every $k \colon z \to y$ and $h \colon y \to x$ with $hk = g$
that the two composites
$$
F(k)\left(F(h)((\phi_{h^{-1}})^*)\right) \circ F(k)\left(F(h)(T(h^{-1})(B))\right)
\circ F(k)((\psi_{k^{-1}})^*) \circ F(k)(T(k^{-1})(C))$$
and
$$(F(g)((\phi_{h^{-1}})^*) \circ F(g)\left(F(h^{-1})(\psi_{k^{-1}})^*\right)
\circ F(g)(T(g^{-1})(C))$$
agree.  We compute for the first one
\begin{eqnarray*}
\lefteqn{F(k)\left(F(h)((\phi_{h^{-1}})^*)\right) \circ F(k)\left(F(h)(T(h^{-1})(B))\right)
\circ F(k)((\psi_{k^{-1}})^*) \circ F(k)(T(k^{-1})(C))} & &
\\
& = &
(F(g)((\phi_{h^{-1}})^*) \circ F(g)(T(h^{-1})(B)) \circ
F(g)\left(F(h^{-1})((\psi_{k^{-1}})^*)\right)
\\
& & \hspace{70mm} \circ   F(g)\left(F(h^{-1})(T(k^{-1})(C))\right).
\end{eqnarray*}
Hence it remains to show that the composites
$$F(g^{-1})(C^*) \xrightarrow{T(g^{-1})(C)} F(g^{-1})(C)^*
\xrightarrow{F(h^{-1})(\psi_{k^{-1}})^*} F(h^{-1})(B)^*
$$
and
\begin{multline*}
F(g^{-1})(C^*) = F(h^{-1})\left(F(k^{-1})(C^*)\right)
 \xrightarrow{F(h^{-1})(T(k^{-1})(C))}
F(h^{-1})(F(k^{-1})(C)^*)
\\
\xrightarrow{F(h^{-1})((\psi_{k^{-1}})^*)} F(h^{-1})(B^*)
\xrightarrow{T(h^{-1})(B)} F(h^{-1})(B)^*
\end{multline*}
agree.  The second one agrees with the composite
\begin{multline*}
F(g^{-1})(C^*)= F(h^{-1})\left(F(k^{-1})(C^*)\right)
 \xrightarrow{F(h^{-1})(T(k^{-1})(C))}
F(h^{-1})(F(k^{-1})(C)^*)
\\
\xrightarrow{T(h^{-1})(F(k^{-1})(C))} F(h^{-1})(F(k^{-1})(C))^*
\xrightarrow{F(h^{-1})(\psi_{k^{-1}})^*} F(h^{-1})(B)^*
\end{multline*}
since $T(h^{-1})$ is a natural transformation
$F(h^{-1}) \circ I_{F(y)} \to I_{F(x)} \circ F(h^{-1})$.
Since
$$(F(h^{-1}),T(h^{-1})) \circ (F(k^{-1}),T(k^{-1})) = (F(k^{-1}h^{-1}),T(k^{-1}h^{-1}))
= (F(g^{-1}),T(g^{-1}))$$
the map $T(g^{-1})(C)$ can be written as the composite
\begin{multline*}
T(g^{-1})(C) \colon F(g^{-1})(C^*) = F(h^{-1})\left(F(k^{-1})(C^*)\right)
\\
\xrightarrow{F(h^{-1})\left(T(k^{-1})(C)\right)} F(h^{-1})\left(F(k^{-1})(C)^*\right)
\\
\xrightarrow{T(h^{-1})\left(F(k^{-1})(C)\right)} F(h^{-1})\left(F(k^{-1})(C)\right)^* = F(g)(C)^*.
\end{multline*}
This finishes the proof that $I_{\intgf{\calg}{F}}$ is a contravariant functor.

The natural equivalence
$$E_{\intgf{\calg}{F}} \colon \id_{\intgf{\calg}{F}} \to
I_{\intgf{\calg}{F}} \circ I_{\intgf{\calg}{F}}$$
assigns to an object $(x,A)$ in $\intgf{\calg}{F}$ the isomorphism
$$\id_x \cdot \left(E_{\calg}(A) \colon A  \xrightarrow{\cong} A^{**}\right) \colon
(x,A) \to (x,A^{**}).$$

We have to check that $E_{\intgf{\calg}{F}}$ is a natural equivalence. Consider a morphism
$\phi = \sum_{f \in \mor_{\calg}(x,y)} f \cdot
\left(\phi_f \colon (x,A) \to (y,B)\right)$ in $\intgf{\calg}{F}$. Then
$\left(I_{\intgf{\calg}{F}} \circ I_{\intgf{\calg}{F}}\right)(\phi)$ has as
the component for $f \colon x \to y$ the composite
\begin{multline*} A^{**} = F(f)\left(F(f^{-1})(A^{**})\right)
\xrightarrow{F(f)\left((T(f^{-1})(A^*)\right)} F(f)\left(F(f^{-1})(A^*)^*\right)
\\
\xrightarrow{F(f)\left(F(f^{-1})((\phi_{f})^*)^*\right)}
F(f)\left(F(f^{-1})(F(f)(B)^*)^*\right)
\\
\xrightarrow{F(f)\left(F(f^{-1})(T(f)(B))^*\right)}
F(f)\left(F(f^{-1})(F(f)(B^*))^*\right) = F(f)(B^{**}).
\end{multline*}

Hence
$\left(I_{\intgf{\calg}{F}} \circ I_{\intgf{\calg}{F}}\right)(\phi) \circ E_{\intgf{\calg}{F}}(x,A)$
has as component for $f\colon x \to y$ the composite
\begin{multline*}
A \xrightarrow{E_{\cala}(A)} A^{**} = F(f)\left(F(f^{-1})(A^{**})\right)
\xrightarrow{F(f)\left((T(f^{-1})(A^*)\right)} F(f)\left(F(f^{-1})(A^*)^*\right)
\\
\xrightarrow{F(f)\left(F(f^{-1})(\phi_{f}^*)^*\right)}
F(f)\left(F(f^{-1})(F(f)(B)^*)^*\right)
\\
\xrightarrow{F(f)\left(F(f^{-1})(T(f)(B))^*\right)}
F(f)\left(F(f^{-1})(F(f)(B^*))^*\right) = F(f)(B^{**}).
\end{multline*}
The component of $E_{\intgf{\calg}{F}}(y,B) \circ \phi$ at $f \colon x \to y
$ is the composite
$$A \xrightarrow{\phi_{f}} F(f)(B)
\xrightarrow{F(f)(E_{\cala}(B))} F(f)(B^{**}).$$
It remains to show that these two morphisms $A \to
F(f)(B^{**})$ agree. The following two diagrams commute since
$E_{\cala}$ and $T(f^{-1})$ are natural transformations
$$
\xymatrix@!C=14em{A \ar[r]^-{E_{\cala}(A)} \ar[d]^-{\phi_{f}}
&
A^{**} = F(f)\left(F(f^{-1})(A^{**})\right) \ar[d]^-{\phi_{f}^{**} =
F(f)\left(F(f^{-1})(\phi_{f}^{**}))\right)}
\\
F(f)(B) \ar[r]^-{E_{\cala}(F(f)(B))}
&
F(f)(B)^{**} = F(f)\left(F(f^{-1})(F(f)(B)^{**})\right)
}
$$
and
$$
\xymatrix@!C=20em{
F(f)\left(F(f^{-1})(A^{**})\right) \ar[r]^-{F(f)(T(f^{-1})(A^*))}
\ar[d]^-{F(f)\left(F(f^{-1})(\phi_{f}^{**}))\right)}
&
F(f)\left(F(f^{-1})(A^*)^*\right)  \ar[d]^-{F(f)\left(F(f^{-1})(\phi_{f}^*)^*\right)}
\\
F(f)\left(F(f^{-1})(F(f)(B)^{**})\right)   \ar[r]^-{F(f)\left(T(f^{-1})(F(f)(B)^*)\right)}
&
F(f)\left(F(f^{-1})(F(f)(B)^*)^*\right).
}
$$
Hence we have to show that
$$F(f)(B) \xrightarrow{F(f)(E_{\cala}(B))} F(f)(B^{**})$$
agrees with the composite
\begin{multline*}
F(f)(B)
\xrightarrow{E_{\cala}(F(f)(B))}
F(f)(B)^{**} = F(f)\left(F(f^{-1})(F(f)(B)^{**})\right)
\\
\xrightarrow{{F(f)\left(T(f^{-1})(F(f)(B)^*)\right)}}
F(f)\left(F(f^{-1})(F(f)(B)^*)^*\right)
\\
\xrightarrow{F(f)\left(F(f^{-1})(T(f)(B))^*\right)}
F(f)\left(F(f^{-1})(F(f)(B^*))^*\right) = F(f)(B^{**}).
\end{multline*}
(Notice that $\phi$ is not involved anymore.)
The following diagram commutes by the axioms
(see \eqref{F(A)_F(A)aastast_F(Aastast)_F(Aast)ast})
$$
\xymatrix@!C=17em{
F(f)(B) \ar[r]^-{E_{\cala}(F(f)(B))}
\ar[d]^-{F(f)(E_{\cala}(B))}
&
F(f)(B)^{**} \ar[d]^-{T(f)(B)^*}
\\
F(f)(B^{**}) \ar[r]^-{T(f)(B^*)}
&
F(f)(B^*)^*
}
$$
Hence it remains to show the commutativity of the following diagram
(which does not involve $\phi$ and $E_{\cala}$ anymore).
$$
\xymatrix@!C=22em{
F(f)(B)^{**} = F(f)\left(F(f^{-1})(F(f)(B)^{*})\right)^*
\ar[r]^-{F(f)\left(T(f^{-1})(F(f)(B)^*)\right)}
\ar[d]^-{T(f)(B)^*}
&
F(f)\left(F(f^{-1})(F(f)(B)^*)^*\right)
\ar[d]_-{F(f)\left(F(f^{-1})(T(f)(B))^*\right)}
\\
F(f)(B^*)^*
&
F(f)(B^{**}) = F(f)\left(F(f^{-1})(F(f)(B^*))^*\right)
\ar[l]^-{T(f)(B^*)}
}
$$
Since $(F(f),T(f)) \circ (F(f^{-1}),T(f^{-1})) = \id$, we have
$$T(f)\left(F(f^{-1})(F(f)(B)^*)\right)~\circ~F(f)\left(T(f^{-1})(F(f)(B)^*)\right) = \id.$$
Hence it suffices to prove the commutativity of the following diagram
$$\xymatrix@!C=22em{
F(f)(B)^{**} = F(f)\left(F(f^{-1})(F(f)(B)^*)\right)^*
\ar[d]^-{T(f)(B)^*}
&
F(f)\left(F(f^{-1})(F(f)(B)^*)^*\right)
\ar[d]_-{F(f)\left(F(f^{-1})(T(f)(B))^*\right)}
\ar[l]_-{T(f)\left(F(f^{-1})(F(f)(B)^*)\right)}
\\
F(f)(B^*)^*
&
F(f)(B^{**}) = F(f)\left(F(f^{-1})(F(f)(B^*))^*\right).
\ar[l]^-{T(f)(B^*)}
}
$$
This follows because this diagram is obtained by applying the
natural transformation $T(f)$ to the morphism
$$F(f^{-1})(F(f)(B^*)) \xrightarrow{F(f^{-1})(T(f)(B))} F(f^{-1})(F(f)(B)^*).$$

The condition~\eqref{E(I(A))_is_I(E(A)-1)} is satisfied for
$(I_{\intgf{\calg}{F}},E_{\intgf{\calg}{F}})$ since it holds for
$(I_{\cala},E_{\cala})$.

We will denote the resulting additive category $\intgf{\cala}{F}$ with
involution $(I_{\intgf{\calg}{F}},E_{\intgf{\calg}{F}})$ by
\begin{eqnarray}
& \intgf{\calg}{(F,T)}. &
\label{int_calg_F,T}
\end{eqnarray}

Let $(F_0,T_0)$ and $(F_1,T_1)$ be two contravariant functors $\calg
\to \addcatinv$. Let $(S,U) \colon (F_0,T_0) \to (F_1,T_1)$ be a
natural transformation of such functors. This means that we for each
object $x$ in $\calg$ we have an equivalence $(S(x),U(x)) \colon
F_0(x) \to F_1(y)$ of additive categories with involution such that
for all $f \colon x \to y$ in $\calg$ the following diagram of
functors of additive categories with involution commutes
\begin{eqnarray}
&
\xymatrix@!C=8em{
F_0(y) \ar[r]^-{(S(y),U(y))} \ar[d]^-{(F_0(f),T_0(f))}
&
F_1(y) \ar[d]^-{(F_1(f),T_1(f))}
\\
F_0(x) \ar[r]^-{(S(x),U(x))}
&
F_1(x)
}
&
\label{(F_1(f),T_1(f))_circ_(S(y),U(y))_is_(S(x),U(x))_circ_(F_0(f),T_0(f))}
\end{eqnarray}
Then both $\intgf{\calg}{(F_0,T_0)}$ and
$\intgf{\calg}{(F_1,T_1)}$ are additive categories with involutions. The functor
of additive categories $\intgf{\calg}{S} \colon \intgf{\calg}{F_0} \to \intgf{\calg}{F_1}$
defined in~\eqref{intgf(calg)(S)} extends to a functor of additive categories with involution
\begin{eqnarray}
& \intgf{\calg}{(S,U)}  \colon \intgf{\calg}{(F_0,T_0)} \to \intgf{\calg}{(F_1,T_1)}
\label{intgf(calg)(S,U)}
\end{eqnarray}
as follows. We have to specify a natural equivalence
$$
\widehat{U} \colon \left(\intgf{\calg}{S}\right)  \circ I_{\intgf{\calg}{(F_0,T_0)}}
\to I_{\intgf{\calg}{(F_1,T_1)}}\circ \intgf{\calg}{S}.
$$
For an object $(x,A)$ in $\intgf{\calg}{F_0}$ the isomorphism
$$\widehat{U}(x,A) \colon \left(\intgf{\calg}{S}\right)  \circ I_{\intgf{\calg}{(F_0,T_0)}}(x,A)
\to I_{\intgf{\calg}{(F_1,T_1)}}\circ \intgf{\calg}{S}(x,A)$$ is given by the isomorphism
$$\id_x \cdot U(x)(A) \colon (x,S(x)(A^*)) \to (x,S(x)(A)^*)$$
in $\intgf{\calg}{F_1}$. Next we check that $\widehat{U}$ is a natural equivalence.

Let $\sum_{f \in \mor_{\calg}(x,y)} f \cdot \phi_f \colon (x,A)
\to (y,B)$ be a morphism in $\intgf{\calg}{F_0}$, where by
definition $\phi_f \colon A \to F(f)(B)$ is a morphism in the
additive category $F_0(x)$. We have to show the commutativity of
the following diagram  in the additive category $\int_{\calg} F_1$
$$\xymatrix@!C=18em{
 \left(\intgf{\calg}{S}\right)  \circ I_{\intgf{\calg}{(F_0,T_0)}}(y,B)
 \ar[r]^-{\left(\intgf{\calg}{S}\right)  \circ I_{\intgf{\calg}{(F_0,T_0)}}(\phi)}
 \ar[d]^-{\widehat{U}(y,B)}
&
\left(\intgf{\calg}{S}\right)  \circ I_{\intgf{\calg}{(F_0,T_0)}}(x,A)
\ar[d]^-{\widehat{U}(x,A)}
\\
 I_{\intgf{\calg}{(F_1,T_1)}}\circ \intgf{\calg}{S}(y,B)
 \ar[r]^-{I_{\intgf{\calg}{(F_1,T_1)}}\circ \intgf{\calg}{S}(\phi)}
&
I_{\intgf{\calg}{(F_1,T_1)}}\circ \intgf{\calg}{S}(x,A)
}
$$
The morphism $I_{\intgf{\calg}{(F_0,T_0)}}(\phi)$ in $\intgf{\calg}{F_0}$ is given by
$$\phi^* = \sum_{f \in \mor_{\calg}(x,y)} f \cdot (\phi^*)_f \colon (y,B^*) \to (x,A^*),$$
where the component $(\phi^*)_f$ is the composite
\begin{multline*}
(\phi^*)_f \colon  B^* = F_0(f)\left(F_0(f^{-1})(B^*)\right)
\xrightarrow{F_0(f)\left(T_0(f^{-1})(B)\right)}
F_0(f)\left(F_0(f^{-1})(B)^*\right)
\\
\xrightarrow{F_0(f)\left((\phi_{f{-1}})^*\right)}
F_0(f)\left(A^*\right).
\end{multline*}
The morphism $\left(\intgf{\calg}{S}\right)
\circ I_{\intgf{\calg}{(F_0,T_0)}}(\phi) \colon (y,B^*) \to (x,A^*)$
in $\intgf{\calg}{F_1}$ is given by
$\sum_{f \in \mor_{\calg}(x,y)} f \cdot \psi_f \colon (y,B^*) \to (x,A^*)$,
where $\psi_f$ is the composite
\begin{multline*}
\psi_f \colon  S(y)(B^*) = S(y)\left(F_0(f)\left(F_0(f^{-1})(B^*)\right)\right)
\\
\xrightarrow{S(y)\left(F_0(f)\left(T_0(f^{-1})(B)\right)\right)}
S(y)\left(F_0(f)\left(F_0(f^{-1})(B)^*\right)\right)
\\
\xrightarrow{S(y)\left(F_0(f)\left((\phi_{f{-1}})^*\right)\right)}
S(y)\left(F_0(f)\left(A^*\right)\right) = F_1(f)\left(S(x)\left(A^*\right)\right).
\end{multline*}
Hence the morphism
$\widehat{U}(x,A) \circ \left(\intgf{\calg}{S}\right)  \circ I_{\intgf{\calg}{(F_0,T_0)}}(\phi)
\colon (y,B^*) \to (x,A^*)$
in $\intgf{\calg}{F_1}$ is given by
$\sum_{f \in \mor_{\calg}(x,y)} f \cdot \mu_f \colon (y,B^*) \to (x,A^*)$,
where $\mu_f$ is the composite in $F_1(y)$
\begin{multline*}
\mu_f \colon  S(y)(B^*) = S(y)\left(F_0(f)\left(F_0(f^{-1})(B^*)\right)\right)
\\
\xrightarrow{S(y)\left(F_0(f)\left(T_0(f^{-1})(B)\right)\right)}
S(y)\left(F_0(f)\left(F_0(f^{-1})(B)^*\right)\right)
\\
\xrightarrow{S(y)\left(F_0(f)\left((\phi_{f{-1}})^*\right)\right)}
S(y)\left(F_0(f)\left(A^*\right)\right) = F_1(f)\left(S(x)\left(A^*\right)\right)
\\
\xrightarrow{F_1(f)\left(U(x)(A)\right)} F_1(f)\left(S(x)(A)^*\right).
\end{multline*}
The morphism $\intgf{\calg}{S}(\phi) \colon (x,S(x)(A)) \to (y,S(y)(B))$ in
$\int_{\calg} F_1$ is given by
$$\sum_{f \in \mor_{\calg}(x,y)} f \cdot \left(S(x)(\phi_f) \colon S(x)(A)
\to S(x)(F_0(f)(B) = F_1(f)(S(y)(B)\right).$$
The morphism $I_{\intgf{\calg}{(F_1,T_1)}}\circ \intgf{\calg}{S}(\phi) \colon
(y,S(y)(B)^*) \to (x,S(x)(A)^*)$ in $\int_{\calg} F_1$
is given by $\sum_{f \in \mor_{\calg}(y,x)} f \cdot \nu_f$,
where $\nu_f$ is the composite in $F_1(y)$.
\begin{multline*}
\nu_f \colon  S(y)(B)^* = F_1(f)\left(F_1(f^{-1})\left(S(y)(B)^*\right)\right)
\\
\xrightarrow{F_1(f)\left(T_1(f^{-1})\left(S(y)(B)\right)\right)}
F_1(f)\left(F_1(f^{-1})\left(S(y)(B)\right)^*\right) =
F_1(f)\left(S(x)\left(F_0(f^{-1})(B)\right)^*\right)
\\
\xrightarrow{F_1(f)\left((S(x)(\phi_{f^{-1}}))^*\right)}
F_1(f)\left(S(x)(A)^*\right).
\end{multline*}
The morphism $I_{\intgf{\calg}{(F_1,T_1)}}\circ \intgf{\calg}{S}(\phi) \circ \widehat{U}(y,B)
\colon (y,S(y)(B^*)) \to (x,S(x)(A)^*)$
in $\int_{\calg} F_1$
is given by $\sum_{f \in \mor_{\calg}(y,x)} f \cdot \omega_f$,
where $\omega_f$ is the composite in $F_1(y)$.
\begin{multline*}
\omega _f \colon  S(y)(B^*) \xrightarrow{U(y)(B)} S(y)(B)^*
= F_1(f)\left(F_1(f^{-1})\left(S(y)(B)^*\right)\right)
\\
\xrightarrow{F_1(f)\left(T_1(f^{-1})\left(S(y)(B)\right)\right)}
F_1(f)\left(F_1(f^{-1})\left(S(y)(B)\right)^*\right)
= F_1(f)\left(S(x)\left(F_0(f^{-1})(B)\right)^*\right)
\\
\xrightarrow{F_1(f)\left((S(x)(\phi_{f^{-1}}))^*\right)}
F_1(f)\left(S(x)(A)^*\right).
\end{multline*}
Hence we have to show for all $f \colon y \to x$ in $\mor_{\calg}(y,x)$
that the two composites in $F_1(y)$
\begin{multline*}
S(y)(B^*) = S(y)\left(F_0(f)\left(F_0(f^{-1})(B^*)\right)\right)
\\
\xrightarrow{S(y)\left(F_0(f)\left(T_0(f^{-1})(B)\right)\right)}
S(y)\left(F_0(f)\left(F_0(f^{-1})(B)^*\right)\right)
\\
\xrightarrow{S(y)\left(F_0(f)\left((\phi_{f{-1}})^*\right)\right)}
S(y)\left(F_0(f)\left(A^*\right)\right) = F_1(f)\left(S(x)\left(A^*\right)\right)
\\
\xrightarrow{F_1(f)\left(U(x)(A)\right)} F_1(f)\left(S(x)(A)^*\right)
\end{multline*}
and
\begin{multline*}
S(y)(B^*) \xrightarrow{U(y)(B)} S(y)(B)^* = F_1(f)\left(F_1(f^{-1})\left(S(y)(B)^*\right)\right)
\\
\xrightarrow{F_1(f)\left(T_1(f^{-1})\left(S(y)(B)\right)\right)}
F_1(f)\left(F_1(f^{-1})\left(S(y)(B)\right)^*\right)
= F_1(f)\left(S(x)\left(F_0(f^{-1})(B)\right)^*\right)
\\
\xrightarrow{F_1(f)\left((S(x)(\phi_{f^{-1}}))^*\right)}
F_1(f)\left(S(x)(A)^*\right).
\end{multline*}
agree. Since $S$ is a natural transformation from $F_0 \to F_1$,
the first composite can be rewritten
as the composite
\begin{multline*}
S(y)(B^*) = F_1(f)\left(S(x)\left(F_0(f^{-1})(B^*)\right)\right)
\\
\xrightarrow{F_1(f)\left(S(x)\left(T_0(f^{-1})(B)\right)\right)}
F_1(f)\left(S(x)\left(F_0(f^{-1})(B)^*\right)\right)
\\
\xrightarrow{F_1(f)\left(S(x)\left((\phi_{f{-1}})^*\right)\right)}
F_1(f)\left(S(x)\left(A^*\right)\right)
\\
\xrightarrow{F_1(f)\left(U(x)(A)\right)} F_1(f)\left(S(x)(A)^*\right).
\end{multline*}
Since $U(x)$ is a natural transformation from $S(x) \circ I_{F_0(x)}$ to $I_{F_1(x)} \circ S(x)$,
this agrees with the composite
\begin{multline*}
S(y)(B^*) = F_1(f)\left(S(x)\left(F_0(f^{-1})(B^*)\right)\right)
\\
\xrightarrow{F_1(f)\left(S(x)\left(T_0(f^{-1})(B)\right)\right)}
F_1(f)\left(S(x)\left(F_0(f^{-1})(B)^*\right)\right)
\\
\xrightarrow{F_1(f)\left(U(x)\left(F_0(f^{-1})(B)\right)\right)}
F_1(f)\left(S(x)\left(F_0(f^{-1})(B)\right)^*\right)
\\
\xrightarrow{F_1(f)\left(\left(S(x)(\phi_{f{-1}})\right)^*\right)} F_1(f)\left(S(x)(A)^*\right).
\end{multline*}
Hence it suffices to show that the following two composites agree
\begin{multline*}
S(y)(B^*) = F_1(f)\left(S(x)\left(F_0(f^{-1})(B^*)\right)\right)
\\
\xrightarrow{F_1(f)\left(S(x)\left(T_0(f^{-1})(B)\right)\right)}
F_1(f)\left(S(x)\left(F_0(f^{-1})(B)^*\right)\right)
\\
\xrightarrow{F_1(f)\left(U(x)\left(F_0(f^{-1})(B)\right)\right)}
F_1(f)\left(S(x)\left(F_0(f^{-1})(B)\right)^*\right)
\end{multline*}
and
\begin{multline*}
S(y)(B^*) \xrightarrow{U(y)(B)} S(y)(B)^* = F_1(f)\left(F_1(f^{-1})\left(S(y)(B)^*\right)\right)
\\
\xrightarrow{F_1(f)\left(T_1(f^{-1})\left(S(y)(B)\right)\right)}
F_1(f)\left(F_1(f^{-1})\left(S(y)(B)\right)^*\right)
= F_1(f)\left(S(x)\left(F_0(f^{-1})(B)\right)^*\right)
\end{multline*}
(Notice that $\phi_{f^{-1}}$ has been eliminated.) This will follow by applying
$F_1(f)$ to the following diagram,
provided we can show that it does commute.
$$\xymatrix@!C=19em{
S(x)\left(F_0(f^{-1})(B^*)\right)= F_1(f^{-1})\left(S(y)(B^*)\right)
\ar[r]^-{S(x)\left(T_0(f^{-1})(B)\right)}
\ar[d]_-{F_1(f^{-1})\left(U(y)(B)\right)}
&
S(x)\left(F_0(f^{-1})(B)^*\right)
\ar[d]_-{U(x)\left(F_0(f^{-1})(B)\right)}
\\
F_1(f^{-1})\left(S(y)(B)^*\right)
\ar[r]^-{T_1(f^{-1})\left(S(y)(B)\right)}
&
S(x)\left(F_0(f^{-1})(B)\right)^*
}
$$
But the latter diagram commutes because we require the
following equality of functors of additive categories
with involution for $f^{-1} \colon x \to y$
(see~\eqref{(F_1(f),T_1(f))_circ_(S(y),U(y))_is_(S(x),U(x))_circ_(F_0(f),T_0(f))})
$$(F_0(f^{-1}),T_0(f^{-1})) \circ (S(x),U(x)) = (S(y),U(y)) \circ (F_1(f^{-1}),T_1(f^{-1})).$$
This finishes the proof that $\widehat{U}$ is a natural equivalence.
One easily checks that condition~\eqref{F(A)_F(A)aastast_F(Aastast)_F(Aast)ast} is satisfied by
$\widehat{U}$ since it holds for $U(x)$ for all objects $x$ in $\calg$.
This finishes the construction of the functor
of additive categories with involution $(S,U)$ (see~\eqref{intgf(calg)(S,U)}).

One easily checks
\begin{eqnarray}
\left(\intgf{\calg}{(S_2,U_2)}\right) \circ \left(\intgf{\calg}{(S_1,U_1)}\right)
& = & \intgf{\calg}{(S_2,U_2) \circ (S_1,U_1)}
\label{int_(S_2_U_2)_circ_int_(S_1,U_1)_is_int_(S_2,U_2)_circ_(S_1,U_1)}
\\
\intgf{\calg}{\id_F} & = & \id_{\intgf{\calg}{F}}.
\label{int_id_is_id_inv}
\end{eqnarray}

Given a functor of groupoids $W \colon \calg_1 \to \calg$  and a
functor $(F,T) \colon \calg \to \addcatinv$, the composition
with $W$ a yields a functor
$(F \circ W, T \circ W)$. Hence
both $\intgf{\calg_1}{(F,T)}  \circ W$ and $\intgf{\calg}{(F,T)}$
are additive categories with involutions. One easily checks that
$I_{\intgf{\calg}{F}} \circ W_* = W_* \circ I_{\intgf{\calg_1}{F \circ W}}$
holds for the functor $W_*$ defined  in~\eqref{W_ast}.
Hence
\begin{eqnarray}
&(W_*,\id) \colon \intgf{\calg_1}{(F,T)}  \circ W  \to \intgf{\calg}{(F,T)}. &
\label{(W_ast,id)}
\end{eqnarray}
is a functor of additive categories with involution.
One easily checks
\begin{eqnarray}
((W_2)_*,\id) \circ ((W_1)_*,\id) & = & ((W_2 \circ W_1)_*,\id);
\label{((W_2)_ast,id)_circ_((W_1)_ast,id)_is_((W_2_circ_W_1)_ast,id)}
\\
(\id_{\calg})_* & = & \id_{\intgf{\calg}{F}}.
\label{(id_calg)_ast_is_id_int_calg_F_inv}
\end{eqnarray}
These two constructions are compatible. Namely, we get
\begin{eqnarray}
\left(\intgf{\calg}{(S,U)} \right) \circ (W_*,\id) & = &
(W_*,\id) \circ \left(\intgf{\calg_1}{(S\circ W,U\circ W)} \right).
\label{compatibility_of_(W_ast,id)_and_int_(S,U)}
\end{eqnarray}

One easily checks
\begin{lemma}
\label{lem:(F_1,T_1)_ast_and_int_calg_S_and_equivalences}
\begin{enumerate}
\item \label{lem:(F_1,T_1)_ast_and_int_calg_S_and_equivalences:(F_1)_ast}
Let $W \colon \calg_1 \to \calg$ be an equivalence of connected groupoids.
Let $(F,T) \colon \calg \to \addcatinv$ be a contravariant functor.  Then
$$W_* \colon \intgf{\calg_1}{(F,T) \circ W} \to \intgf{\calg}{(F,T)}$$
is an equivalence of additive categories with involution.

\item \label{lem:(F_1,T_1)_ast_and_int_calg_S_and_equivalences:int_calg_S}
Let $\calg$ be a connected groupoid. Let $S \colon (F_1,T_1) \to (F_2,T_2)$ be a transformation
of contravariant functors $\calg \to \addcatinv$ such that for every
object $x$ in $\calg$ the functor $S(x) \colon F_0(x) \to F_1(x)$ is
an equivalence of additive categories. Then
$$\intgf{\calg}{S} \colon \intgf{\calg}{(F_1,T_1)} \to \intgf{\calg}{(F_2,T_2)}$$
is an equivalence of additive categories  with involution.
\end{enumerate}
\end{lemma}

%%%%%%%%%%%%%%%%%%%%%%%%%%%%%%%%%%%%%%%%%%%%%%%%%%%%%%%%%%%%%%%%%%%%%%%%%%%%%%%%%%%%%%%%%
%%   From crossed product rings with involution to additive categories with involution %%
%%%%%%%%%%%%%%%%%%%%%%%%%%%%%%%%%%%%%%%%%%%%%%%%%%%%%%%%%%%%%%%%%%%%%%%%%%%%%%%%%%%%%%%%%

\typeout{---  Section 8:  From crossed product rings with
involution to additive categories with involution
  ----}

\section{From crossed product rings with involution to additive categories with involution}
\label{sec:From_crossed_product_rings_with_involution_to_additive_categories_with_involution}

Next we want to extend Example~\ref{exa:R_ast_c,tau_as_add_cat_with_G-action} and
Lemma~\ref{lem:int_G_R_FGF(R)_c,tau_is_R_ast_c,tauG-FGF}
to rings and additive categories with involutions.
Let $R$ be a ring and let $G$ be a group. Suppose that we are given
maps of sets
\begin{eqnarray*}
c\colon G & \to & \aut(R), \quad g \mapsto c_g;
\\
\tau\colon G \times G & \to & R^{\times};
\\
w \colon G & \to & R,
\end{eqnarray*}
satisfying conditions~\eqref{tau_and_c_group_homo},
\eqref{tau_and_c_cocykel}, \eqref{c_e-Is_id},
\eqref{tau(g,e)_is_1}, \eqref{tau(e,g)_is_1}, \eqref{w(e)_is_1},
\eqref{w(gh)}, \eqref{overlinew(g)},  and~\eqref{overlinec_g(r)}.
We have constructed in Section~\ref{sec:Crossed_product_rings_and_involutions}
an involution on the crossed product $R\ast G = R \ast_{c,\tau} G$.
We have denoted this ring with involution by $R\ast G = R \ast_{c,\tau,w} G$
(see~\eqref{r_ast_c,tau,w_G}). The additive category
$\FGF{R\ast G}$ inherits the structure of an additive category with
involution (see Example~\ref{exa:R-mod_as_add_cat}).

We have introduced notion of an additive $G$-category with
involution in
Definition~\ref{def:additive_G-category_with_involution} and
constructed an explicit example $\FGP{R}_{c,\tau,w}$
in~\eqref{FGP(R)_(c,tau,w)}. All these constructions restrict to
the subcategory $\FGF{R} \subseteq \FGP{R}$ of finitely generated
free $R$-modules. Thus we obtain the additive $G$-category  with
involution
\begin{eqnarray}
 & \FGP{R}_{c,\tau,w} &  \label{FGF(R)_(c,tau,w)}
\end{eqnarray}

\begin{lemma} \label{lem:int_G_R_FGF(R)_c,tau_w_is_R_ast_c,tau,w_G-FGF}
Consider the data $(R,c,\tau,w)$ and the additive $G$-category with involution
$\FGF{R}_{c,\tau,w}$ of~\eqref{FGF(R)_(c,tau,w)}. Let $\intgf{G}{\FGF{R}_{c,\tau,w}}$
be the additive category with involution
defined in~\eqref{int_calg_F,E_int_calg_F}. Since $G$ regarded as a groupoid has
precisely one object, we can (and will) identify the set of objects in
$\intgf{G}{\FGF{R}_{c,\tau,w}}$ with the set of objects in $\FGF{R}_{c,\tau,w}$
which consists of pairs $(M,g)$ for $M$ a finitely generated free $R$-module and $g \in G$.
 Denote by
$\left(\intgf{G}{\FGF{R}_{c,\tau,w}}\right)_e$ the full subcategory of
$\intgf{G}{FGF{R}_{c,\tau,w}}$ consisting of objects of the shape
$(M,e)$ for $e \in G$ the unit element. Denote by $R \ast G = R
\ast_{c,\tau,w} G$ the ring with involution given by the crossed product ring
(see~\ref{r_ast_c,tau,w_G}). Then

\begin{enumerate}
\item
\label{lem:int_G_R_FGF(R)_c,tau_w_is_R_ast_c,tau_w_G-FGF:R_ast_G} There
is an equivalence of additive categories with involution
$$
(\alpha,\beta) \colon \left(\intgf{G}{\FGF{R}_{c,\tau,w}}\right)_e  \to  \FGF{R \ast_{c,\tau,w}G};
$$

\item  \label{lem:int_G_R_FGF(R)_c,tau_w_is_R_ast_c,tau_wG-FGF:e}
The inclusion
$$\left(\intgf{G}{\FGF{R}_{c,\tau,w}}\right)_e \to \intgf{G}{\FGF{R}_{c,\tau,w}}$$
is an equivalence of additive categories with involution.
\end{enumerate}

\end{lemma}
\begin{proof}\ref{lem:int_G_R_FGF(R)_c,tau_w_is_R_ast_c,tau_w_G-FGF:R_ast_G}
We have already constructed an equivalence of categories
$$
\alpha \colon \left(\intgf{G}{\FGF{R}_{c,\tau}}\right)_e  \to  \FGF{R \ast_{c,\tau}G};
$$
in
Lemma~\ref{lem:int_G_R_FGF(R)_c,tau_is_R_ast_c,tauG-FGF}~\ref{lem:int_G_R_FGF(R)_c,tau_is_R_ast_c,tauG-FGF:R_ast_G}.
We want to show that $\alpha$ is compatible with the involution, i.e.,
there is a functor of categories with involutions
$$
(\alpha,\beta) \colon \left(\intgf{G}{\FGF{R}_{c,\tau,w}}\right)_e  \to  \FGF{R \ast_{c,\tau,w}G}.
$$
The natural equivalence $\beta \colon \alpha \circ I_{\left(\intgf{G}{\FGF{R}_{c,\tau,w}}\right)_e}
\to
I_{\FGF{R \ast_{c,\tau,w}G}} \circ \alpha$ assigns to an object $(M,e)$ in
$\left(\intgf{G}{\FGF{R}_{c,\tau,w}}\right)_e$ the $R \ast_{c,\tau}G$-isomorphism
$$\beta(M,e) \colon R \ast_{c,\tau,w}G \otimes_R M^* \xrightarrow{\cong}
\left(R \ast_{c,\tau,w}G \otimes_R M\right)^*$$
given by $\beta(M,e)(u \otimes f)(v \otimes m) = vf(m)\overline{u}$
for $f \in M^*$, $u \in R \ast_{c,\tau}G$ and $m \in M$. Obviously
$\beta$ is compatible with the structures of additive categories.

Next we check that $\beta$ is a natural transformation. We have to show for a morphism
$\phi \colon (M,e) \to (N,e)$ in $\left(\intgf{G}{\FGF{R}_{c,\tau,w}}\right)_e$
that the following diagram commutes
$$\xymatrix{
R \ast G \otimes_R N^*
\ar[d]^-{\beta(N,e)} \ar[r]^-{\alpha(\phi^*)}
&
R \ast G \otimes_R M^*
\ar[d]^-{\beta(M,e)}
\\
\left(R \ast G \otimes_R N\right)^*
\ar[r]^-{\alpha(\phi)^*}
&
\left(R \ast G \otimes_R M\right)^*
}
$$
Recall that a morphism
$$\phi = \sum_{g \in G} g \cdot \phi_g \colon (M,e) \to (N,e)$$
in $\left(\intgf{G}{\FGF{R}_{c,\tau,w}}\right)_e$
is given by a collection of morphisms $\phi_g \colon (M,e) \to R_g(N,e) = (N,g)$
in $\FGF{R}_{c,\tau,w}$ for $g \in G$, where $\phi_g$ is a $R$-homomorphism $M \to \res_{c_g} N$.
We want to unravel what the dual morphism
$$\phi^* = \sum_{g \in G} g \cdot (\phi^*)_g \colon (N,e)^* = (N^*,e)\to (M,e)^* = (M^*,e)$$
in $\left(\intgf{G}{\FGF{R}_{c,\tau,w}}\right)_e$ is. It is given by a collection of morphisms
$\{(\phi^*)_g \colon (N^*,e) \to R_g(M^*,e) = (M^*,g)\mid g \in G\}$
in $\FGF{R}_{c,\tau,w}$, where $(\phi^*)_g$ is a $R$-homomorphism $N^* \to \res_{c_g} M^*$.
In $\FGF{R}_{c,\tau,w}$ the morphism $(\phi^*)_g$ is given by the composite
$$(N^*,e) = (N^*,g^{-1})\cdot g = R_g(N^*,g^{-1})\xrightarrow{R_g((\phi_{g^{-1}})^*)}
R_g(M^*,e) = (M^*,g).$$
The morphism $(\phi_{g^{-1}})^*$ is given by the composite
$$\res_{c_{g^{-1}}} N^* \xrightarrow{t_{g^{-1}}(N)} \left(\res_{c_{g^{-1}}} N\right)^*
\xrightarrow{\left(\phi_{g^{-1}}\right)^*} M^*.$$
Explicitly this is the map
$$N^* \to M^*, \quad f(x) \mapsto
c_{g^{-1}}^{-1}\left(f\circ \phi_{g^{-1}}(x)\right)\left(w(g^{-1})\tau(g,g^{-1})\right)^{-1}.$$
The morphism $R_g((\phi_{g^{-1}})^*)$ is the composite
$$
N^* = \res_{c_{g^{-1}g}} N^* \xrightarrow{L_{\tau(g^{-1},g)}} \res_{g} \res_{c_{g^{-1}}} N^*
\xrightarrow{\res_{c_g} (\phi_{g^{-1}})^*}  \res_{c_g} M^*.
$$
Hence the $R$-linear map $(\phi^*)_g \colon N^* \to \res_{c_g} M^*$ sends $f \in N^*$
to the element in $M^*$ given by
$$x~\mapsto~c_{g^{-1}}^{-1}\left(f\circ \phi_{g^{-1}}(x)\overline{\tau(g^{-1},g)}\right)
\left(w(g^{-1})\tau(g,g^{-1})\right)^{-1}.$$
This implies that the $R \ast G$-homomorphism
$$\alpha(\phi^*) \colon R \ast G \otimes_R N^* \to R \ast G \otimes_R M^*$$
sends $u \otimes f$ for $u \in R\ast G$ and
$f \in N^*$ to the $R$-linear map $M \to R$ given by
$$\sum_{g \in G} u \cdot \tau(g^{-1},g)^{-1} \cdot g^{-1} \otimes
\left(c_{g^{-1}}^{-1}\circ f\circ \phi_{g^{-1}}\right)
c_{g^{-1}}^{-1}\left(\overline{\tau(g^{-1},g)}\right)
\left(w(g^{-1})\tau(g,g^{-1})\right)^{-1}.$$
We conclude that the composite $\beta(M,e) \circ \alpha(\phi^*)$ sends
$u \otimes f$ for $u \in R\ast G$ and $f \in N^*$
to the $R$-linear map $R\ast G\otimes_R M \to R \ast G$
which maps $v \otimes x$ for $v \in R \ast G$ and $x \in M$ to the element in $R\ast G$
\begin{multline*}
\sum_{g \in G} v \cdot
\left(c_{g^{-1}}^{-1}\circ f\circ \phi_{g^{-1}}\right)(x) c_{g^{-1}}^{-1}
\left(\overline{\tau(g^{-1},g)}\right)
\left(w(g^{-1})\tau(g,g^{-1})\right)^{-1}
\\
\cdot \overline{u \cdot\tau(g^{-1},g)^{-1} \cdot g^{-1}}.
\end{multline*}
We compute that the composite $\alpha(\phi)^*\circ \beta(N,e)$ sends
$u \otimes f$ for $u \in R\ast G$ and $f \in N^*$
to the $R$-linear map $R\ast G\otimes_R M \to R \ast G$
which maps $v \otimes x$ for $v \in R \ast G$ and $x \in M$ to the element in $R\ast G$
\begin{eqnarray*}
\lefteqn{\beta(N,e)(u \otimes f)\left(\alpha(\phi)(v \otimes x)\right)}
& &
\\
& = &
\beta(N,e)(u \otimes f)\left(\sum_{g \in G}
v \cdot \tau(g^{-1},g)^{-1} \cdot g^{-1} \otimes \phi_g(x)\right)
\\
& = &
\sum_{g \in G} \beta(N,e)(u \otimes f)
\left(v \cdot \tau(g^{-1},g)^{-1} \cdot g^{-1} \otimes \phi_g(x)\right)
\\
& = &
\sum_{g \in G} v \cdot \tau(g^{-1},g)^{-1} \cdot g^{-1}
\cdot f\left(\phi_g(x)\right) \cdot\overline{u}.
\end{eqnarray*}
Hence it suffices to show for each $g \in G$, $u,v \in R \ast G$ and $x \in M$
\begin{multline*}
v \left(c_{g^{-1}}^{-1}\circ f\circ \phi_{g^{-1}}\right)(x) c_{g^{-1}}^{-1}
\left(\overline{\tau(g^{-1},g)}\right)
\left(w(g^{-1})\tau(g,g^{-1})\right)^{-1}
\\
\cdot \overline{u \cdot\tau(g^{-1},g)^{-1} \cdot g^{-1}}
\\
=
v \cdot \tau(g,g^{-1})^{-1} \cdot g \cdot f\left(\phi_{g^{-1}}(x)\right) \cdot\overline{u}.
\end{multline*}
Since
\begin{eqnarray*}
\overline{u \cdot\tau(g^{-1},g)^{-1} \cdot g^{-1}}
& = &
w(g^{-1}) c_{g}(\overline{\tau(g^{-1},g)^{-1}}) \cdot g \cdot \overline{u};
\\
g \cdot f\left(\phi_{g^{-1}}(x)\right) \cdot\overline{u}
& = &
c_g\left(f\left(\phi_{g^{-1}}(x)\right)\right) \cdot g \cdot\overline{u},
\end{eqnarray*}
it remains to show for all $g \in G$ and $x \in M$
\begin{multline*}
\left(c_{g^{-1}}^{-1}\circ f\circ \phi_{g^{-1}}\right)(x) c_{g^{-1}}^{-1}
\left(\overline{\tau(g^{-1},g)}\right)
\left(w(g^{-1})\tau(g,g^{-1})\right)^{-1}
\\
\cdot w(g^{-1}) c_{g}(\overline{\tau(g^{-1},g)^{-1}})
\\
=
\tau(g,g^{-1})^{-1} \cdot c_g\left(f\left(\phi_{g^{-1}}(x)\right)\right).
\end{multline*}
If we put $r = f \circ \phi_{g^{-1}}(x)$, this becomes equivalent
to showing for all $g \in G$ and $r \in R$
\begin{multline*}
c_{g^{-1}}^{-1}(r) c_{g^{-1}}^{-1}
\left(\overline{\tau(g^{-1},g)}\right)
\left(w(g^{-1})\tau(g,g^{-1})\right)^{-1}
\cdot w(g^{-1}) c_{g}(\overline{\tau(g^{-1},g)^{-1}})
\\
=
\tau(g,g^{-1})^{-1} \cdot c_g(r).
\end{multline*}
This is equivalent to showing
$$\tau(g,g^{-1}) c_{g^{-1}}^{-1}\left(r\overline{\tau(g^{-1},g)}\right)
\tau(g,g^{-1})^{-1}
=
c_{g}(r\overline{\tau(g^{-1},g)}).$$
From~\eqref{tau_and_c_group_homo} and~\eqref{c_e-Is_id} we conclude for $x \in R$
$$\tau(g,g^{-1})  c_{g^{-1}}^{-1}(x) \tau(g,g^{-1})^{-1}
=
c_g(x),$$
and the claim follows, i.e., $\beta$ is a natural equivalence.

It remains to check that the following diagram
(see~\eqref{F(A)_F(A)aastast_F(Aastast)_F(Aast)ast})
commutes for every object $(M,e)$ in $\FGF{R}_{c,\tau}[G]^e$.
$$
\xymatrix@!C=18em{R\ast G \otimes _R M \ar[r]^-{E_{\FGF{R \ast_{c,\tau,w} G}}(R\ast G \otimes _R M)}
\ar[d]^-{\alpha(E_{\FGF{R}_{c,\tau,w}}(M,e))}
& \left(R\ast G \otimes _R M\right)^{**} \ar[d]^-{\left(\beta(M,e)\right)^*}
\\
R\ast G \otimes _R M^{**}  \ar[r]^-{\beta((M,e)^*)}
&\left(R\ast G \otimes _R M^{*}\right)^*
}
$$
We consider an element $u \otimes x$ in the left upper corner for
$u \in R \ast G$ and $x \in M$. It is sent by the upper horizontal arrow
to the element in $\left(R\ast G \otimes _R M\right)^{**}$ which maps
$h \in \left(R\ast G \otimes_R M\right)^{*}$ to $\overline{h(u \otimes x)}$.
This element is mapped by the right vertical arrow to the element
in $\left(R\ast G \otimes _R M^{*}\right)^*$ which sends
$v \otimes f$ for $v \in R \ast G$ and $f \in M^*$ to
$$
\overline{\beta(M,e)(v \otimes f)(u \otimes x)} =
\overline{uf(x)\overline{v}} = v\overline{f(x)}\overline{u}.
$$
The left vertical arrow sends $u \otimes x$ to $u \otimes
I_{\FGF{R}}(x)$, where $I_{\FGF{R}}(x)$ sends $f \in M^*$ to $f(x)$.
This element is mapped by the lower horizontal arrow to the element in
$\left(R\ast G \otimes _R M^{*}\right)^*$ which sends $v \otimes f$
for $v \in R \ast G$ and $f \in M^*$ to
\begin{eqnarray*}
vI_{\FGF{R}}(x)(f)\overline{u}
& = & v\overline{f(x)}\overline{u}.
\end{eqnarray*}
This finishes the proof that
$$(\alpha,\beta) \colon \left(\intgf{G}{\FGF{R}_{c,\tau,w}}\right)_e  \to  \FGF{R \ast_{c,\tau,w}G}$$
is an equivalence of additive category with involutions.
\\[1mm]\ref{lem:int_G_R_FGF(R)_c,tau_w_is_R_ast_c,tau_wG-FGF:e}
This has already been proved in
Lemma~\ref{lem:int_G_R_FGF(R)_c,tau_is_R_ast_c,tauG-FGF}~\ref{lem:int_G_R_FGF(R)_c,tau_is_R_ast_c,tauG-FGF:e}.
\end{proof}

%%%%%%%%%%%%%%%%%%%%%%%%%%%%%%%%%%%%%%%%%%%%%%%%%%%%%%%%%%%%%%%%%%%%%%%%%%%%%%%%%%%%%%%%%
%%%%%%%%%%%%%%%%%%%%%%%%%%%% $G$-homology theories %%%%%%%%%%%%%%%%%%%%%%%%%%%%%%%%%%%%%%
%%%%%%%%%%%%%%%%%%%%%%%%%%%%%%%%%%%%%%%%%%%%%%%%%%%%%%%%%%%%%%%%%%%%%%%%%%%%%%%%%%%%%%%%%

\typeout{--------------------------  Section 9: $G$-homology theories ----------}

\section{$G$-homology theories}
\label{sec:G-homology_theories}

In this section we construct $G$-homology theories and discuss induction.

\begin{definition}[Transport groupoid]
\label{def:transport_groupoid}
Let $G$ be a group and let $\xi$ be a $G$-set. Define the \emph{transport groupoid}
$\calg^G(\xi)$ to be the following groupoid. The set of objects is $\xi$ itself.
For $x_1, x_2  \in \xi$ the set of morphisms from $x_1$ to $x_2$ consists of those elements $g$
in $G$ for which $gx_1 = x_2$ holds. Composition of morphisms
comes from the group multiplication in $G$.

A $G$-map $\alpha\colon \xi \to \eta$ of $G$-sets induces a covariant functor
$\calg^G(\alpha) \colon \calg^G(\xi) \to \calg^G(\eta)$ by sending an object $x \in \xi$
to the object $\alpha(x) \in \eta$. A morphism $g\colon x_1 \to x_2$ is sent to the morphism
$g \colon \alpha(x_1) \to \alpha(x_2)$.
\end{definition}

Fix a functor
\begin{eqnarray*}
& \bfE \colon \addcatinv \to \Spectra &
\end{eqnarray*}
which sends weak equivalences of additive categories with
involutions to weak homotopy
equivalences of spectra.

Let $G$ be a group.  Let $\Groupoidsover{G}$ be the category of
connected groupoids over $G$ considered as a groupoid with one object,
i.e., objects a covariant functors $F \colon \calg \to G$ with a
connected groupoid as source and $G$ as target and a morphism from
$F_0 \colon \calg_0 \to G$ to $F_1 \colon \calg_1 \to G$ is a
covariant functor $W \colon \calg_0 \to \calg_1$ satisfying $F_1 \circ
W = F_0$. For a $G$-set $S$ let
$$\pr_G \colon \calg^G(S) \to \calg(G/G) = G$$
be the functor induced by the projection $S \to G/G$.
The transport category yields a functor
$$\calg^G \colon \Or{G} \to \Groupoidsover{G}$$
by sending $G/H$ to
$\pr_G \colon \calg^G(G/H) \to \calg^G(G/G) = G$.

Let $\cala$ be an additive $G$-category with involution in the sense
of Definition~\ref{def:additive_G-category_with_involution}. We obtain
a functor
\begin{eqnarray}
& \bfE_{\cala} \colon \Or{G} \to \Spectra, \quad G/H \mapsto
\bfE\left(\intgf{\calg(G/H)}{\cala \circ \pr_G}\right). &
\label{bfE_calaG}
\end{eqnarray}
Associated to it there is a $G$-homology theory in the sense
of~\cite[Section~1]{Lueck(2002b)}
\begin{eqnarray}
&H^G_*(-;\bfE_{\cala})&
\label{HG_ast(-;bfEG_cala)}
\end{eqnarray} such that
$H^G_*(G/H;\bfE_{\cala}) \cong \pi_n\left(\bfE_{\cala}(G/H)\right)$
holds for every $n \in \IZ$ and every subgroup $H \subseteq G$.
Namely, define for a $G$-$CW$-complex $X$
$$H_n^G(X;\bfE_{\cala}) = \pi_n\left(\map_G(G/?,X)_+ \wedge_{\OrG{G}}
  \bfE_{\cala}(G/?)\right).$$
For more details about spectra and
spaces over a category and associated homology theories we refer
to~\cite{Davis-Lueck(1998)}. (Notice that there $\wedge_{\OrG{G}}$ is
denoted by $\otimes_{\OrG{G}}$.)

\begin{lemma}\label{lem:cala_to_calb}
  Let $f \colon \cala \to \calb$ be a weak equivalence of additive
  $G$-categories with involution.  Then the induced map
  $$H_n^G(X;\bfE_{f}) \colon H_n^G(X;\bfE_{\cala}) \xrightarrow{\cong}
  H_n^G(X;\bfE_{\cala})$$
  is a bijection for all $n \in \IZ$.
\end{lemma}
\begin{proof}
  This follows from
  Lemma~\ref{lem:(F_1,T_1)_ast_and_int_calg_S_and_equivalences}
  and~\cite[Lemma~4.6]{Davis-Lueck(1998)}.
\end{proof}

Let $\phi\colon K \to G$ be a group homomorphism. Given a
$K$-$CW$-complex $X$, let $G \times_{\phi} X$ be the $G$-$CW$-complex
obtained from $X$ by induction with $\phi$.  If $\calh^G_*(-)$ is a
$G$-homology theory, then $\calh^G(\phi_*(-))$ is a $K$-homology
theory.  The next result is essentially the same as the proof of the
existence of an induction structure
in~\cite[Lemma~6.1]{Bartels-Echterhoff-Lueck(2007colim)}.

\begin{lemma} \label{lem:homology_and_induction}
  Let $\phi\colon K \to G$ be a group homomorphism. Let $\cala$ be an
  additive $G$-category with involution in the sense of
  Definition~\ref{def:additive_G-category_with_involution}.  Let
  $\res_{\phi} \cala$ be the  additive $K$-category with
  involution obtained from $\cala$ by restriction with $\phi$.

  Then there is a transformation of $K$-homology theories
  $$\sigma_* \colon H^K_*(-;\bfE_{\res_{\phi} \cala}) \to H^G_*(\phi_*
  (-);\bfE_{\cala})$$
  If $X$ is a $K$-$CW$-complex on which
  $\ker(\phi)$ acts trivially, then
  $$\sigma_n \colon H^K_n(X;\bfE_{\res_{\phi} \cala})
  \xrightarrow{\cong} H^G_n(\phi_*X;\bfE_{\cala})$$
  is bijective
  for all $n \in \IZ$.
\end{lemma}
\begin{proof}
  We have to construct for every $K$-$CW$-complex $X$ a natural
  transformation
\begin{multline}
  \map_K(K/?,X)_+ \wedge_{\OrG{K}}
  \bfE\left(\intgf{\calg^K(K/?)}{\res_{\phi} \cala \circ \pr_K}\right)
  \\
  \to \map_G(G/?,\phi_*X)_+ \wedge_{\OrG{G}}
  \bfE\left(\intgf{\calg^G(G/?)}{\cala \circ \pr_G}\right).
\label{desired_map_of_spectra)}
\end{multline}
The group homomorphism $\phi$ induces for every transitive $K$-set
$\xi$ a functor, natural in $\xi$,
$$\calg^{\phi}(\xi) \colon \calg^K(\xi) \to \calg^G(\phi_*\xi)$$
which
sends an object $x \in \xi$ to the object $(e,x)$ in $G \times_{\phi}
\xi$ and sends a morphism given by $k \in K$ to the morphism given by
$\phi(k)$. We obtain for every transitive $K$-set $\xi$ a functor of
additive categories with involutions, natural in $\xi$
(see~\eqref{(W_ast,id)})
$$
\calg^{\phi}(\xi)_* \colon \intgf{\calg^K(\xi)}{\res_{\phi} \cala
  \circ \pr_K} = \intgf{\calg^K(\xi)}{\cala \circ \pr_G \circ
  \calg^{\phi}(\xi)} \to \intgf{\calg^G(\phi_* \xi)}{\cala \circ
  \pr_G}.
$$
Thus we obtain a map of spectra
\begin{multline*}
  \map_K(K/?,X)_+ \wedge_{\OrG{K}}
  \bfE\left(\intgf{\calg^K(K/?)}{\res_{\phi} \cala \circ \pr_K}\right)
  \\
  \to \map_K(K/?,X)_+ \wedge_{\OrG{K}}
  \bfE\left(\intgf{\calg^G(\phi_*(K/?))}{\cala \circ \pr_G}\right).
\end{multline*}
>From the adjunction of induction and restriction with the functor
$$\OrG{\phi} \colon \OrG{K} \to \OrG{G}, \quad K/H \mapsto \phi_*
K/H,$$
and the canonical map of contravariant $\Or(G)$-spaces
$$\OrG{\phi}_*\left(\map_K(K/?,X)\right) \to \map_G(G/?,\phi_*X),$$
which is an isomorphism for a $K$-$CW$-complexes $X$, we obtain
maps of spectra
\begin{multline*}
  \map_K(K/?,X)_+ \wedge_{\OrG{K}}
  \bfE\left(\intgf{\calg^G(\phi_*(K/?))}{\cala \circ \pr_G}\right)
  \\
  \cong \map_K(K/?,X)_+ \wedge_{\OrG{K}}
  \OrG{\phi}^*\left(\bfE\left(\intgf{\calg^G(G/?)}{\cala \circ
        \pr_G}\right)\right)
  \\
  \cong \OrG{\phi}_*\left(\map_K(K/?,X)\right)_+ \wedge_{\OrG{G}}
  \bfE\left(\intgf{\calg^G(G/?)}{\cala \circ \pr_G}\right)
  \\
  \cong \map_G(G/?,\phi_*X)_+ \wedge_{\OrG{G}}
  \bfE\left(\intgf{\calg^G(G/?)}{\cala \circ \pr_G}\right).
\end{multline*}
Now the desired map of spectra~\eqref{desired_map_of_spectra)} is the
composite of the two maps above.

The proof that $\tau_n(X)$ is bijective if $\ker(\phi)$ acts freely on
$X$ is the same as the one
of~\cite[Lemma~1.5]{Bartels-Echterhoff-Lueck(2007colim)}.
\end{proof}

%%%%%%%%%%%%%%%%%%%%%%%%%%%%%%%%%%%%%%%%%%%%%%%%%%%%%%%%%%%%%%%%%%%%%%%%%%%%%%%%%%%%%%%%%
%%%%%%%%%%%% Z-categories and additive categories with involutions %%%%%%%%%%%%%%%%%%%%%%
%%%%%%%%%%%%%%%%%%%%%%%%%%%%%%%%%%%%%%%%%%%%%%%%%%%%%%%%%%%%%%%%%%%%%%%%%%%%%%%%%%%%%%%%%

\typeout{------  Section 10: Z-categories and additive categories with involutions ------}

\section{$\IZ$-categories and additive categories with involutions}
\label{sec:Z-categories_and_additive_categories_with_involutionsG-homology_theories_and_restriction}

For technical reason it will be useful that $\cala$ comes with a
(strictly associative) functorial direct sum.  It will be used in
the definition of the category $\ind_{\phi}\cala $
in~\eqref{ind_phi_cala} and in functorial constructions about
categories arising in controlled topology. (See for 
instance~\cite[Section~2.2]{Bartels-Farrell-Jones-Reich(2004)},
\cite[Section~3]{Bartels-Lueck-Reich(2007hyper)}.)

\begin{definition}[$\IZ$-category (with involution)]
\label{def:Z-category_with_inv}
A \emph{$\IZ$-category} $\cala$ is  an additive
category except that we drop the condition that finite direct sums
do exists. More precisely, a  $\IZ$-category $\cala$ is a small
category such that for two objects $A$ and $B$ the morphism set
$\mor_{\cala}(A,B)$ has the structure of an abelian group and
composition yields bilinear maps $\mor_{\cala}(A,B) \times
\mor_{\cala}(B,C) \to \mor_{\cala}(B,C)$.

The notion of a \emph{$\IZ$-category
with involution} $\cala$ is defined analogously. Namely, we require the
existence of the pair $(I_{\cala},E_{\cala})$ with the same axioms
as in Section~\ref{sec:Additive_categories_with_involution} except
that we forget everything about finite direct sums.
\end{definition}

Of course an additive category (with involution) is a $\IZ$-category (with involution),
just forget the existence of the direct sum of two objects.

Given  a $\IZ$-category  $\cala$, we can enlarge it to an additive category
$\cala_{\oplus}$ with a functorial direct sums as follows.   The
objects in $\cala_{\oplus}$ are $n$-tuples $\underline{A} = (A_1,A_2,
\ldots, A_n)$ consisting of objects $A_i$ in $\cala$ for $i = 1,
2, \ldots, n$ and $n = 0,1,2, \ldots$, where we think of the empty set as $0$-tuple which we
denote by $0$.  the $\IZ$-module of morphisms from $\underline{A}
= (A_1, \ldots , A_m)$ to $\underline{B} = (B_1, \ldots , B_n)$ is
given by
$$\mor_{\cala_{\oplus}}(\underline{A},\underline{B}) := \bigoplus_{1
  \le i \le m, 1 \le j \le n} \; \mor_{\cala}(A_i,B_j).$$
Given a
morphism $f\colon \underline{A} \to \underline{B}$, we denote by
$f_{i,j} \colon A_i \to B_j$ the component which belongs to $i \in
\{1, \ldots, m\}$ and $j \in \{1, \ldots, n\}$. If $A$ or $B$ is the
empty tuple, then $\mor_{\cala_{\oplus}}(A,B)$ is defined to be the
trivial $\IZ$-module. The composition of $f \colon \underline{A}
\to \underline{B}$ and $\underline{g} \colon \underline{B} \to
\underline{C}$ for objects $\underline{A} = (A_1, \ldots, A_m)$,
$\underline{B} = (B_1, \ldots, B_n)$ and $\underline{C} = (C_1,
\ldots, C_p)$ is defined by
$$(g \circ f)_{i,k} := \sum_{j=1}^n g_{j,k} \circ f_{i,j}.$$
The sum
on $\cala_{\oplus}$ is defined on objects by sticking the tuples
together, i.e., for $\underline{A} = (A_1, \ldots, A_m)$ and
$\underline{B} = (B_1, \ldots, B_n)$ define
$$\underline{A} \oplus \underline{B} := (A_1, \ldots ,A_m,B_1, \ldots
,B_n).$$ The definition of the sum of two morphisms is now
obvious. The zero object is given by the empty tuple $0$.
The construction is strictly associative.  These
data define the structure of an additive category with functorial direct sum on
$\cala_{\oplus}$. Notice that this is more than an additive category since
for an additive category the existence of the direct sum of two objects is required
but not a functorial model.

In the sequel functorial direct sum is always to be understood to be
strictly associative, i.e., we have for three objects $A_1$, $A_2$ and
$A_3$ the equality $(A_1 \oplus A_2) \oplus A_3 = A_1 \oplus (A_2
\oplus A_3)$ and we will and can omit the brackets from now on in the
notion. We have constructed a functor from the category of
$\IZ$-categories to the category of additive categories with
functorial direct sum
$$\oplus \colon \Zcat \to \addcat_{\oplus}, \quad \cala \mapsto
\cala_{\oplus}.$$
Let
$$\forget \colon \addcat_{\oplus} \to \Zcat$$
be the forgetful
functor.

\begin{lemma}
\label{lem:adjoint_pair_(oplus,forget)}
\begin{enumerate}

\item \label{lem:adjoint_pair_(oplus,forget):adjoint}
We obtain an adjoint pair of functors $(\oplus,\forget)$.

\item \label{lem:adjoint_pair_(oplus,forget):equivalence}
We get for every $\IZ$-category $\cala$ a  functor
of $\IZ$-categories
$$Q_{\cala} \colon \cala \to \forget(\cala_{\oplus})$$
which is natural in $\cala$.

If $\cala$ is already an additive category, $Q_{\cala}$ is an equivalence
of additive categories.
\end{enumerate}
\end{lemma}
\begin{proof}\ref{lem:adjoint_pair_(oplus,forget):adjoint} We have to construct
  for every $\IZ$-category $\cala$ and every additive category $\calb$
  with functorial direct sum to one another inverse maps
  $$\alpha \colon \func_{\addcat_{\oplus}}(\cala_{\oplus},\calb) \to
  \func_{\Zcat}(\cala,\forget(\calb))$$
  and
  $$\beta \colon \func_{\Zcat}(\cala,\forget(\calb)) \to
  \func_{\addcat_{\oplus}}(\cala_{\oplus},\calb).
  $$
  Given $F \colon \cala_{\oplus} \to \calb$, define $\alpha(F)
  \colon \cala \to \calb$ to be the composite of $F$ with the obvious
  inclusion $Q_{\cala} \colon \cala \to \cala_{\oplus}$ which sends $A$ to
  $(A)$.  Given $F \colon \cala \to \forget(\calb)$, define $\beta(F)
  \colon \cala_{\oplus} \to \calb$ by sending $(A_1,A_2, \ldots ,
  A_n)$ to $F(A_1) \oplus F(A_2) \oplus \cdots \oplus F(A_n)$.
  \\[2mm]\ref{lem:adjoint_pair_(oplus,forget):equivalence} We have
  defined $Q_{\cala}$ already above. It is the adjoint of the identity on
  $\cala_{\oplus}$. Obviously $Q_{\cala}$ induces a bijection
  $\mor_{\cala}(A,B) \to \mor_{\cala_{\oplus}}(Q_{\cala}(A),Q_{\cala}(B))$
  for every objects $A,B \in  \cala$.
  Suppose that $\cala$ is an additive category. Then every object $(A_1,A_2, \ldots, A_n)$
  in $\cala_{\oplus}$ is isomorphic to an object in the image of $P_{\cala}$, namely
  to $P_{\cala}(A_1 \oplus A_2\oplus \cdots A_n) = (A_1 \oplus A_2\oplus \cdots \oplus A_n)$.
  Hence  $Q_{\cala}$ is an equivalence of additive categories.
\end{proof}

\begin{definition}[Additive category with functorial direct sum and involution]
\label{def:additive_category_with_oplus_and_inv}
An \emph{additive category with functorial sum and involution}
is an additive category with (strictly associative) functorial sum $\oplus$
and involution $(I,E)$ which are strictly compatible with one another, i.e.,
if $A_1$ and $A_2$ are two objects in $\cala$, then
$I(A_1 \oplus A_2) = I(A_1) \oplus I(A_2)$ and $E(A_1 \oplus A_2) = E(A_1) \oplus E(A_2)$
hold.
\end{definition}

One easily checks that if the $\IZ$-category $\cala$ comes with an
involution $(I_{\cala},E_{\cala})$, the additive category
$\cala_{\oplus}$ constructed above inherits the structure of an
additive category with functorial direct sum and involution in the
sense of Definition~\ref{def:additive_category_with_oplus_and_inv}.
Namely, define
\begin{eqnarray*}
I_{\cala_{\oplus}}\left((A_1,A_2, \ldots , A_n)\right)
& = &
\left(I_{\cala}(A_1), I_{\cala}(A_1), \ldots , I_{\cala}(A_1)\right);
\\
E_{\cala_{\oplus}}\left((A_1,A_2, \ldots , A_n)\right)
& = &
E_{\cala}(A_1) \oplus E_{\cala}(A_2) \oplus \cdots \oplus  E_{\cala}(A_1).
\end{eqnarray*}

We obtain a functor from the category of $\IZ$-categories with
involution to the category of additive categories with functorial
direct sum and involution
$$\oplus \colon \Zcatinv \to \addcatinv_{\oplus}, \quad \cala \mapsto
\cala_{\oplus}.$$
Let
$$\forget \colon \addcatinv_{\oplus} \to \Zcatinv$$
be the forgetful functor. One easily extends the proof of
Lemma~\ref{lem:adjoint_pair_(oplus,forget)} to the case with involution.

\begin{lemma}
\label{lem:adjoint_pair_(oplus,forget)_with_inv}
\begin{enumerate}

\item \label{lem:adjoint_pair_(oplus,forget)_with_inv:adjoint} We
  obtain an adjoint pair of functors $(\oplus,\forget)$.

\item \label{lem:adjoint_pair_(oplus,forget)_with_inv:equivalence} We get for
  every $\IZ$-category with involution $\cala$ a functor of
  $\IZ$-categories with involution
  $$Q_{\cala} \colon \cala \to \forget(\cala_{\oplus})$$
  which is
  natural in $\cala$.

  If $\cala$ is already an additive category with involution, then
  $Q_{\cala}$ is an equivalence of additive categories with
  involution.
\end{enumerate}
\end{lemma}

\begin{definition}\label{def:Z-G-category_with_involution}
  A \emph{$\IZ$-$G$-category with involution} $\cala$ is the
  same as an additive $G$-category in the sense of
  Definition~\ref{def:additive_G-category_with_involution} except that
  one forgets about the direct sum.
\end{definition}

\begin{definition}[Additive $G$-category with functorial sum and (strict) involution]
\label{def:additive_G-category_with_oplus_and_inv}
An \emph{additive $G$-category with functorial sum and involution}
is an additive $G$-category with (strictly associative) functorial sum $\oplus$
and involution $(I,E)$ which are strictly compatible with one another, i.e.,
we have:

\begin{enumerate}

\item If  $A_1$ and $A_2$ are two objects in $\cala$, then
$I(A_1 \oplus A_2) = I(A_1) \oplus I(A_2)$ and $E(A_1 \oplus A_2) = E(A_1) \oplus E(A_2)$
hold;

\item If  $A_1$ and $A_2$ are two objects in $\cala$ and $g \in G$, then
$R_g(A_1) \oplus R_g(A_2) = R_g(A_1 \oplus A_2)$ holds;

\item If  $A$ is an object in $\cala$, then
$I(R_g(A)) = R_g(I(A))$ and $E(R_g(A)) = R_g(E(A))$ hold.

\end{enumerate}

If the involution is strict in the sense of
Section~\ref{sec:Additive_categories_with_involution}, i.e., $E = \id$
and $I \circ I = \id$, we call $\cala$ an
\emph{additive $G$-category with functorial sum and strict involution}.

Define a $\IZ$-$G$-category with (strict) involution analogously,
just forget the direct sum.
\end{definition}

 We obtain a functor from the category of $\IZ$-$G$-categories with
involution to the category of additive categories with functorial
direct sum and involution
$$\oplus \colon \ZGcatinv \to \addGcatinv_{\oplus}, \quad \cala \mapsto
\cala_{\oplus}.$$
Let
$$\forget \colon \addGcatinv_{\oplus} \to \Zcatinv$$
be the forgetful functor. One easily extends the proof of
Lemma~\ref{lem:adjoint_pair_(oplus,forget)} to the case with
$G$-action and involution.

\begin{lemma}
\label{lem:adjoint_pair_(oplus,forget)_with_inv_and_G}
\begin{enumerate}

\item \label{lem:adjoint_pair_(oplus,forget)_with_inv_and_G:adjoint} We
  obtain an adjoint pair of functors $(\oplus,\forget)$.

\item \label{lem:adjoint_pair_(oplus,forget))_with_inv_and_G:equivalence} We get for
  every $\IZ$-$G$-category with involution $\cala$ a functor of
  $\IZ$-categories with involution
  $$Q_{\cala} \colon \cala \to \forget(\cala_{\oplus})$$
  which is
  natural in $\cala$.

  If $\cala$ is already an additive $G$-category with involution, then
  $Q_{\cala}$ is an equivalence of additive $G$-categories with
  involution;

\item \label{lem:adjoint_pair_(oplus,forget))_with_inv_and_G:strict}
The corresponding definitions and results carry over to the case of strict involutions.
\end{enumerate}
\end{lemma}

\begin{remark} \label{rem:comparison_int_ast}
Given an additive $G$-category $\cala$ and a $G$-set $T$, we have
constructed the additive $G$-category 
$\left(\intgf{\calg(T)}{\cala \circ \pr_G}\right)_{\oplus}$.
Let $\cala \ast_G T$ be the additive $G$-category defined 
in~\cite[Definition~2.1]{Bartels-Reich(2005)}. We obtain a functor of $\IZ$-categories
$$\rho(T) \colon \intgf{\calg(T)}{\cala \circ \pr_G} \to \cala \ast_G T$$
by sending an object $(x,A)$ to the object $\{B_t \mid t \in T\}$ for
which $B_x = A$ if $x = t$ and $B_x = 0$ if $x \not= t$. It induces
a functor of additive categories  with functorial direct sum
$$\rho(T)_{\oplus} \colon 
\left(\intgf{\calg(T)}{\cala \circ\pr_G}\right)_{\oplus}  
\to \left(\cala \ast_G T\right)_{\oplus}.$$
Recall that we have the functor of $\IZ$-categories
$$Q_{\cala \ast_G T} \colon \cala \ast_G T 
\to \left(\cala \ast_G  T\right)_{\oplus}.$$
One easily checks that both $\rho(T)_{\oplus}$ and 
$Q_{\cala \ast_G T}$ are equivalences of additive categories and
natural in $T$.

If $\cala$ is an additive $G$-category with strict involution,
then we obtain on the source and the target of $\rho(T)_{\oplus}$
and of $Q_{\cala \ast_G T}$  strict involutions
such that both  $\rho(T)_{\oplus}$ and $Q_{\cala \ast_G T}$ are
equivalences of additive categories with strict involution.

This implies that the $G$-homology theories constructed for $K$- and
$L$-theory here and in~\cite[Definition~2.1]{Bartels-Reich(2005)} are
naturally isomorphic and lead to isomorphic assembly maps.

\end{remark}

%%%%%%%%%%%%%%%%%%%%%%%%%%%%%%%%%%%%%%%%%%%%%%%%%%%%%%%%%%%%%%%%%%%%%%%%%%%%%%%%%%%%%%%%%
%%%%%%%%%%%% $G$-homology theories and restriction %%%%%%%%%%%%%%%%%%%%%%%%%%%%%%%%%%%%%%
%%%%%%%%%%%%%%%%%%%%%%%%%%%%%%%%%%%%%%%%%%%%%%%%%%%%%%%%%%%%%%%%%%%%%%%%%%%%%%%%%%%%%%%%%

\typeout{------------------  Section 11: $G$-homology theories and restriction ----------}

\section{$G$-homology theories and restriction}
\label{sec:G-homology_theories_and_restriction}

Fix a functor
\begin{eqnarray}
& \bfE \colon \addcatinv \to \Spectra &
\label{bfE_addcatinf_to_spectra}
\end{eqnarray}
which sends weak equivalences of additive categories with involutions to weak homotopy
equivalences of spectra. We call it \emph{compatible with direct sums}
if for any family of additive categories with involutions
$\{\cala_i\mid i \in I\}$   the map induced by the canonical inclusions
$\cala_i \to  \bigoplus_{i\in I} \cala_i$ for $i \in I$
$$\bigvee_{i \in I}  \bfE(\cala_i)  \to \bfE\left(\bigoplus_{i\in I} \cala_i\right)$$
is a weak homotopy equivalence of spectra.

\begin{example} \label{rem:K-and_L_are_compatible_with_direct_sums}
  The most important examples for $\bfE$ will be for us the functor
  which sends an additive category $\cala$ to its non-connective
  algebraic $K$-theory spectrum $\bfK_{\cala}$ in the sense of
  Pedersen-Weibel~\cite{Pedersen-Weibel(1985)}, and the functor which
  sends an additive category with involution $\cala$ to its algebraic
  $L^{-\infty}$-spectrum $\bfL^{-\infty}_{\cala}$ in the sense of
  Ranicki (see~\cite{Ranicki(1988)}, \cite{Ranicki(1992)}
  and~\cite{Ranicki(1992a)}).  Both functors send weak equivalences to
  weak homotopy equivalences and are compatible with direct sums. The
  latter follows from the fact that they are compatible with finite
  direct sums and compatible with directed colimits. This is proven
  for rings in~\cite[Lemma~5.2]{Bartels-Echterhoff-Lueck(2007colim)},
  the proof carries over to additive categories with involution.
\end{example}

Given a $G$-$CW$-complex $X$ and a group homomorphism $\phi \colon K
\to G$, let $\phi^* X$ be the $K$-$CW$-complex obtained from $X$ by
restriction with $\phi$.  Given a $K$-homology theory $\calh_*^K$, we
obtain a $G$-homology theory by sending a $G$-$CW$-complex $X$ to
$\calh^K_*(\phi^* X)$.  Recall that we have assigned to an additive
$G$-category $\cala$ with involution a $G$-homology theory
$H_*^G(-;\bfE_{\cala})$ in~\eqref{HG_ast(-;bfEG_cala)}.  The main
result of this section is
\begin{theorem}
\label{the:homology_theories-and_restriction}
Suppose that the functor $\bfE$ of~\eqref{bfE_addcatinf_to_spectra} is
compatible with direct sums.  Let $\phi \colon K \to G$ be a group
homomorphism.  Let $\cala$ be a $\IZ$-$K$-category with involution
in the sense of Definition~\ref{def:Z-G-category_with_involution}.
Let $\ind_{\phi} \cala$ be the $G$-$\IZ$-category with
involution defined in~\eqref{ind_phi_cala}.

Then there is a natural equivalence of $G$-homology theories
$$\tau_* \colon \calh^K\bigl(\phi^*(-);\bfE_{\cala_{\oplus}}\bigr) \xrightarrow{\cong}
\calh^G\bigl(-;\bfE_{(\ind_{\phi} \cala)_{\oplus}}\bigr).$$
\end{theorem}
Its proof needs some preparation.

Given a contravariant functor $F \colon \calg \to \addcatinv$ from a
groupoid into the category $\addcatinv$ of additive categories with
involution, we have defined an additive category with involution $\intgf{\calg}{F}$
in~\eqref{int_calg_F,T}, provided  that $\calg$ is connected.
We want to drop the assumption that $\calg$ is connected.
The connectedness of $\calg$ was only used in the construction of the direct sum of two
objects in $\intgf{\calg}{F}$. Hence everything goes through if we
refine us to the construction of $\IZ$-categories with involution.
Namely, if we drop the connectivity assumption on $\calg$,
all constructions and all the functoriality properties
explained in Section~\ref{sec:Connected_groupoids_and_additive_categories_with_involutions}
remain true if we work within the category $\Zcatinv$ instead of $\addcatinv$.

Let $G$ and $K$ be groups. Consider a (left) $K$-set $\xi$ and a
$K$-$G$-biset $\eta$.  Then $G$ acts from the right on the transport
groupoid $\calg^K(\eta)$. Namely, for an element $g \in G$ the map
$R_g \colon \eta \to \eta, \; x \mapsto xg$ is $K$-equivariant and
induces a functor $\calg^K(R_g) \colon \calg^K(\eta) \to
\calg^K(\eta)$.

Consider a $K$-$\IZ$-category with involution $\cala$. Let $\pr_K
\colon \calg^K(\eta) \to \calg^K(K/K) = K$ be the functor induced by
the projection $\eta \to K/K$.  Then $\cala \circ \pr_K$ is a
contravariant functor $\calg^K(\eta) \to \Zcatinv$.  We obtain a
$\IZ$-category with involution $\intgf{\calg^K(\eta)}{\cala \circ
  \pr_K}$ (compare~\eqref{int_calg_F,T}). Given $g \in G$, the functor
$\calg^G(R_g) \colon \calg^K(\eta) \to \calg^K(\eta)$ induces a
functor of $\IZ$-categories with involution
(compare~\eqref{(W_ast,id)})
$$\calg^G(R_g) \colon \intgf{\calg^K(\eta)}{\cala \circ \pr_K} =
\intgf{\calg^K(\eta)}{\cala \circ \pr_K \circ \calg^K(R_g)} \to
\intgf{\calg^K(\eta)}{\cala \circ \pr_K},$$
which strictly commutes with the involution. Thus
$\intgf{\calg^K(\eta)}{\cala \circ \pr_K}$ becomes
a $\IZ$-$G$-category with involution in the sense of
Definition~\ref{def:Z-G-category_with_involution}.
We conclude that $\left(\intgf{\calg^K(\eta)}{\cala \circ \pr_K}\right) \circ
\pr_G$ is a contravariant functor $\calg^G(\xi) \to \Zcatinv$.  We
obtain a $\IZ$-category with involution (compare~\eqref{int_calg_F,T})
$$\intgf{\calg^G(\xi)}{\left(\intgf{\calg^K(\eta)}{\cala \circ
      \pr_K}\right) \circ \pr_G}.$$
Consider $\eta \times \xi$ as a
left $G \times K$ set by $(g,k) \cdot (y,x) = (kyg^{-1},gx)$.
Then $\cala \circ \pr_{G \times K}$ is a contravariant functor
$\calg^{G \times K}(\eta \times \xi) \to \Zcatinv$.  We obtain a
$\IZ$-category with involution (compare~\eqref{int_calg_F,T})
$$\intgf{\calg^{G \times K}(\eta \times \xi)}{\cala \circ \pr_{G\times K}}.$$

\begin{lemma}
\label{lem:indent:int_G_times_K_and_int_G_int_K}
There is an isomorphism of $\IZ$-categories with involution
$$\omega \colon
\intgf{\calg^G(\xi)}{\left(\intgf{\calg^K(\eta)}{\cala \circ \pr_K}\right) \circ \pr_G}
\xrightarrow{\cong}
\intgf{\calg^{G \times K}(\eta \times \xi)}{\cala \circ \pr_{G\times K}}$$
which is natural in both $\xi$ and $\eta$.
\end{lemma}
\begin{proof}
An object in
$\intgf{\calg^G(\xi)}{\left(\intgf{\calg^K(\eta)}{\cala \circ \pr_K}\right) \circ \pr_G}$
is given by $(x,(y,A))$, where $x \in \xi$ is an object in $\calg^G(\xi)$
and $(y,A)$ is an object in $\intgf{\calg^K(\eta)}{\cala \circ \pr_K}$
which is given by an object $y \in \eta$ in $\calg^K(\eta)$ and an object $A$ in $\cala$.
The object $(x,(y,A))$ is sent under $\omega$ to the object $((y,x),A)$ given by the object
$(y,x)$ in $\calg^{G \times K}(\eta \times \xi)$ and the object $A \in \cala$.

A morphism $\phi$ in
$\intgf{\calg^G(\xi)}{\left(\intgf{\calg^K(\eta)}{\cala \circ \pr_K}\right) \circ \pr_G}$
from $(x_1,(y_1,A_1))$ to $(x_1,(y_2,A_2))$ is given by $g \cdot \psi$ for a morphism
$g \colon x_1 \to x_2$ in $\calg^G(\xi)$ and a morphism
$\psi \colon (y_1,A_1) \to \calg^K(R_g)_*(y_2,A_2)$.
The morphism $\psi$ itself is given by $k \cdot \nu$ for a morphism
$k \colon y_1 \to y_2g$ in
$\calg^K(\eta)$ and a morphism $\nu \colon A \to r_k(A)$ in $\cala$.
Define the image of $\phi$ under $\omega$
to be the morphism in $\intgf{\calg^{G \times K}(\eta \times \xi)}{\cala \circ \pr}$ given by
the morphism $(g^{-1},k) \colon (y_1,x_1) \to (y_2,x_2)$ in $\calg^{G \times K}(\eta \times \xi)$
and the morphism $\phi \colon A \to r_k(A)$. This makes sense since
$r_k(A)$ is the image of $A$ under the functor $\cala \circ \pr(g^{-1},k)$.

One easily checks that $\omega$ is an isomorphism of $\IZ$-categories with involution
and natural with respect to $\xi$ and $\eta$.
\end{proof}

Let $\phi \colon K \to G$ be a group homomorphism and $\xi$ be a
$G$-set.  Let $\phi^* \xi$ be the $K$-set obtained from the $G$-set
$\xi$ by restriction with $\phi$.  Consider a
$K$-$\IZ$-category with involution $\cala$ in the sense of
Definition~\ref{def:Z-G-category_with_involution}. Let
$\phi^* G$ be the $K$-$G$-biset for which multiplication with $(k,g)
\in K \times G$ sends $x \in G$ to $\phi(k)xg^{-1}$.  We have
explained above how $\intgf{\calg^K(\phi^*G)}{\cala}$ can be
considered as a $G$-$\IZ$-category with involution. We will denote it
by
\begin{eqnarray}
\ind_{\phi} \cala & := & \intgf{\calg^K(\phi^*G)}{\cala}.
\label{ind_phi_cala}
\end{eqnarray}

\begin{lemma}
\label{lem:int_phiast_xi_A_is_int_xi_ind_A}
For every $G$-set $\xi$ there is a natural equivalence of $\IZ$-categories with involutions
$$\tau  \colon
\intgf{\calg^G(\xi)}{\ind_{\phi} \cala} \xrightarrow{\simeq}
\intgf{\calg^K(\phi^*\xi)}{\cala \circ \pr_K}.$$
It is natural in $\xi$.
\end{lemma}
\begin{proof}
Because of Lemma~\ref{lem:indent:int_G_times_K_and_int_G_int_K} it suffices to construct a natural
equivalence
$$\tau \colon  \intgf{\calg^{G \times K}(\phi^* G \times \xi)}{\cala \circ \pr}
\xrightarrow{\simeq} \intgf{\calg^K(\phi^*\xi)}{\cala\circ \pr_K}.$$
Consider the functor
$$W \colon \calg^{G \times K}(\phi^* G \times \xi) \to \calg^K(\phi^*\xi)$$
sending an object $(x,y) \in G \times \xi$ in $\calg^{G \times K}(\phi^* G \times \xi)$
to the object $xy \in \xi$ in $\calg^K(\phi^*\xi)$ and a morphism
$(g,k) \colon (x_1,y_1) \to (x_2,y_2)$ to the morphism
$k \colon xy_1 \to xy_2$. Now define $\tau$ to be
$$W_* \colon
\intgf{\calg^{G \times K}(\phi^* G \times \xi)}{\cala \circ \pr}
= \intgf{\calg^{G \times K}(\phi^* G \times \xi)}{\cala \circ \pr_K \circ W }
\xrightarrow{\simeq} \intgf{\calg^K(\phi^*\xi)}{\cala}$$
(see~\eqref{(W_ast,id)}). Since $W$ is a weak equivalence of
groupoids, $\tau$ is a weak equivalence of additive categories with involution by
Lemma~\ref{lem:(F_1)_ast_and_int_calg_S_and_equivalences}~\ref{lem:(F_1)_ast_and_int_calg_S_and_equivalences:(F_1)_ast}.
One easily checks that this construction is natural in $\xi$.
\end{proof}
Now we can give the proof of Theorem~\ref{the:homology_theories-and_restriction}.
\begin{proof}
In the sequel we write $E_{\oplus} \colon \Zcatinv \to \Spectra$
for the composite of the functor $\bfE$
of~\eqref{bfE_addcatinf_to_spectra} and the functor $\Zcatinv \to
\addcatinv$ sending $\cala$ to $\cala_{\oplus}$. Given  a
$G$-$CW$-complex $X$, we have to define a weak equivalence of
spectra
\begin{multline*}
\map_K(K/?,\phi^*X)_+ \wedge_{\OrG{K}} \bfE_{\oplus}\left(\intgf{\calg^K(K/?)}{\cala \circ \pr_K} \right)
\\
\to
\map_G(G/?,X)_+ \wedge_{\OrG{G}}
\bfE_{\oplus}\left(\intgf{\calg^G(G/?)}{\ind_{\phi}\cala \circ \pr_G} \right).
\end{multline*}
The left  hand side can be rewritten as
\begin{eqnarray*}
\lefteqn{\map_K(K/?,\phi^*X)_+ \wedge_{\OrG{K}}
\bfE_{\oplus}\left(\intgf{\calg^K(K/?)}{\cala  \circ \pr_K}\right)}
& &
\\ & = &
\map_G(\phi_*(K/?),X)_+ \wedge_{\OrG{K}} \bfE_{\oplus}\left(\intgf{\calg^K(K/?)}{\cala \circ \pr_K}\right)
\\
& = &
\map_G(G/?,X)_+ \wedge_{\OrG{G}} \map_K(\phi_*(K/??),G/?)_+
\wedge_{\OrG{K}} \bfE_{\oplus}\left( \intgf{\calg^K(K/??)}{\cala \circ \pr_K}\right)
\\
& = &
\map_G(G/?,X)_+ \wedge_{\OrG{G}} \map_K(K/??,\phi^*G/?)_+
\wedge_{\OrG{K}} \bfE_{\oplus}\left(\intgf{\calg^K(K/??)}{\cala \circ \pr_K}\right).
\end{eqnarray*}
Because of Lemma~\ref{lem:int_phiast_xi_A_is_int_xi_ind_A}
the right hand side can be identified with
\begin{eqnarray*}
\lefteqn{\map_G(G/?,X)_+ \wedge_{\OrG{G}}
\bfE_{\oplus}\left(\intgf{\calg^G(G/?)}{\ind_{\phi}\cala \circ \pr_G}\right)}
& &
\\
& = &
\map_G(G/?,X)_+ \wedge_{\OrG{G}}  \bfE_{\oplus}\left(\intgf{\calg^K(\phi^*G/?)}{\cala \circ \pr_K}\right).
\end{eqnarray*}
Hence we need to construct for every $K$-set $\xi$ a weak homotopy equivalence, natural in $\xi$
$$
\rho(\xi) \colon \map_K(K/??,\xi)_+ \wedge_{\OrG{K}}
\bfE_{\oplus}\left(\intgf{\calg^K(K/??)}{\cala  \circ \pr_K}\right)
\to
\bfE_{\oplus}\left(\intgf{\calg^K(\xi)}{\cala \circ \pr_K} \right).
$$
The map $\rho(\xi)$ sends an element in the source given by $(\phi,z)$ for
a $K$-map $\phi \colon K/?? \to \xi$ and
$z \in \bfE_{\oplus}\left(\intgf{\calg^K(K/??)}{\cala  \circ \pr_K}\right)$
to $\bfE_{\oplus}\left(\calg^K(\phi)_*\right)(z)$, where
$$\calg^K(\phi)_* \colon \intgf{\calg^K(K/??)}{\cala  \circ \pr_K}
 = \intgf{\calg^K(K/??)}{\cala  \circ \pr_K \circ \calg^K(\phi)}
\to \intgf{\calg^K(\xi)}{\cala  \circ \pr_K}$$
has been defined in~\eqref{(W_ast,id)}.
Obviously it is natural in $\xi$ and is an isomorphism if $\xi$ is a $K$-orbit.
For a family of $K$-sets $\{\xi_i \mid i \in I\}$ there is a natural isomorphism of spectra
\begin{multline*}
\bigvee_{i \in I} \left(\map_K(K/??,\xi_i)_+ \wedge_{\OrG{K}}
\bfE_{\oplus}\left(\intgf{\calg^K(K/??)}{\cala \circ \pr_K}\right)\right)
\\
\xrightarrow{\cong}
\map_K\left(K/??,\coprod_{i \in I} \xi_i\right)_+ \wedge_{\OrG{K}}
\bfE_{\oplus}\left(\intgf{\calg^K(K/??)}{\cala \circ \pr_K}\right).
\end{multline*}
We have
\begin{eqnarray*}
\coprod_{i \in I} \calg^K(\xi_i) & \cong &\calg^K\left(\coprod_{i \in I} \xi_i\right);
\\
\coprod_{i \in I} \intgf{\calg^K(\xi_i)}{\cala \circ \pr_k}
& \cong &
\intgf{\coprod_{i \in I} \calg^K(\xi_i)}{\cala \circ \pr_k};
\\
\bigoplus_{i \in I} \left(\intgf{\calg^K(\xi_i)}{\cala \circ \pr_k}\right)_{\oplus}
& \cong &
\left(\coprod_{i \in I}  \intgf{\calg^K(\xi_i)}{\cala \circ \pr_k}\right)_{\oplus}.
\end{eqnarray*}
By assumption $\bfE$ is compatible with direct sums.  Hence we obtain
a weak equivalence
$$\bigvee_{i \in I} \bfE_{\oplus}\left(\intgf{\calg^K(\xi_i)}{\cala
    \circ \pr_K}\right) \xrightarrow{\simeq}
\bfE_{\oplus}\left(\intgf{\calg^K\left(\coprod_{i \in I}
      \xi_i\right)}{\cala \circ \pr_K}\right).$$
We conclude that
$\rho\left(\coprod_{i \in I} \xi_i\right)$ is a weak homotopy
equivalence if and only if $\bigvee_{i \in I} \rho(\xi_i)$ is a weak
homotopy equivalence. Since a $K$-set is the disjoint union of its
$K$-orbits and a wedge of weak homotopy equivalences of spectra is
again a weak homotopy equivalence, $\rho(\xi)$ is a weak homotopy
equivalence for every $K$-set $\xi$.  This finishes the proof of
Theorem~\ref{the:homology_theories-and_restriction}.
\end{proof}

%%%%%%%%%%%%%%%%%%%%%%%%%%%%%%%%%%%%%%%%%%%%%%%%%%%%%%%%%%%%%%%%%%%%%%%%%%%%%%%%%%%%%%%%%
%%%%%%%%%%%% Proof of the main theorems %%%%%%%%%%%%%%%%%%%%%%%%%%%%%%%%%%%%%%
%%%%%%%%%%%%%%%%%%%%%%%%%%%%%%%%%%%%%%%%%%%%%%%%%%%%%%%%%%%%%%%%%%%%%%%%%%%%%%%%%%%%%%%%%

\typeout{------------------  Section 12: Proof of the main theorems ----------}

\section{Proof of the main theorems}
\label{sec:Proof_of_the_main_theorems}

In this section we can finally give the proofs
of Theorem~\ref{the:FJC_for_crossed_products},
Theorem~\ref{the:fibered_versus_unfibered} and
Theorem~\ref{the:strict_coefficients}.

\begin{proof}[Proof of Theorem~\ref{the:FJC_for_crossed_products}]
This follows from
Lemma~\ref{lem:(F_1,T_1)_ast_and_int_calg_S_and_equivalences} and
Lemma~\ref{lem:int_G_R_FGF(R)_c,tau_w_is_R_ast_c,tau,w_G-FGF}.
\end{proof}

\begin{proof}[Proof of Theorem~\ref{the:fibered_versus_unfibered}]
Let $\phi \colon K \to G$ be a group
homomorphism and let $\calb$ be a
additive $K$-category with involution. We have to show that the following assembly map is
bijective
\begin{eqnarray*}
& \asmb^{K,\calb}_n \colon H_*^K(\EGF{K}{\phi^*\VCyc};\bfL_{\calb}) \to
H_n^K(\pt;\bfL_{\calb}) =  L_n\left(\intgf{K}{\calb}\right).
&
\end{eqnarray*}
Since $\phi^*\EGF{G}{\VCyc}$ is a model for $\EGF{K}{\phi^*\VCyc}$,
this follows from the commutative diagram
\begin{eqnarray*}
&\xymatrix@!C=19em{
H_n^K(\EGF{K}{\phi^*\VCyc};\bfL_{\calb})
\ar[d]^-{H^K_n(\id;\bfL_{Q_{\calb}})}_-{\cong}
\ar[r]^-{H_n^K(\pr;\bfL_{\calb}) }
&
H_n^K(\pt;\bfL_{\calb}) = L_n\left(\intgf{K}{\calb}\right)
\ar[d]^-{H_n^K(\id;\bfL_{Q_{\calb}})}_-{\cong}
\\
H_n^K(\EGF{K}{\phi^*\VCyc};\bfL_{\calb_{\oplus}})
\ar[d]^-{\tau_n^{\phi}(\EGF{G}{\VCyc})}_-{\cong}
\ar[r]^-{H_n^K(\pr;\bfL_{\calb_{\oplus}}) }
&
H_n^K(\pt;\bfL_{\calb_\oplus}) = L_n\left(\intgf{K}{\calb_{\oplus}}\right)
\ar[d]^-{\tau_n^{\phi}(\pt)}_-{\cong}
\\
H_n^G\bigl(\EGF{G}{\VCyc};\bfL_{(\ind_{\phi} \calb)_{\oplus}}\bigr)
\ar[r]^-{H_n^G\bigl(\pr;\bfL_{(\ind_{\phi} \calb)_{\oplus}}\bigr)}
&
H_n^G\bigl(\pt;\bfL_{(\ind_{\phi} \calb)_{\oplus}}\bigr) =
L_n\left(\intgf{G}{(\ind_{\phi} \calb)_{\oplus}}\right)
}
&
\end{eqnarray*}
where $\pr$ denotes the projection onto the one-point-space $\pt$ and
$Q_{\calb} \colon \calb \to \calb_{\oplus}$ is the natural equivalence
coming from Lemma~\ref{lem:adjoint_pair_(oplus,forget)_with_inv_and_G}
and the vertical arrows are isomorphisms because of
Lemma~\ref{lem:cala_to_calb}  and
Theorem~\ref{the:homology_theories-and_restriction}.
\end{proof}

\begin{proof}[Proof of Theorem~\ref{the:strict_coefficients}]

  Given an additive category $\cala$ with involution $\cala$, we can
  consider it as an additive category with $(\IZ/2,v)$-operation as
  explained in Example~\ref{exa:additive_categories_with_involution}.
  If we apply Lemma~\ref{lem:adjoint_pair_(cals,forget)}, we obtain an
  additive category with strict involution $\cals^{\IZ/2}(\cala)$
  together with a weak equivalence of additive categories with
  involutions
  $$P_{\cala} \colon \cala \to \cals^{\IZ/2}(\cala)$$
  If $\cala$ is an
  additive $G$-category with involution in the sense of
  Definition~\ref{def:additive_G-category_with_involution}, then
  $P_{\cala}$ is an equivalence of additive $G$-categories with
  involution.

  If we apply
  Lemma~\ref{lem:adjoint_pair_(oplus,forget)_with_inv_and_G}, we
  obtain an additive $G$-category with functorial direct sum and strict
  involution $\cals^{\IZ/2}(\cala)_{\oplus}$ in the sense of
  Definition~\ref{def:additive_G-category_with_oplus_and_inv} and an
  equivalence of additive $G$-categories with strict involution
  $$Q_{\cals^{\IZ/2}(\cala)} \colon \cals^{\IZ/2}(\cala) \to
  \cals^{\IZ/2}(\cala)_{\oplus}$$
  The composite
  $$f:= Q_{\cals^{\IZ/2}(\cala)} \circ P_{\cala} \colon \cala \to
  \cals^{\IZ/2}(\cala)_{\oplus}$$
  is a weak equivalence of additive
  $G$-category with involution. Now the claim follows from the
  following commutative diagram
\begin{eqnarray*}
&\xymatrix@!C=19em{
H_n^G(\EGF{G}{\VCyc};\bfL_{\cala})
\ar[d]^-{H_n^K(\id;\bfL_{f})}_-{\cong}
\ar[r]^-{H_n^K(\pr;\bfL_{\cala}) }
&
H_n^G(\pt;\bfL_{\cala}) = L_n\left(\intgf{G}{\cala}\right)
\ar[d]^-{H_n^K(\id;\bfL_{f})}_-{\cong}
\\
H_n^G\bigl(\EGF{G}{\VCyc};\bfL_{\cals^{\IZ/2}(\cala)_{\oplus}}\bigr)
\ar[r]^-{H_n^K\bigl(\pr;\bfL_{\cals^{\IZ/2}(\cala)_{\oplus}}\bigr)}
&
H_n^G\bigl(\pt;\bfL_{\cala}\bigr) = L_n\left(\intgf{G}{\cals^{\IZ/2}(\cala)_{\oplus}}\right)
}
&
\end{eqnarray*}
whose vertical arrows are isomorphisms by
Lemma~\ref{lem:cala_to_calb}.
\end{proof}

%%%%%%%%%%%%%%%%%%%%%%%%%%%%%%%%%%%%%%%%%%%%%%%%%%%%%%%%%%%%%%%%%%%%%%%%%%%%%%%%%%%%%%%%%
%%%%%%%%%%%%%%%%%%%%%%%%%%%%%%%%%%%%%% References  %%%%%%%%%%%%%%%%%%%%%%%%%%%%%%%%%%%%%%
%%%%%%%%%%%%%%%%%%%%%%%%%%%%%%%%%%%%%%%%%%%%%%%%%%%%%%%%%%%%%%%%%%%%%%%%%%%%%%%%%%%%%%%%%

\typeout{--------------------------------  References  ---------------------------------}

%\bibliographystyle{abbrv}
%\bibliography{dbpub,dbpre,dbbl-coeffinvextra}

\begin{thebibliography}{10}

\bibitem{Bartels-Echterhoff-Lueck(2007colim)}
A.~Bartels, S.~Echterhoff, and W.~L\"uck.
\newblock Inheritance of isomorphism conjectures under colimits.
\newblock Preprintreihe SFB 478 --- Geometrische Strukturen in der Mathematik,
  Heft 452, M\"unster, arXiv:math.KT/0702460, 2007.

\bibitem{Bartels-Farrell-Jones-Reich(2004)}
A.~Bartels, T.~Farrell, L.~Jones, and H.~Reich.
\newblock On the isomorphism conjecture in algebraic {$K$}-theory.
\newblock {\em Topology}, 43(1):157--213, 2004.

\bibitem{Bartels-Lueck-Reich(2007hyper)}
A.~Bartels, W.~L\"uck, and H.~Reich.
\newblock The ${K}$-theoretic {F}arrell-{J}ones {C}onjecture for hyperbolic
  groups.
\newblock Preprintreihe SFB 478 --- Geometrische Strukturen in der Mathematik,
  Heft 434, M\"unster, arXiv:math.GT/0609685, 2007.

\bibitem{Bartels-Reich(2005)}
A.~Bartels and H.~Reich.
\newblock {C}oefficients for the {F}arrell-{J}ones conjecture.
\newblock Preprintreihe SFB 478 --- Geometrische Strukturen in der Mathematik,
  Heft 402, M\"unster, arXiv:math.KT/0510602, to appear in Advances, 2005.

\bibitem{Davis-Lueck(1998)}
J.~F. Davis and W.~L{\"u}ck.
\newblock Spaces over a category and assembly maps in isomorphism conjectures
  in ${K}$- and ${L}$-theory.
\newblock {\em $K$-Theory}, 15(3):201--252, 1998.

\bibitem{Farrell-Jones(1993a)}
F.~T. Farrell and L.~E. Jones.
\newblock Isomorphism conjectures in algebraic ${K}$-theory.
\newblock {\em J. Amer. Math. Soc.}, 6(2):249--297, 1993.

\bibitem{Farrell-Linnell(2003a)}
F.~T. Farrell and P.~A. Linnell.
\newblock {$K$}-theory of solvable groups.
\newblock {\em Proc. London Math. Soc. (3)}, 87(2):309--336, 2003.

\bibitem{Hambleton-Pedersen-Rosenthal(2007)}
I.~Hambleton, E.~K. Pedersen, and D.~Rosenthal.
\newblock Assembly maps for group extensions in {$K$}- and {$L$}-theory.
\newblock Preprint, arXiv: math.KT/0709.0437v1, 2007.

\bibitem{Lueck(2002b)}
W.~L{\"u}ck.
\newblock Chern characters for proper equivariant homology theories and
  applications to ${K}$- and ${L}$-theory.
\newblock {\em J. Reine Angew. Math.}, 543:193--234, 2002.

\bibitem{Lueck-Reich(2005)}
W.~L{\"u}ck and H.~Reich.
\newblock The {B}aum-{C}onnes and the {F}arrell-{J}ones conjectures in {$K$}-
  and {$L$}-theory.
\newblock In {\em Handbook of $K$-theory. Vol. 1, 2}, pages 703--842. Springer,
  Berlin, 2005.

\bibitem{Pedersen-Weibel(1985)}
E.~K. Pedersen and C.~A. Weibel.
\newblock A non-connective delooping of algebraic ${K}$-theory.
\newblock In {\em Algebraic and Geometric Topology; proc. conf. Rutgers Uni.,
  New Brunswick 1983}, volume 1126 of {\em Lecture notes in mathematics}, pages
  166--181. Springer, 1985.

\bibitem{Ranicki(1988)}
A.~A. Ranicki.
\newblock Additive ${L}$-theory.
\newblock {\em $K$-Theory}, 3(2):163--195, 1989.

\bibitem{Ranicki(1992)}
A.~A. Ranicki.
\newblock {\em Algebraic ${L}$-theory and topological manifolds}.
\newblock Cambridge University Press, Cambridge, 1992.

\bibitem{Ranicki(1992a)}
A.~A. Ranicki.
\newblock {\em Lower ${K}$- and ${L}$-theory}.
\newblock Cambridge University Press, Cambridge, 1992.

\bibitem{Thomason(1979)}
R.~W. Thomason.
\newblock Homotopy colimits in the category of small categories.
\newblock {\em Math. Proc. Cambridge Philos. Soc.}, 85(1):91--109, 1979.

\end{thebibliography}

\end{document}